\documentclass[final,letterpaper]{siamart250211}
\usepackage{stmaryrd}
\SetSymbolFont{stmry}{bold}{U}{stmry}{m}{n}
\usepackage{amsmath,bbm,mathbbol,mathtools}
\usepackage{tabularx,lmodern}
\usepackage{color,soul,enumitem}
\usepackage{epsfig,epstopdf}
\usepackage{algorithm,algpseudocode}
\usepackage{flushend}
\usepackage{verbatim}
\usepackage{caption}
\usepackage{graphicx,subfig}
\usepackage{arydshln}
\usepackage[level]{datetime}
\usepackage{float,multirow}
\usepackage{booktabs}
\usepackage{relsize}
\usepackage{setspace}
\usepackage{epstopdf}
\usepackage{appendix}
\usepackage{bm}
\usepackage{adjustbox}
\usepackage{array}
\usepackage{empheq}
\usepackage{hhline}
\usepackage{makecell}
\usepackage{hyperref}
\numberwithin{equation}{section}

\newcommand\twoitems[2]{%
\item#1%
\hspace{15pt}%
\hspace{\labelsep}#2
}

\title{Computing dynamics of spin-1 Bose-Einstein condensates with arbitrary-angle rotation via shear decomposition
\thanks{This work was partially supported by the National Key Research and Development Program of China No. 2024YFA1012803, the basic research fund of Tianjin University No. 2025XJ21-0010
and the National Natural Science Foundation of China No. 12271400 (L. Wei and Y. Zhang),
the National Natural Science Foundation of China No. 12561076 (X. Zhang),
and the National Natural Science Foundation of China No. 12601848 (S. Zhang).}}
\author{Lianbo Wei\thanks{National Key Laboratory of Information Systems Engineering, National University of Defense Technology, Changsha, 410073, China. ({\tt wlb@tju.edu.cn}).}\and
Shaobo Zhang\thanks{Corresponding author. School of Mathematics, Sichuan University, No.24 South Section 1, Yihuan Road,
Chengdu, 610065, China.
 ({\tt shaobo\_zhang@scu.edu.cn}).} \and
Xiaohua Zhang\thanks{Yunnan Key Laboratory of Modern Analytical Mathematics and Applications, Yunnan Normal University, Kunming, 650500, China.
 ({\tt zhangxiaohua@ynnu.edu.cn}).} \and
Yong Zhang \thanks{ Center for Applied Mathematics and KL-AAGDM, Tianjin University, Tianjin, 300072, China; State Key Laboratory of Synthetic Biology, Tianjin University, Tianjin, 300072, China. ({\tt zhang\_yong@tju.edu.cn}).}
}

\begin{document}

\newcommand{\be}{\begin{equation}}
\newcommand{\ee}{\end{equation}}
\newcommand{\ba}{\begin{array}}
\newcommand{\ea}{\end{array}}
\newcommand{\bea}{\begin{eqnarray}}
\newcommand{\eea}{\end{eqnarray}}
\newcommand{\beas}{\begin{eqnarray*}}
\newcommand{\eeas}{\end{eqnarray*}}

\newcommand{\fl}[2]{\frac{#1}{#2}}

\newcommand{\bx}{{\textbf x}}
\newcommand{\by}{{\textbf y}}
\newcommand{\bz}{{\textbf z}}
\newcommand{\bn}{{\textbf n}}
\newcommand{\br}{{\textbf r}}
\newcommand{\im}{i}

\newcommand{\bu}{{\textbf u}}
\newcommand{\bv}{{\textbf v}}
\newcommand{\bbf}{{\textbf f}}
\newcommand{\bbg}{{\textbf g}}
\newcommand{\bh}{{\textbf h}}

\newcommand{\Th}{{\mathcal T}_{\bf h}}

\newcommand{\dif}{{\textrm d}}

\newcommand{\gm}{{\gamma}}

\newcommand\bl{{\textbf B}_\textbf{L}}
\newcommand{\bk}{{\textbf k}}
\newcommand\bR{\mathbb{R}}
\newcommand{\hh}{\widetilde{h}}
\renewcommand{\theequation}{\arabic{section}.\arabic{equation}} %
\newtheorem{exmp}{Example}
\numberwithin{exmp}{section}   
\newtheorem{remark}{Remark}[section]
\newtheorem{prop}{Proposition}
\newtheorem{thm}{Theorem}
\newtheorem{lem}{Lemma}
\newtheorem{prf}{proof}

\DeclareRobustCommand{\rchi}{{\mathpalette\irchi\relax}}
\newcommand{\irchi}[2]{\raisebox{\depth}{$#1\chi$}} 

\maketitle
\begin{abstract}
We propose high-order splitting Fourier spectral methods to compute the dynamics of spin-1
Bose-Einstein condensates  with arbitrary-angle rotation via shear decomposition.
The Hamiltonian is split into a linear part,
which consists of Laplacian and arbitrary-angle rotation terms,
and a nonlinear part that includes the remaining terms.
For the linear subproblem,
we introduce an arbitrary-angle rotating coordinate transformation to eliminate the rotation term.
The transformed subproblem contains only the Laplacian term and can be solved exactly in Fourier space.
To numerically access this coordinate transformation accurately and efficiently,
we propose a spectrally accurate and fully explicit method
based on shear decomposition
and it holds true for arbitrary time step.
The nonlinear subproblem can be solved exactly in physical space.
Then, we can construct efficient and unconditionally stable numerical schemes by
composing the above  subproblem integrators via operator splitting,
achieving spectral accuracy in space and high-order in time.
Extensive numerical results are presented to confirm the accuracy and efficiency, verify the dynamics properties and investigate vortex-state evolution.
\end{abstract}


\pagestyle{myheadings}\thispagestyle{plain}
\markboth{L. Wei, S. Zhang, X. Zhang and Y. Zhang}
{dynamics of spin-1 BECs with arbitrary-angle rotation}

\section{Introduction}
The first experimental realization of Bose-Einstein condensates (BECs) in 1995 \cite{BECstart, NaBEC} opened a new era for exploring the quantum world from a macroscopic perspective,
leading to numerous significant advances \cite{BaoSpinRev,phase_rotatingBEC,  TiltModeVortices} that span atomic and molecular physics as well as quantum optics.
In the early stages, atoms were confined in magnetic traps, which froze their spin degrees of freedom.
With the application of advanced optical trapping techniques to confine ultracold atoms, the spin degrees of freedom became unfrozen, giving rise to spinor BECs, e.g., spin-1 BEC experimentally realized in $^{23}$Na \cite{NaSpin}.
Compared with scalar BECs, spinor BECs exhibit a richer variety of symmetry-breaking mechanisms and phases \cite{ phaseSpin,SpontSymmBreakSpinorBEC_Nature},
and reveal intriguing spin structures during the spin-mixing dynamics \cite{physical_spin}.

Quantized vortices \cite{ObseVortex, phase_rotatingBEC}, marking the superfluid nature of BECs, have made rotating BECs a subject of widespread interest and intensive study.
While early experiments focused on rotation along the $z$-axis \cite{early_rot_BEC}, real-world applications may demand rotation at arbitrary angles.
A recent study \cite{ArbiRotBEC} investigated a BEC confined by a rotating harmonic trap where the rotation axis is misaligned with all principal axes.
Further, Ref. \cite{ArbRotPCGGS} studied how the bent vortex line patterns evolve with the angle between the
rotation axis and the elongated direction (a principal axis) of the trap.
Thus, spin-1 BECs with arbitrary-angle rotation (ArbRot) represents a fascinating subject in the field of quantum
physics, promising the emergence of further novel quantum phenomena.

\begin{figure*}[htbp]
\centering
\includegraphics[height=4cm,width=4cm]{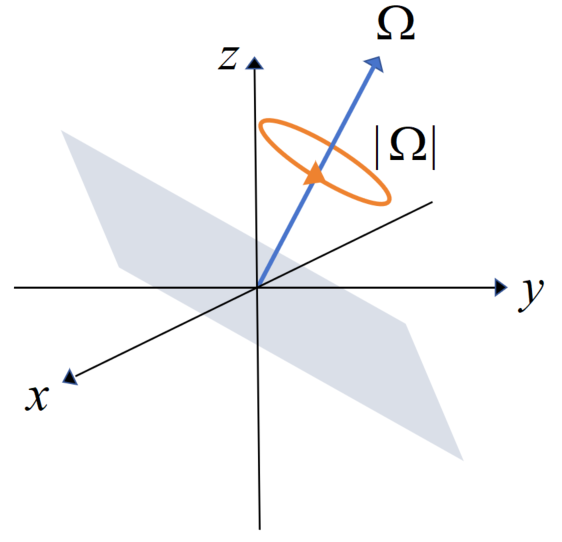}
\caption{Schematic diagram of the rotation with vector $\bm{\Omega}$.}
\label{fig_rotation}
\end{figure*}

At temperatures $T$ lower than the critical temperature $T_c$,
the macroscopic behavior of the spin-1 BECs with ArbRot can be well described by the complex-valued three-component wave function $\Psi = (\psi_1, \psi_0, \psi_{-1})^{\top}$,
whose evolution is governed by the dimensionless coupled Gross-Pitaevskii equations (CGPEs)
\cite{spin_1_Genrot_properties, spin1_BEC_model, spin_1BEC_Genrot, ArbRotPCGGS}
\begin{equation}\label{GGPE}
i\partial_t \Psi(\bx,t) = \Big[\Big(-\frac{1}{2}\Delta -\bm{\Omega} \cdot  \bm{L}  +V(\bx)+c_0 \rho\Big)  I_{3}+c_1 \mathbf{F}\cdot \mathbf{f} \Big]\Psi(\bx,t),
\end{equation}
where $t$ is the time variable, $\Delta=\nabla^2$ is the Laplacian operator with respect
to the spatial variable $\bx = (x, y, z)^{\top} \in \mathbb{R}^3$ and  $I_{3}$ is the $3\times 3 $ identity matrix.
$V(\bx)$ is a given
real-valued external trapping potential, which is determined by the type of system under study.
In most BECs experiments, a harmonic potential is employed to trap the condensate, i.e.,
\bea \label{harmonicPot}
V(\bx) = \frac{1}{2} (\gamma_x^2 x^2+ \gamma_y^2 y^2 +\gamma_z^2 z^2),
\eea
where $\gamma_{\alpha}$ $(\alpha=x, y, z)$ is the trapping frequency in the $\alpha$-direction.
The rotation vector
$
\bm{\Omega} = (\omega_x, \omega_y, \omega_z)^{\top} \in \mathbb{R}^3
$
represents the angular velocity vector,
and the angular momentum operator is given by
$
\bm{L} = \bx \wedge (-\im \nabla) = (L_x, L_y, L_z)^{\top}
$
with
\beas
L_x = -\im (y \partial_z - z\partial_y), \qquad L_y = -\im (z \partial_x - x \partial_z),
\qquad L_z = -\im (x \partial_y - y\partial_x).
\eeas
For readers' convenience, we present the graphical illustration of the rotation with vector $\bm{\Omega}$ in Figure \ref{fig_rotation}.

The dimensionless constants $c_0$ and $c_1$ describe the strengths of mean-field and  spin-exchange interactions respectively.
Here,
$\rho = \sum^{1}_{v=-1} \rho_v$ is the total density with $\rho_v= |\psi_v|^2$ being the $v$-th component density.
The spin-1 matrices $\mathbf{f}:=(f_x,f_y,f_z)^{\top}$ are given as
\begin{equation*}
f_x=\frac{1}{\sqrt{2}}\begin{pmatrix}
0&1&0\\
1&0&1\\
0&1&0
\end{pmatrix},\qquad f_y=\frac{i}{\sqrt{2}}\begin{pmatrix}
0&-1&0\\
1&0&-1\\
0&1&0\end{pmatrix},\qquad f_z=\begin{pmatrix}
1&0&0\\
0&0&0\\
0&0&-1\end{pmatrix}.
\end{equation*}
The corresponding spin vector $\mathbf{F}:=(F_x,F_y,F_z)^{\top}$, defined by
$F_v=\Psi^{\mathsf H} f_v \Psi$, can be written explicitly as
\beas	
F_x &=& \frac{1}{\sqrt{2}} [\bar{\psi}_1\psi_0+\bar{\psi}_0(\psi_1+\psi_{-1})+\bar{\psi}_{-1}\psi_0 ],\\	
F_y &=& \frac{i}{\sqrt{2}} [-\bar{\psi}_1\psi_0+\bar{\psi}_0(\psi_1-\psi_{-1})+\bar{\psi}_{-1}\psi_0],\\[0.15em]
F_z &=& |\psi_1|^2-|\psi_{-1}|^2.
\eeas
Then we can obtain
\begin{small}
\beas
\mathbf{F}\cdot \mathbf{f}:=F_xf_x + F_yf_y+F_zf_z =
\begin{pmatrix}
F_z                    & \frac{1}{\sqrt{2}}F_{-} & 0    \\
\frac{1}{\sqrt{2}}F_{+}& 0                       &\frac{1}{\sqrt{2}}F_{-} \\
0                      &\frac{1}{\sqrt{2}}F_{+}  &-F_z
\end{pmatrix},\eeas
\end{small}
where $F_{+} = F_x + iF_y = \bar{F_{-}} =\sqrt{2}\,(\bar{\psi}_1\psi_0+\bar{\psi}_0\psi_{-1})$ and  $\bar{\psi_{v}}$ is the conjugate of $\psi_{v}$.

As for  scalar BECs, there have been extensive theoretical and numerical studies on the dynamics,  and we refer to
\cite{Antoine2013Dynamic, BaoBECRev,BaoRotBECPolarCylind, Bao2009dynamicLaguerrFourSpin1, RotLagBao,RotADIBao,  Besse2017ExponIntegRotBEC, DyRot} for more details.
Among them, the time-splitting Fourier spectral method \cite{Antoine2013Dynamic} is one of the most successful methods:
it  is stable, highly efficient, and achieves spectral accuracy in space.
For the rotating case where the rotation axis is aligned with the $z$-axis, i.e., $\omega_x=\omega_y=0$,
several successful methods have been developed.
One efficient and spectrally accurate approach proposed in \cite{RotADIBao} relies on a proper application of the alternating direction implicit technique (ADI) to handle the coupling of variable coefficients and differential operators in the rotation term.
However, extending this method to construct a high-order scheme is somewhat tedious and complicated.
The exact splitting method (ESM) \cite{Bernier2021ExtSplit, Bernier2021ExtSplitAppl} takes a different strategy, decomposing  the linear part into a product of several simple sub-operators
that are exactly solvable in physical or Fourier space.
Alternatively, the rotating Lagrangian coordinate method (RLC) \cite{RotLagBao} eliminates the rotation term by introducing
time-dependent rotating Lagrangian coordinates \cite{Besse2017ExponIntegRotBEC},
which readily facilitates the implementation of high-order time marching schemes.
When real-time dynamics are required, one has to recover the wave function in Cartesian coordinates at each time step,
and the Fourier spectral interpolation yields spectral accuracy but at the price of large complexity of $O(N^{d+2})$ with naive summation implementation.
Actually, as pointed out in \cite{Threeshear},
such a coordinate transformation for classical rotation around the $z$-axis can be accelerated via one-dimensional Fast Fourier Transform (FFT) based on shear decomposition, and we choose to follow this line for the arbitrary rotation case.

Several methods originally developed for scalar BECs have been successfully extended to spinor BECs \cite{BaoSpinRev,DynLawspin1,Threeshear, CompactSplitRotSpin1}.
Symes \textit{et al.} derived an exact evaluation of the nonlinear subproblem
in physical space for spin-1 \cite{ExtNonSpin1} and spin-2 \cite{ExtNonSpin2} BECs.
Liu \textit{et al.} \cite{CompactSplitRotSpin1} solved the spin-1 rotating BECs 
using the ESM and derived its explicit splitting coefficients.
However, the direct extension of ESM to spin-1 BECs with ArbRot is non-trivial and somewhat challenging numerically.
To the best of our knowledge, numerical studies of the dynamics of spin-1 BECs with ArbRot
remain  relatively limited.
The primary difficulty lies in the rotation term, which consists of variable coefficients multiplying differential operators and cannot be diagonalized in either Fourier or physical space. This makes it challenging to simultaneously achieve high accuracy and computational efficiency.
Such difficulty is particularly pronounced for ArbRot, which is more complex than typical rotation that is aligned along the $z$-axis.

Inspired by the RLC method, this paper develops efficient high-order
Fourier spectral methods to compute the dynamics of spin-1 BECs with ArbRot. We
first apply an operator splitting technique to decouple the problem into a Laplacian-ArbRot subproblem, i.e., $-\frac{1}{2}\Delta-\bm{\Omega}\cdot \bm{L}$, and a nonlinear subproblem (all the other terms). For the Laplacian-ArbRot
subproblem, we introduce a rotating coordinate transformation to eliminate the ArbRot term.
The transformed subproblem can be solved exactly in Fourier space.
Based on the shear decomposition, the coordinate transformation can be implemented
by only one-dimensional FFT/iFFT.
Furthermore, to address the time-step constraint arising from the splitting coefficient singularity of shear decomposition,
we employ a time step decomposition technique to guarantee that it retains both accuracy and efficiency for arbitrary time step.
The nonlinear part is solved exactly in physical space. This facilitates the construction of
high-order operator splitting schemes with high computational efficiency.

The rest of the paper is organized as follows.
Some dynamical properties are shown in Section \ref{sec:dyn}.
In Section \ref{sec:NumeMeth}, we present high-order splitting Fourier spectral methods for simulating the dynamics, including a detailed description of the fast rotating-coordinate transformation algorithm and stability analysis.
In Section \ref{sec:NumeResu}, we provide  spatial/temporal accuracy confirmation
and efficiency tests, together with extensive numerical results for confirming
dynamical laws and investigating vortex-state dynamics.
Finally, some conclusions are drawn in Section \ref{conclusions}.

\section{Dynamical properties}
\label{sec:dyn}
The CGPEs \eqref{GGPE} conserve mass, energy and magnetization,
i.e.,
{\it mass}
\beas\label{mass}
{N}(t)={N}(\Psi(\cdot,t)):=\sum^{1}_{v=-1}\displaystyle\int_{\mathbb{R}^{3}} \left| \psi_{v}(\bx, t)\right|^{2}\mathrm{d} \bx \equiv {N}(\Psi(\cdot,0))=1, \quad t \ge 0, \notag
\eeas
{\it energy}
\beas
{E}(t)
 \!=\! {E}(\Psi(\cdot,t)) \!&:=&\!\!\int_{\mathbb{R}^3} \bigg[\sum^{1}_{v=-1}\left(\frac{1}{2}|\nabla\psi_v|^2+V |\psi_v|^2  - \bar{\psi}_v (\bm{\Omega} \cdot \bm{L} )  \psi_v\right) + \frac{c_0}{2}\rho^2+\frac{c_1}{2}|	\mathbf{F}|^2 \bigg]{\rm d}\bx  \\
& \equiv& {E}(\Psi(\cdot,0)), ~~ t \ge 0,
\eeas
and {\it magnetization}
\beas
{M}(t)={M}(\Psi(\cdot,t)):=\sum^{1}_{v=-1}\displaystyle\int_{\mathbb{R}^{3}} v \left| \psi_{v}(\mathbf{x},t)\right|^{2}\mathrm{d}\bx  ~ \equiv {M}(\Psi(\cdot,0)), \quad t \ge 0.
\eeas

The \textit{center of mass} of a condensate is defined as
\bea
\label{Mass-Center}
\bx_c(t)  = \int_{\mathbb{R}^3} \bx \rho ~{\rm d} \bx, \quad t\ge 0,
\eea
which satisfies the property stated in Lemma \ref{lem-CenterMass}.
\begin{lem}
\label{lem-CenterMass}
For CGPEs \eqref{GGPE} with a harmonic potential \eqref{harmonicPot}, the dynamics
of center of mass is governed by the following second-order ODE system:
\bea \label{centerODE}
&&\ddot{\mathbf{x}}_c(t) - 2 A ~ \dot{\mathbf{x}}_c(t)+ (A^2+S) \mathbf{x}_c(t) = 0, \quad t>0, \\
&&\mathbf{x}_c(0) = \int_{\mathbb{R}^3} \mathbf{x} \rho^0  {\rm d} \mathbf{x}, \quad \dot{\mathbf{x}}_c(0) =A ~\mathbf{x}_c(0)+ 
\sum^{1}_{v=-1}\int_{\mathbb{R}^3} \mathrm{Im}(\overline{\psi_{v}^{0}} \nabla \psi_{v}^{0}){\rm d} \mathbf{x}
\label{centerODEInitial}
\eea
where
$
A = \begin{psmallmatrix}
0         &  \omega_z & -\omega_y \\
-\omega_z & 0         & \omega_x \\
\omega_y  & -\omega_x & 0
\end{psmallmatrix}
$
and
$
S = {\rm diag} \{\gamma_{x}^2, \gamma_{y}^2, \gamma_{z}^2 \}.
$
\end{lem}

\begin{remark}
\label{rem:boundedness-jordan}
The boundedness of the solution $\mathbf{x}_c(t)$ to ODEs \eqref{centerODE}
is governed by initial data \eqref{centerODEInitial} and the eigenvalues $\lambda \in \mathbb{C}$ of the matrix
\[
M = \begin{pmatrix}
0 & I_3 \\
-(A^2 + S) & 2A
\end{pmatrix}.
\]
\begin{itemize}
    \item \textbf{Boundedness}: If all eigenvalues of $M$ lie in the left plane (i.e., $\operatorname{Re}(\lambda)\le 0$),
     and all eigenvalues lying on the $y$-axis are semi-simple (that is, its geometric multiplicity equals the algebraic multiplicity),
        then $\mathbf{x}_c(t)$ is bounded for any initial data \eqref{centerODEInitial} and all $t\ge 0$.
             For example, in an isotropic trap where $\gamma_x = \gamma_y = \gamma_z$, the associated matrix $M$ satisfies the above spectral properties.
    \item \textbf{Unboundedness}: Otherwise, there exists some initial data \eqref{centerODEInitial} for which the contribution from the Jordan blocks associated with eigenvalues satisfying
        $\operatorname{Re}(\lambda)>0$ or $\operatorname{Re}(\lambda)=0$ (non-semi-simple) is active, making $\mathbf{x}_c(t)$  unbounded.
\end{itemize}
\end{remark}

\

\begin{proof}
The second-order ODEs \eqref{centerODE} can  be rewritten as the first-order system $\dot{\mathbf{y}}=M\mathbf{y}$
with $\mathbf{y}=(\mathbf{x}_c,\dot{\mathbf{x}}_c)^\top$,
whose solution reads as follows:
\[
\mathbf{y}(t)=e^{tM}\mathbf{y}_0=P e^{tJ}P^{-1}\mathbf{y}_0 \quad \mbox{where}\quad \mathbf{y}_0=\mathbf{y}(0).
\]
Here
$
M = P J P^{-1},
$
$P$ is an invertible matrix and $J$ is the Jordan normal form of $M$. Specifically,
\[
J = \operatorname{diag}\bigl(J_1, J_2, \dots, J_k\bigr),
\]
where each $J_i$ is a Jordan block of size $m_i$ associated with an eigenvalue $\lambda_i\in\mathbb{C}$ of $M$ and has the form
\[
J_i =
\begin{pmatrix}
\lambda_i & 1         &        &           \\
          & \lambda_i & \ddots &           \\
          &           & \ddots & 1         \\
          &           &        & \lambda_i
\end{pmatrix}_{m_i\times m_i}.
\]

The growth of $\mathbf{y}(t)$
is determined by $\mathbf{y}_0$ and the entries of $
e^{tJ} = \operatorname{diag}\bigl(e^{tJ_1}, e^{tJ_2}, \dots, e^{tJ_k}\bigr),
$
where the matrix exponential $e^{tJ_i}$ is explicitly given by
\[
e^{tJ_i}=e^{\lambda_i t}
\begin{pmatrix}
1 & t & \frac{t^{2}}{2!} & \cdots & \frac{t^{m_i-1}}{(m_i-1)!}\\
0 & 1 & t & \cdots & \frac{t^{m_i-2}}{(m_i-2)!}\\
\vdots & & \ddots & \ddots & \vdots\\
0 & \cdots & 0 & 1 & t\\
0 & \cdots & 0 & 0 & 1
\end{pmatrix}.
\]

Specifically,
\begin{enumerate}
\item If $\operatorname{Re}(\lambda_i)<0$, the exponential decay $e^{\operatorname{Re}(\lambda_i)t}$ dominates any polynomial factor $t^{p}$ for $p<m_i$.
    Hence, the action of this Jordan block on its associated vector components decays exponentially to zero and remains bounded.

\item If $\operatorname{Re}(\lambda_i)=0$ and $m_i=1$ (i.e., the eigenvalue is \textit{semi-simple}), then $e^{tJ_i}=e^{\lambda_i t}$ with $\lambda_i\in \mathrm{i}\mathbb{R}$ and $|e^{\lambda_i t}|=1$ for all $t$.
    Hence,  the action of this Jordan block on its associated vector components remains bounded.

\item If $\operatorname{Re}(\lambda_i)=0$ and $m_i\ge 2$, the entry containing $t^{m_i-1}e^{\lambda_i t}$ grows polynomially in $t$.
Hence,  the action of this Jordan block on its associated nonzero vector components is unbounded.

\item If $\operatorname{Re}(\lambda_i)>0$, the factor $e^{\operatorname{Re}(\lambda_i)t}$ grows exponentially. Hence,
 the action of this Jordan block on its associated nonzero vector components is unbounded.
\end{enumerate}
\end{proof}

\section{Numerical method}
\label{sec:NumeMeth}
In this section, we propose high-order splitting Fourier spectral methods
to simulate the dynamics of spin-1 BECs with ArbRot \eqref{GGPE}. Due to the presence of a confining potential $V(\bx)$,
the wave function is usually smooth and decays exponentially.
Thus, we can truncate the whole space into a bounded rectangular domain $\mathbf{R}_L= [-L, L]^3$ such that the truncation error is negligible
and impose periodic boundary conditions on the wave function. The domain is discretized uniformly with mesh size $h = 2L/N$
($N \in 2 \mathbb{Z}^{+} $),
and the sets of indices and grid points are defined as
\beas
\mathcal{I}_N &:=& \left\{\textbf{j}:=(j, k, l) \in \mathbb{N}^3 \Big|~ 0 \le j, k, l \le N-1 \right\}, \\
\mathcal{K}_N &:=& \left\{\textbf{p}:=(p, q, s) \in \mathbb{Z}^3 \Big|~ -N/2 \le p, q, s \le N/2-1 \right\},\\
\mathcal{T}_N &:=& \left\{  \bx_{\textbf{j}}: = (x_j, y_k, z_{l})^{\top} := (-L+jh, -L+kh, -L+l h)^{\top}, \quad \textbf{j} \in
  \mathcal{I}_N  \right\}.
\eeas
By introducing the  Fourier basis function
\beas
M(\bx)=e^{\im \nu_p^x (x+L)}e^{\im \nu_q^y (y+L)}e^{\im \nu_s^z (z+L)},
\eeas
with $\nu_p^x = \pi p/L$, $\nu_q^y = \pi q/L$ and $\nu_s^z = \pi s/L$,
the Fourier spectral approximation \cite{Shen2011Spectral} of wave function $\psi_v~(v = -1,0,1)$  is given by
\bea
\label{Fourier}
\psi_v(\bx) \approx \sum_{\textbf{p} \in \mathcal{K}_N}
 \widehat{\psi}_{v,\textbf{p}} ~M(\bx),
\eea
with Fourier coefficients
\bea \label{FourierCoef}
\widehat{\psi}_{v,\textbf{p}}~=~ \frac{1}{N^3} \sum_{\textbf{j}\in \mathcal{I}_N} \psi_v(\bx_{\textbf{j}})~\overline{M}(\bx_{\textbf{j}}).
\eea
The above summation can be accelerated by FFT within $\mathcal{O}(3N^3 \log N)$ operations.
In fact, the approximation \eqref{Fourier} is actually an interpolation function at grid points $\bx_{\textbf{j}}\in \mathcal{T}_N$.

\

We choose a time step $\tau>0$ and denote $t_n = n \tau$ for $n = 0,1, \dots$.  To evolve from $t_n$ to $t_{n+1}$,
the CGPEs \eqref{GGPE} is split into a linear subproblem
\bea
\label{LinearSubprob}
\im \partial_t \Psi(\bx,t) = \Big[-\frac{1}{2} \Delta -  \bm{\Omega} \cdot \bm{L} \Big] \Psi(\bx, t) := \mathcal{H}_{\rm{lin}} \Psi(\bx,t),
\eea
and a nonlinear subproblem
\bea
\label{NonLinearSubprob}
i\partial_t \Psi(\bx,t) = \big[(V+c_0 \rho)  I_{3}+c_1 \mathbf{F}\cdot \mathbf{f} \big]\Psi(\bx,t) := \mathcal{H}_{\rm{non}} \Psi(\bx,t).
\eea
Given the initial value $\Psi^{n}(\bx)\!:=\Psi(\bx,t_n)$, 
the solutions to linear and nonlinear subproblems at time $t_{n+1}$ are denoted by $\Psi^{n+1} = e^{-\im \tau \mathcal{H}_{\rm{lin}}} \Psi^{n}$ and $\Psi^{n+1} = e^{-\im \tau \mathcal{H}_{\rm{non}}} \Psi^{n}$ respectively.
Thus the  high-order splitting schemes can be constructed as
\cite{YoshHigherInteg}
\bea\label{splitting}
\Psi^{n+1} = \Big[\prod_{j=1}^{m}e^{-\im a_{j} \tau \mathcal{H}_{\rm{non}}} e^{-\im b_{j} \tau \mathcal{H}_{\rm{lin}}} \Big] \Psi^{n},
\eea
where $a_j, b_j\in \mathbb C$ are suitably chosen coefficients.
For the second-order Strang splitting scheme ($m =2$), the coefficients are chosen as $a_1 = a_2 = \frac{1}{2}$, $b_1 = 1$  and  $b_2 = 0$.
For the fourth-order Yoshida splitting scheme ($m=4$),
the coefficients are given by
$ a_1 = a_4 = \frac{1}{2(2-2^{1/3})}$,  $a_2 = a_3 = \frac{1-2^{1/3}}{2(2-2^{1/3})}$,
$ b_1 = b_3 = \frac{1}{2-2^{1/3}}$, $ b_2 = -\frac{2^{1/3}}{2-2^{1/3}}$ and $b_4 = 0$.
As for the computation of subproblems \eqref{LinearSubprob} and \eqref{NonLinearSubprob}, we provide the details in the following subsection.


\subsection{Numerical integration of the linear subproblem}

It is crucial to note that the rotation term indeed can be reformulated as follows
\bea
\label{matrixK}
\bm{\Omega} \cdot \bm{L}=-iK\bx\cdot\nabla, \qquad K= \scalebox{1.1}{$\begin{psmallmatrix} 0& -\omega_z&\omega_y\\ \omega_z &0& -\omega_x \\-\omega_y &\omega_x &0\\\end{psmallmatrix}$}.
\eea
To handle such rotation term, inspired by the RLC method \cite{RotLagBao},  we introduce the rotating coordinate  transformation
\begin{equation} \label{fun-rot-def2}
\tilde{\bx}=e^{(t-t_n)K}\bx:=\mathcal{R}(t-t_n)\bx, ~~
\phi_v (\tilde{\bx},t)  := \psi_v(\mathcal{R}^{-1}(t-t_n)\, \tilde{\bx},t)=\psi_v(\bx,t),
\end{equation}
with $v=1,0,-1$. The matrix $K$ is antisymmetric, and therefore, it is easy to verify that $\mathcal{R}^{-1}(t) = \mathcal{R}^{\top}(t)$ for any time $t$.
Using the chain rule, we derive
\beas
      &&\partial_t \phi_v(\tilde{\textbf{x}},t)=\partial_t \psi_v(\textbf{x},t) - \nabla\psi_v(\textbf{x},t)\cdot(K \bx),\\
      &&\nabla\phi_v(\tilde{\textbf{x}},t) = \mathcal{R}(t-t_n) \nabla\psi_v(\textbf{x},t) , \quad \Delta\phi_v(\tilde{\textbf{x}},t) =\Delta\psi_v(\textbf{x},t),
\eeas
and the linear subproblem \eqref{LinearSubprob} is then reformulated as follows
\bea
\label{LinearSubprobtran}
\left\{\begin{array}{l}
\im \partial_t \Phi(\tilde{\bx},t) = -\frac{1}{2} \Delta ~ \Phi(\tilde{\bx}, t),  \quad t_n \leq t\leq t_{n+1},
 \\[0.5em]
\Phi(\tilde{\bx}, t_n) = \Phi^{n}(\tilde{\bx}),  \\
\end{array}\right.
\eea
where $\Phi = (\phi_1, \phi_0, \phi_{-1})^{\top}$.
Compared to subproblem \eqref{LinearSubprob}, Eq. \eqref{LinearSubprobtran} does not contain the
ArbRot term, which allows  us to develop simple and efficient numerical methods for
simulating the dynamics.

For each component $\phi_v(\tilde{\bx},t)$, we approximate it using the Fourier approximation \eqref{Fourier}. Substituting this into Eq. \eqref{LinearSubprobtran}, we obtain
\bea
i\partial_t \widehat{\phi}_{v,\textbf{p}}(t)
=Q~\widehat{\phi}_{v,\textbf{p}}(t), \quad Q=\big[(\nu_p^{\tilde{x}})^2+(\nu_q^{\tilde{y}})^2+(\nu_s^{\tilde{z}})^2\big]/2,
\eea
and thus
\bea\label{SolLinear}
\phi_v(\tilde{\bx},t)=\sum_{\textbf{p} \in \mathcal{K}_N}
 e^{-i(t-t_n)Q}\widehat{\phi}_{v,\textbf{p}}(t_n) ~M(\tilde{\bx}).
\eea

It is worth noting that the rotation matrix $\mathcal{R}(t-t_n)$
ensures $\tilde{\bx}=\bx$ at time $t=t_n$, thus  $\Phi^{n}(\tilde{\bx})=\Psi^{n}(\tilde{\bx})$ in \eqref{LinearSubprobtran}.
Hence, apart from the computation of  \eqref{SolLinear}, the only additional operation required to evaluate original linear subproblem
\eqref{LinearSubprob} in practice is the implementation of a rotating coordinate transformation at $t=t_{n+1}$, i.e.,
\bea\label{rotCoordinateTransform}
\psi^{n+1}_v(\bx)=\phi^{n+1}_v(\mathcal{R}(\tau)\bx).
\eea

Given the grid function $\phi_v^{n+1}$, a straightforward method  to evaluate $\psi_v^{n+1}$ \eqref{rotCoordinateTransform} accurately and efficiently is to adopt the Fourier spectral interpolation.
However, since the target points $\mathcal{R}(\tau)\bx$ are no longer uniformly distributed in each spatial direction, the FFT cannot be applied directly,
making the computation time-consuming.
To further enhance efficiency, an accurate and efficient implementation of this transformation is of essential importance for the simulation.

\subsubsection{Fast rotation based on shear decomposition}
In fact,
the rotating coordinate transformation \eqref{rotCoordinateTransform}
can be reformulated into a computation-friendly form defined as the action of exponential operators on a function.
Furthermore, such reformulation allows for an efficient implementation that can be realized purely by one-dimensional FFT/iFFT.

Specifically, for any $d \times d$ matrix $B$,
the identity
\be\label{iden_te}
u(e^{tB} \bx) = e^{tB\mathbf{x} \cdot \nabla} u(\mathbf{x})
\ee
holds. This follows from the transport equation
\begin{equation}
    \begin{cases}
        \partial_t w = B\mathbf{x} \cdot \nabla w, ~~ \mathbf{x} \in \mathbb{R}^d, \\
        w(\mathbf{x},t=0) = u(\mathbf{x}),
    \end{cases}
\end{equation}
which admits a unique solution $w(\bx,t)=u(e^{tB} \bx)$ and the solution can also be equivalently written as
$
w(\bx,t) = e^{tB\mathbf{x} \cdot \nabla} u(\mathbf{x}).
$

Taking $t=1,d=3$, and substituting $u(\mathbf{x}) = \phi^{n+1}_v(\mathbf{x})$ and  $B = \tau K$ ($K$ is shown in \eqref{matrixK}) into Eq. \eqref{iden_te}, we can rewrite coordinate transformation \eqref{rotCoordinateTransform} as
\be\phantomsection \label{OperFun-RotMap}
\psi^{n+1}_v(\mathbf{x})
= \phi^{n+1}_v\left(e^{\tau K}\mathbf{x}\right)
= e^{\tau K\mathbf{x}\cdot\nabla}\phi^{n+1}_v(\mathbf{x}).
\ee

However, efficient implementation of the last formula is non-trivial.
When $\omega_x \!= \!\omega_y \!= \!0$, to address this challenge, Liu et al. \cite{Threeshear} proposed the Rotation-Shear-Decomposition-Acceleration (RSDA) method based on the shear decomposition.
Inspired by this, we turn to looking for suitable shear decomposition of the arbitrary-angle rotation matrix
\bea \label{MatrDecomExp}
\qquad\quad \mathcal{R}(\tau)
&=&e^{\tau K}=
\scalebox{1.3}{$
\begin{psmallmatrix} \cos(|\bm{\Omega}|\tau)+\omega_x^2 q & -\omega_z p+\omega_x\omega_y q & \omega_y p+\omega_x\omega_z q \\[0.2em]
\omega_z p+\omega_x\omega_y q  & \cos(|\bm{\Omega}|\tau)+\omega_y^2 q & -\omega_x p+\omega_y\omega_z q  \\[0.2em]
-\omega_y p+\omega_x\omega_z q & \omega_x p+\omega_y\omega_z q  &\cos(|\bm{\Omega}|\tau)+\omega_z^2 q  \\\end{psmallmatrix}$}  \\ [0.5em] \notag
&=&
\scalebox{1.3}{$\begin{psmallmatrix}
1 & \lambda_1 & \lambda_2\\
0& 1  & 0\\
0& 0 & 1 \\
\end{psmallmatrix}$}
\scalebox{1.3}{$\begin{psmallmatrix}
1 & 0 & 0\\
\lambda_3& 1  & \lambda_4\\
0& 0 & 1 \\
\end{psmallmatrix}$}
\scalebox{1.3}{$\begin{psmallmatrix}
1 & 0 & 0\\
0& 1  & 0\\
\lambda_5& \lambda_6 & 1 \\
\end{psmallmatrix}$}
\scalebox{1.3}{$\begin{psmallmatrix}
1 & \lambda_7 & \lambda_8\\
0& 1  & 0\\
0& 0 & 1 \\
\end{psmallmatrix}$}
= e^{M_1} e^{M_2} e^{M_3} e^{M_4},
\eea
where $p=\sin(|\bm{\Omega}|\tau)/|\bm{\Omega}|$, $q=(1-\cos(|\bm{\Omega}|\tau))/|\bm{\Omega}|^2$ and $M_1,\cdots,M_4$ are all nilpotent matrices defined below
\beas
M_{1} = \scalebox{1.3}{$\begin{psmallmatrix}
0 & \lambda_1 & \lambda_2\\
0& 0  & 0\\
0& 0 & 0 \\
\end{psmallmatrix}$}, \quad
M_{2} = \scalebox{1.3}{$\begin{psmallmatrix}
0 & 0 & 0\\
\lambda_3& 0  & \lambda_4\\
0& 0 & 0 \\
\end{psmallmatrix}$}, \quad
M_{3} = \scalebox{1.3}{$\begin{psmallmatrix}
0 & 0 & 0\\
0& 0  & 0\\
\lambda_5& \lambda_6 & 0 \\
\end{psmallmatrix}$}, \quad
M_{4} = \scalebox{1.3}{$\begin{psmallmatrix}
0 & \lambda_7 & \lambda_8\\
0& 0  & 0\\
0& 0 & 0 \\
\end{psmallmatrix}$},
\eeas
with coefficients $\lambda_s,~ s=1,\cdots, 8$.
Matching both sides of Eq. \eqref{MatrDecomExp} entry by entry,
the coefficients are given by
\bea
\label{CoeeExt1}
&& \qquad~~
\lambda_1 = \frac{(n_x^2 + n_z^2) s}{g(n_x, n_y, n_z)},~~
\lambda_2 = \frac{-g_1(n_x, n_y, n_z) s}{g(n_x, n_y, n_z) ~ g(n_x, n_z, n_y)},~~
\lambda_3 = -2 s g(n_x, n_y, n_z)
,  \\[0.5em]
&& \qquad~~
\lambda_4 = \frac{2 n_x n_y s}{g(n_x, n_z, n_y)},~~
\lambda_5 = 2sg(n_x, n_z, n_y), \qquad \qquad \quad~~
\lambda_6 = \frac{-2 n_x n_z s}{g(n_x, n_y, n_z)},
\label{CoeeExt2}
\\[0.5em]
&& \qquad~~
\lambda_7 = \frac{g_1(n_x, n_z, n_y) s}{g(n_x, n_y, n_z) ~ g(n_x, n_z, n_y)}, \quad
\lambda_8 = \frac{-(n_x^2 + n_y^2) s}{g(n_x, n_z, n_y)},
\label{CoeeExt3}
\eea
where
$n_\alpha := \omega_\alpha/|\bm{\Omega} |$ for $\alpha = x, y, z$, $n_x^2+n_y^2+n_z^2=1$,
$g(n_x, n_y, n_z) := n_x n_y s   - n_z m,$
~
$g_1(n_x, n_y, n_z) := n_x n_y\left(1+n_z^2\right)s + n_z\left(n_x^2 - n_y^2\right)m$
with
$m = \cos(|\bm{\Omega} | \tau/2)$ and  $s = \sin(|\bm{\Omega} | \tau/2)$.

\

Based on the shear decomposition \eqref{MatrDecomExp}, the rotating coordinate transformation \eqref{OperFun-RotMap} can be computed as
\bea\notag
\psi_v^{n+1}(\bx)
&=& \phi_v^{n+1}\left(e^{M_1} e^{M_2} e^{M_3} e^{M_4}\bx\right) \\[0.1em]
& =&  e^{M_4 \bx \cdot \nabla} e^{M_3 \bx \cdot \nabla} e^{M_2 \bx \cdot \nabla} e^{M_1 \bx \cdot \nabla}\phi_v^{n+1}(\bx) \notag\\[0.6em]
&=&e^{\left(\lambda_7y+\lambda_8 z\right) \partial_{x} } e^{\left(\lambda_5x+\lambda_6 y\right) \partial_{z}  } e^{\left(\lambda_3x+\lambda_4 z\right) \partial_{y} } e^{\left(\lambda_1y+\lambda_2 z\right) \partial_{x} } \phi_v^{n+1}(\bx).
\label{RotComptForwfourshear}
\eea
The sub-operators in \eqref{RotComptForwfourshear} are well approximated using Fourier spectral method, which reads as follows
\begin{small}
\bea \label{Rotexe}
\phi_{v,~\!\textbf{j}}^{(1)}
&=& \sum_{p\in \mathcal{J}}   e^{ \left(\lambda_1y_k+\lambda_2 z_l\right)  ( \im \nu_p^x)}  ~\Big(\widehat{\phi_{v, kl}^{n+1}}\Big)^x_{p}  ~ e^{\im \nu_p^x (x_j+L)}, \\
\phi_{v,~\!\textbf{j}}^{(2)}
&=& \sum_{q\in \mathcal{J}}  e^{ \left(\lambda_3x_j+\lambda_4 z_l\right) ( \im \nu_q^y)} ~ \Big(\widehat{\phi_{v,jl}^{(1)}}\Big)^y_{q} ~ e^{\im \nu_q^y (y_k+L)}, \notag \\
\phi_{v,~\!\textbf{j}}^{(3)}
&=& \sum_{s\in \mathcal{J}}  e^{ \left(\lambda_5x_j+\lambda_6 y_k\right) ( \im \nu_s^z)} ~ \Big(\widehat{\phi_{v,jk}^{(2)}}\Big)^z_{s} ~ e^{\im \nu_s^z (z_l+L)}, \notag \\
\psi_{v,~\!\textbf{j}}^{n+1}
&=& \sum_{p\in \mathcal{J}}   e^{ \left(\lambda_7y_k+\lambda_8 z_l\right)  ( \im \nu_p^x)}  ~\Big(\widehat{\phi_{v, kl}^{(3)}}\Big)^x_{p}  ~ e^{\im \nu_p^x (x_j+L)}, \notag
\eea
\end{small}
where  $\mathcal{J} = \{-N/2, \dots, N/2-1\}$ and $\psi_{v,~\!\textbf{j}}^{n+1} $ denotes the numerical approximation of $\psi_v(x_j, y_k, z_{l},t_{n+1})$.
Here, $\widehat{f}^x_p$, $\widehat{f}^y_q$ and  $\widehat{f}^z_s$ represent the Fourier coefficients of $f$ along the $x$, $y$ and $z$ directions respectively.
The above computation \eqref{Rotexe} can be implemented by one-dimensional FFT/iFFT.

\begin{remark}
All the exponential coefficients in \eqref{Rotexe} depend only on time step $\tau$ and rotation vector $\bm{\Omega}$,
therefore, they can be computed  once for all in the precomputation step.
\end{remark}

In summary, we can solve linear subproblem \eqref{LinearSubprob} by first obtaining the solution to transformed equation \eqref{LinearSubprobtran} and then applying rotating coordinate transformation \eqref{rotCoordinateTransform}.
With the fast rotation algorithm \eqref{RotComptForwfourshear},
the solution to \eqref{LinearSubprob} can be expressed  as
\begin{equation}
\label{ExtSplit3D}
\begin{aligned}
\psi^{n+1}_v(\mathbf{x})
&= e^{i\tau\left[ \bm{\Omega} \cdot \mathbf{L} + \Delta/2 \right]}\psi^{n}_v(\mathbf{x})
=e^{i\tau \bm{\Omega} \cdot  \bm{L}}{e^{i\tau\Delta/2}}\psi^{n}_v(\mathbf{x})\\
&= e^{\left(\lambda_7y+\lambda_8 z\right) \partial_{x}}
  e^{\left(\lambda_5x+\lambda_6 y\right) \partial_{z}}
  e^{\left(\lambda_3x+\lambda_4 z\right) \partial_{y}}
  e^{\left(\lambda_1y+\lambda_2 z\right) \partial_{x}}
  e^{i\tau \Delta/2}
  \psi_v^n(\mathbf{x}).
\end{aligned}
\end{equation}

\subsubsection{Singularities induced by shear decomposition}

The scheme \eqref{ExtSplit3D} exhibits singularities at $\tau_s$ when either
$g(n_x,n_y,n_z)=0$ or $g(n_x,n_z,n_y)=0$,
where some coefficients in \eqref{CoeeExt1}-\eqref{CoeeExt3} blow up as $\tau \rightarrow\tau_s$.
Without loss of generality, we assume $\omega_x$, $\omega_y$, $\omega_z\geq0$.
Indeed, for any given $\bm{\Omega}$, the problem \eqref{LinearSubprob} can be reduced, through a suitable variable substitution, to the case where all components of the rotation vector
$\bm{\Omega}$ are non-negative.
Table \ref{tab_substitution} in Appendix \ref{substitution} summarizes the specific variable substitutions for various $\bm{\Omega}$, respectively.
Since $\mathcal{R}(\tau)$ \eqref{MatrDecomExp} is $2\pi/|\bm{\Omega}|$-periodic in $\tau$,
we can restrict our discussion of singularities $\tau_s$ to the interval $[0,2\pi]/|\bm{\Omega}|$.

Setting $g(n_x, n_y, n_z) = 0$,
we obtain
\beas
n_x n_y s  - n_z m =
\sqrt{{n_x^2 n_y^2} + {n_z^2}}~ \sin\Big( {|\bm{\Omega} |\tau}/{2} - \arctan\big({n_z^2}/{w}\big) \Big) = 0,
\label{eq:condition_5}
\eeas
where $w:=n_x n_y n_z$ and  we set $\arctan(+\infty)=\pi/2$ by continuous extension, that is,
\bea
\sqrt{{n_x^2 n_y^2} + {n_z^2}}=0 \qquad \mbox{or} \qquad \sin\Big( {|\bm{\Omega} |\tau}/{2} - \arctan\big({n_z^2}/{w}\big) \Big)=0.
\eea
Similarly, by setting $g(n_x, n_z, n_y) = 0$,
we obtain
\bea
\sqrt{{n_x^2 n_z^2} + {n_y^2}}=0 \qquad \mbox{or} \qquad \sin\Big( {|\bm{\Omega} |\tau}/{2} - \arctan\big({n_y^2}/{w}\big) \Big)=0.
\eea
If rotation parameter $\bm{\Omega}$ satisfies $\sqrt{{n_x^2 n_y^2} + {n_z^2}}=0$  or $\sqrt{{n_x^2 n_z^2} + {n_y^2}}=0$,
i.e., when any two of $\omega_x$, $\omega_y$, $\omega_z$ are zero,
some coefficients in \eqref{CoeeExt1}-\eqref{CoeeExt3}  blow up for almost all $\tau$.
In fact, the rotation  in such cases reduces to a rotation about a coordinate axis. 
For example, when $\bm{\Omega}=(0,0,\omega_z)^\top$, the rotation term becomes $-\omega_z L_z$, i.e.,
a classical rotation about the $z$-axis,
which can be treated using the method in \cite{Fshear}.
The cases with $\bm{\Omega}=(\omega_x,0,0)^\top$ and $\bm{\Omega}=(0,\omega_y, 0)^\top$ can be handled analogously.
When $\sqrt{{n_x^2 n_y^2} + {n_z^2}}\neq0$  and $\sqrt{{n_x^2 n_z^2} + {n_y^2}}\neq0$, 
i.e., $\bm{\Omega}$ has at least two nonzero components,
the singularities $\tau_s$ are exactly given by
$\frac{2}{|\bm{\Omega}|}\arctan\!\left( {n_z^2}/{\omega}\right)$ and $\frac{2}{|\bm{\Omega}|} \arctan\!\left( {n_y^2}/{\omega} \right)$.

We define
\[
f_1 = 2\arctan\!\left({n_x^2}/{\omega}\right),\quad
f_2 = 2\arctan\!\left({n_y^2}/{\omega}\right),\quad
f_3 = 2\arctan\!\left({n_z^2}/{\omega}\right).
\]
When restricted to the interval $[0,2\pi]$,  we have $0\leq f_1,~f_2,~f_3\leq\pi$ and $f_1+f_2+f_3\neq0$.
Considering that
\beas
\tan(f_1/2 + f_2/2+ f_3/2)
=0,
\eeas
we obtain
\bea
f_1 + f_2 + f_3 = 2\pi.
\eea
This further yields
\bea
{\mathrm{median}\{f_1,f_2,f_3\}}\in [\pi/2, \pi] \quad  \mbox{and} \quad
{\max\{f_1,f_2,f_3\}}\in [\pi/2, \pi].
\eea

Therefore, we state the following lemma:
\begin{lem}\label{tauc}
When $\bm{\Omega}$ has at least two nonzero components and $\omega_x$  is the smallest among $\omega_{\alpha}$ for $\alpha=x,y,z$,
the singularities $\tau_s$  in scheme \eqref{ExtSplit3D} are  given by
\be\label{tau_c_five}
\left\{ \frac{2}{|\bm{\Omega}|}\arctan\!\left( {n_z^2}/{w} \right),
\; 
\frac{2}{|\bm{\Omega}|}
\arctan\!\left( {n_y^2}/{w} \right) \right\} \subset [\pi/2, \pi]/|\bm{\Omega}|.
\ee
\end{lem}

For the cases where $\omega_x$ is not the smallest among $\omega_{\alpha}$ ($\alpha=x,y,z$), we can obtain the following two schemes by adopting variable substitutions to \eqref{LinearSubprob}.
When $\omega_z$ is the smallest, the scheme is given by
\bea
\label{ExtSplit3D_shiftzxy}
\qquad \quad \psi^{n+1}_v(\bx)
&=&
e^{\left(\tilde{\lambda}_7x+\tilde{\lambda}_8 y\right) \partial_{z} }
e^{\left(\tilde{\lambda}_5z+\tilde{\lambda}_6 x\right) \partial_{y} }
e^{\left(\tilde{\lambda}_3z+\tilde{\lambda}_4 y\right) \partial_{x} }
e^{\left(\tilde{\lambda}_1x+\tilde{\lambda}_2 y\right) \partial_{z} }
e^{\im\tau \Delta/2 }
\psi_v^n(\bx),
\eea
which is obtained by taking
$\tilde{\bx}  = P_1 \bx =(z,x,y)^\top, ~\tilde{\bm{\Omega}}  = P_1 \bm{\Omega} = (\omega_z,\omega_x,\omega_y)^\top$
and
$\tilde{\bm{\lambda}}
=
(\tilde{\lambda}_1, \tilde{\lambda}_2, \tilde{\lambda}_3,\tilde{\lambda}_4, \tilde{\lambda}_5,$
$ \tilde{\lambda}_6, \tilde{\lambda}_7,\tilde{\lambda}_8 )
=
\bm{\lambda}(\tilde{\bm{\Omega}})$. The associated $\tau_s$ are given by
\be
\label{tau_c_five_zxy}
\left\{\frac{2}{|\bm{\Omega}|}
\arctan\!\left( {n_y^2}/{w} \right),~
\frac{2}{|\bm{\Omega}|}
\arctan\!\left( {n_x^2}/{w} \right)\right\} \subset [\pi/2, \pi]/|\bm{\Omega}|.
\ee
When $\omega_y$ is the smallest, the scheme reads as
\bea
\label{ExtSplit3D_shiftyzx}
\qquad\quad \psi^{n+1}_v(\bx)
&=&
e^{\left(\tilde{\lambda}_7z+\tilde{\lambda}_8 x\right) \partial_{y} }
e^{\left(\tilde{\lambda}_5y+\tilde{\lambda}_6 z\right) \partial_{x} }
e^{\left(\tilde{\lambda}_3y+\tilde{\lambda}_4 x\right) \partial_{z} }
e^{\left(\tilde{\lambda}_1z+\tilde{\lambda}_2 x\right) \partial_{y} }
e^{\im\tau \Delta/2 }
\psi_v^n(\bx),
\eea
which is obtained by taking
$\tilde{\bx}  = P_2 \bx =(y,z,x)^\top$,
$\tilde{\bm{\Omega}}= P_2\bm{\Omega}$ and the corresponding $\tau_s$ are given by
\be
\left\{\frac{2}{|\bm{\Omega}|}
\arctan\!\left( {n_x^2}/{w} \right),~
\frac{2}{|\bm{\Omega}|}
\arctan\!\left( {n_z^2}/{w} \right)\right\} \subset [\pi/2, \pi]/|\bm{\Omega}|.
\ee
Thus, for a given $\bm{\Omega}$ with at least two nonzero components,  we can always select an appropriate scheme that has singularities $\tau_s \subset \left[{\pi}/{2},\,\pi\right]/|\bm{\Omega}|$,
as detailed in Algorithm  \ref{alg:selection}.
As $\tau$ approaches the singularities $\tau_s$,
   the scheme suffers a significant accuracy loss. To address this issue,  we introduce  the time step decomposition technique in subsection \ref{timedecomp}.
\begin{algorithm}
\caption{Selection of the optimal scheme}
\label{alg:selection}
\begin{algorithmic}[0]
\Require The rotation vector $\bm{\Omega} = (\omega_x, \omega_y, \omega_z)^\top$ with $\omega_x$, $\omega_y$, $\omega_z\geq0$.
\vspace{0.5em}
\State Compute $m = \min\{\omega_x,\ \omega_y,\ \omega_z\}$.
\vspace{0.5em}
\If{$m = \omega_x$} {select scheme} \eqref{ExtSplit3D}
\vspace{0.3em}
\ElsIf{$m = \omega_y$} {select scheme} \eqref{ExtSplit3D_shiftyzx}
\vspace{0.3em}
\Else ~~{select scheme} \eqref{ExtSplit3D_shiftzxy}
\EndIf
\vspace{0.5em}
\Ensure The optimal scheme.
\end{algorithmic}
\end{algorithm}

\

\subsubsection{Extension to arbitrary time step}\label{timedecomp}
Although the numerical results  become inaccurate when $\tau$  approaches the singularities $\tau_s$,
this limitation for schemes \eqref{ExtSplit3D}, \eqref{ExtSplit3D_shiftzxy} and \eqref{ExtSplit3D_shiftyzx} can be effectively overcome using a time step decomposition technique.

Based on the  periodicity of $\mathcal{R}(\tau)$,
we first decompose the time step as
\bea
\label{decom1}
\tau =   k_1 \frac{2\pi}{|\bm{\Omega}|}+ \tau_1 \quad \mbox{with} \quad  0\leq \tau_1 <\frac{2\pi}{|\bm{\Omega}|},
\eea
where $k_1=\lfloor \frac{\tau|\bm{\Omega}|}{2\pi} \rfloor  \in \mathbb{Z}$ and $\lfloor \cdot \rfloor $ denotes the floor function.
Consequently, it is sufficient to consider $\tau_1 \in [0,{2\pi})/|\bm{\Omega}|$.
If $\tau_1=0$, we obtain $\mathcal{R}(\tau)=I_3$, which implies that no rotation is required.
Our method is mathematically correct when $\tau_1 \neq \tau_{s}$,
but its numerical accuracy deteriorates severely as $\tau_1$ approaches $\tau_{s}$.
Given that the singularities $\tau_{s}$ lie in $\left[{\pi}/{2},\,\pi\right]/|\bm{\Omega}|$, we numerically extend this interval outward by $\frac{\pi}{2|\bm{\Omega}|}$ on each side for robustness.
Then we denote the resulting extended interval $I_f:=(0,{3\pi}/{2})/|\bm{\Omega}|$
as the ``failure interval", and $I_s=[{3\pi}/{2},2\pi]/|\bm{\Omega}|$ as  the ``security interval".
When $\tau_1 \!\in \!I_s$, the rotating coordinate transformation can be evaluated directly.
Since $\mathcal{R}(\tau)$ is a rotation matrix with rotation angle $|\bm{\Omega}|\tau$,
the cases where $\tau_1 \!\in \!I_f$ can be accurately implemented by making multiple calls to the case $\tau_1 \in I_s$.

Specifically, if $\tau_1 \in I_s$,  we obtain
\be\label{tau_1_operator}
e^{i\tau[\bm{\Omega} \cdot \bm{L}+\Delta/2 ]}
= e^{i(k_1\frac{2\pi}{|\bm{\Omega}|} + \tau_1)[\bm{\Omega} \cdot \bm{L}+\Delta/2 ]}
=e^{i\tau_1 \bm{\Omega} \cdot  \bm{L}}{e^{i\tau\Delta/2}}.
\ee
Otherwise, we denote $\tau_2=\tau_1-2\pi/|\bm{\Omega}|$ and we can always decompose  $\tau_2$ as
\be\label{decom2_EEI5}
\tau_2 = k_2
\tau_{\star} + \tilde{\tau}, \qquad \tilde{\tau} \in I_s-2\pi/|\bm{\Omega}|,
\ee
where $\tau_{\star} = \frac{-\pi}{2|\bm{\Omega}|}$.
Thus we obtain
\be
e^{i\tau[\Delta/2 + \bm{\Omega} \cdot  \bm{L}]}
=
\underbrace{e^{i\tau_{\star}\bm{\Omega} \cdot  \bm{L}} \cdots e^{i\tau_{\star}\bm{\Omega} \cdot  \bm{L}}}_{k_2 \text{ operators}}
e^{i\tilde{\tau}\bm{\Omega} \cdot  \bm{L}}
{e^{i\tau\Delta/2}},\label{tau_2_operator}
\ee
with $1\leq k_2 \leq 3$.
To sum up, from time $t_n$ to $t_{n+1}$,
we can solve the linear subproblem \eqref{LinearSubprob} by using  Algorithm \ref{LinearSol}.
\begin{algorithm}
\caption{Fast solver for the linear subproblem \eqref{LinearSubprob}}
\label{LinearSol}
\begin{algorithmic}[0]
\Require Wave function $\psi^n_{v,~\!\textbf{j}}$, time step $\tau$ and rotation vector $\bm{\Omega}$.
\vspace{0.15cm} 
\State \textbf{1.} Apply Algorithm~\ref{alg:selection} and select the optimal scheme.
\State \textbf{2.} Decompose $\tau$ according to Eq.~\eqref{decom1} to obtain the residual step $\tau_1$.
\begin{itemize}
                        \item If $\tau_1 \in I_s$, compute $\psi^{n+1}_{v,~\!\textbf{j}}$ via Eq.~\eqref{tau_1_operator}.
                         \item If $\tau_1 \in I_f$, compute  $\psi^{n+1}_{v,~\!\textbf{j}}$ via Eq.~\eqref{tau_2_operator}.
\end{itemize}
\vspace{0.1cm}
\Ensure Updated wave function $\psi^{n+1}_{v~\!\textbf{j}}$.
\end{algorithmic}
\end{algorithm}

\begin{thm}\label{thm:applicability}
For any time step $\tau$, the proposed fast algorithm for the rotating coordinate transformation is spectrally accurate.
As a solver for the linear subproblem \eqref{LinearSubprob}, Algorithm \ref{LinearSol} is temporally exact and spectrally accurate in space.
\end{thm}

\

Here we verify the accuracy of Algorithm \ref{LinearSol} 
by comparing it with a solver that does not employ the time step decomposition technique. For convenience,
we refer to Algorithm \ref{LinearSol} as ``improved scheme" and the other as ``original scheme".
The numerical error is calculated as follows
\beas
    e_v^{h,\tau}(t) &:=& \lVert \psi_v^{\text{ref}}(t)-\psi_v^{h,\tau}(t)\rVert_{2} ~/~\lVert \psi_v^{\text{ref}}(t) \rVert_{2},
\eeas
where  $\|\cdot \|_2$ is the standard $l_2$-norm,
$\psi_v^{h,\tau}(t)$ is the numerical approximation at time $t$ obtained with the mesh size $h$ and time step $\tau$, and $\psi_v^{\text{ref}}(t)$ is the reference solution.
All subsequent error calculations will adhere to this definition.
For Eq. \eqref{LinearSubprob}, we consider the following initial value
$
 \psi_v^0(\mathbf{x}) = \frac{1}{\pi^{3/4}} ~x^2~ e^{-\frac{1}{2}|\mathbf{x}|^2},
$
which admits an exact solution
\cite{CompactSplitRotSpin1}
 \[
\psi_v(\mathbf{x},t)= \frac{1}{\pi^{3/4}}~\frac{x_1^2(t)-t^2+\im t}{(1+\im t)^{7/2}}~ e^{-\frac{1}{2 (\im t+1)}|\mathbf{x}|^2},
\]
with $x_1(t) = x\cos(|\bm{\Omega}| t)
+ \frac{\sin(|\bm{\Omega}| t)}{|\bm{\Omega}|} (\omega_y z-\omega_z y)+ \frac{1-\cos(|\bm{\Omega}| t)}{|\bm{\Omega}|^2} \omega_x \bm{\Omega} \cdot \mathbf{x} $.

\begin{exmp}
\label{Exmp:shear001}
For the linear problem \eqref{LinearSubprob}, we solve it for the parameters
$\bm{\Omega} = (0,10^{-10},10)^{\top}$ $(\mathbf{Case ~\uppercase\expandafter{\romannumeral1}})$,~
$\bm{\Omega} = (0,8,10)^{\top} (\mathbf{Case ~\uppercase\expandafter{\romannumeral2}})$ ~and~
$\bm{\Omega} = (6,8,10)^{\top} (\mathbf{Case ~\uppercase\expandafter{\romannumeral3}})$.
\end{exmp}

\begin{figure}[h]
    \centering
  \includegraphics[width=6cm, height=0.6cm]{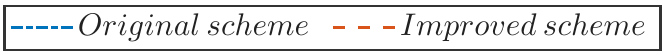}\\
  \centering
  \includegraphics[width=4.4cm, height=3.3cm]{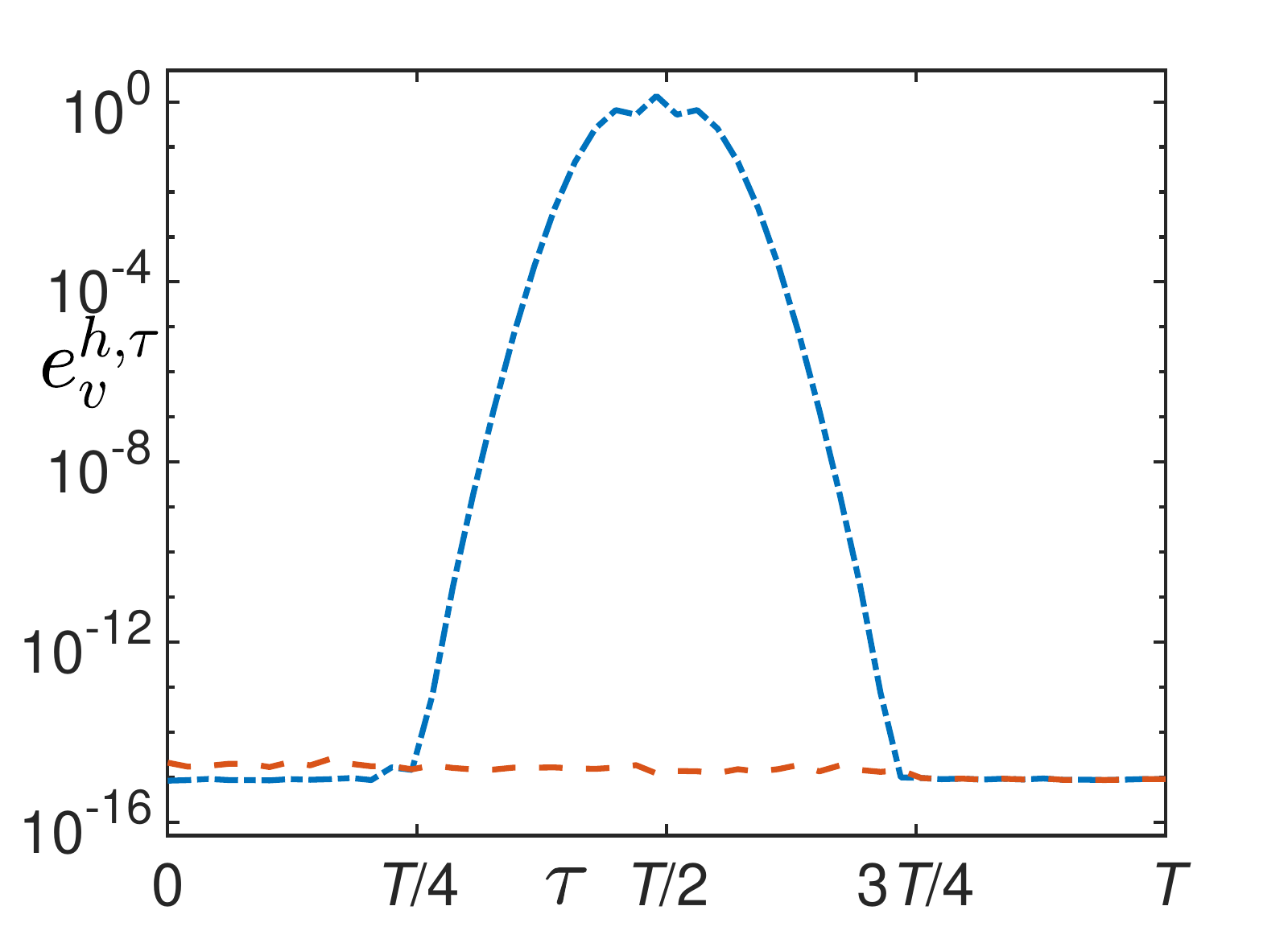}
  \hspace{-0.4cm}
    \includegraphics[width=4.4cm, height=3.3cm]{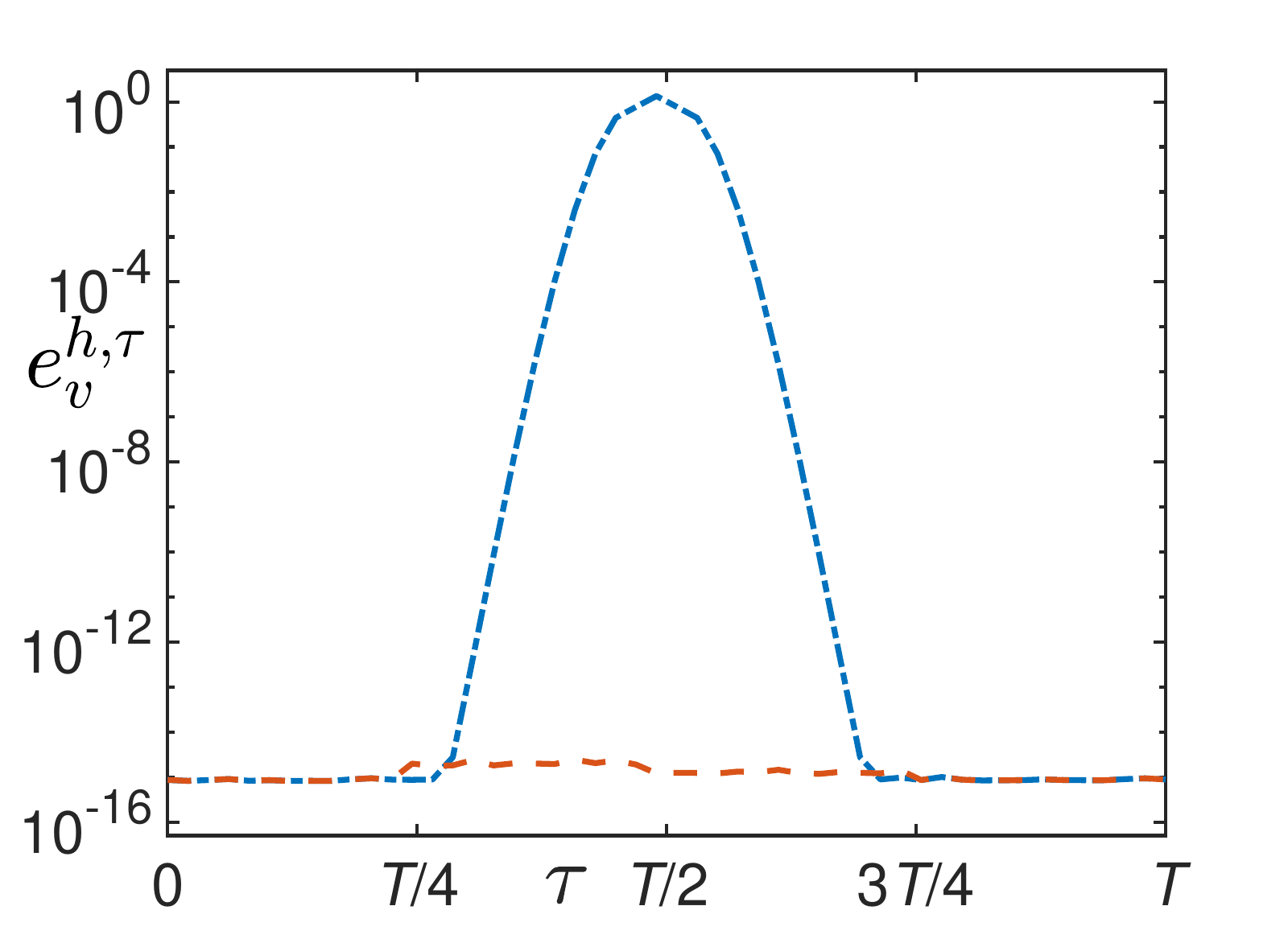}
  \hspace{-0.4cm}
  \includegraphics[width=4.4cm, height=3.3cm]{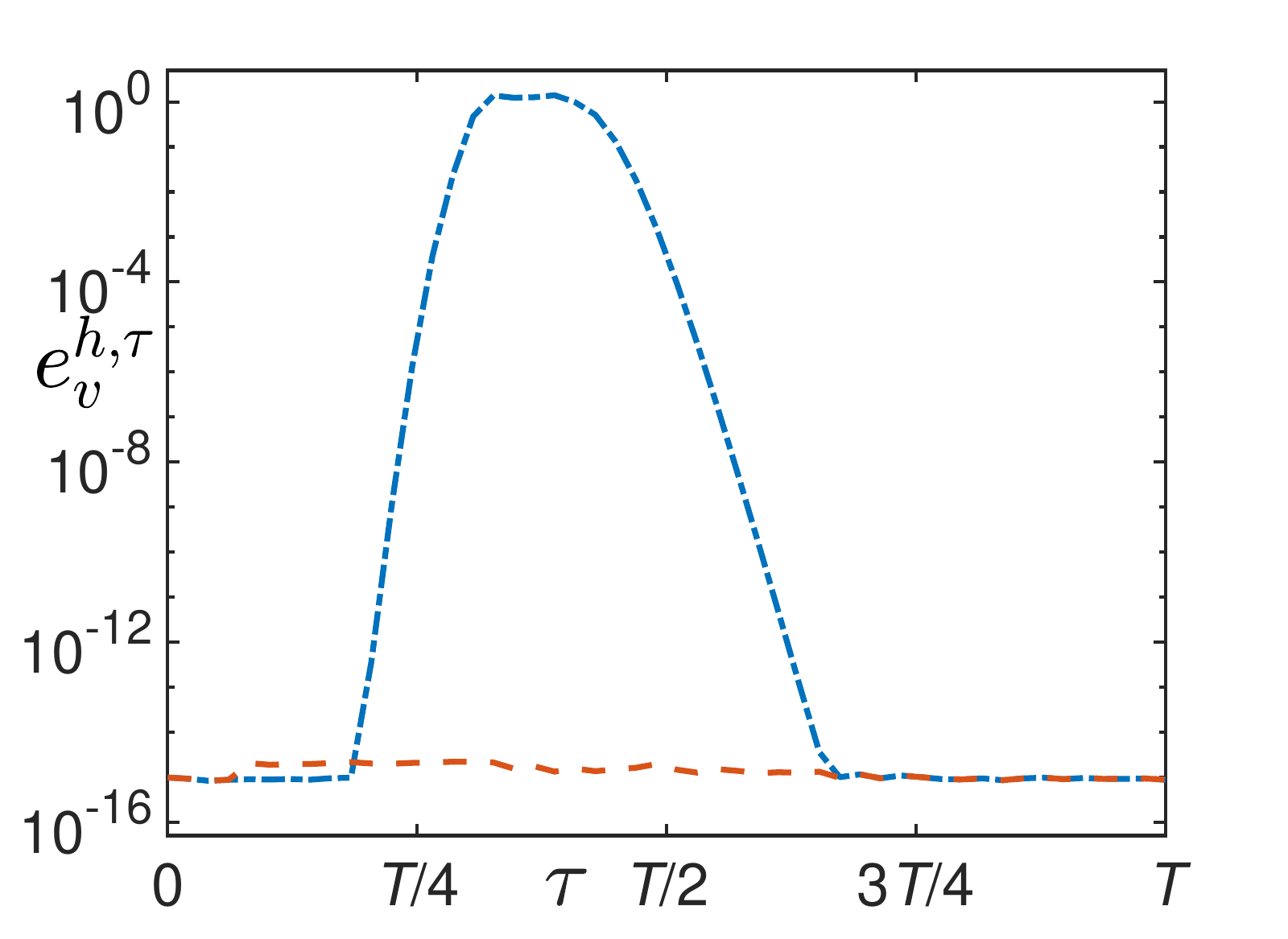}
  \caption{The errors of  \textbf{\textrm{Case \uppercase\expandafter{\romannumeral1}-\uppercase\expandafter{\romannumeral3}}} (left to right) in \textbf{Example \ref{Exmp:shear001}}. }
  \label{compare_error_001}
\end{figure}
Here we  take the computational domain as $[-32, 32]^3$ with mesh size $h = 1/8$, and the final time $t=\tau \in[0,T]$
where the period $T = 2\pi/|\bm{\Omega}|$.
Figure~\ref{compare_error_001} shows the errors over one period for three representative rotation vectors,
including cases where one or more components of $\bm{\Omega}$ are zero.
The results numerically confirm the effectiveness of our method (Algorithm \ref{LinearSol}).
%

\begin{exmp}
\label{Exmp:shear}
For the linear problem \eqref{LinearSubprob},  we consider the following parameters
\begin{itemize}[itemsep=2pt, topsep=8pt]
\item[] $\mathbf{Case ~\uppercase\expandafter{\romannumeral1}}$:~ $\bm{\Omega} = (2,6,8)^\top$ with period
$T=\frac{2\pi}{|\bm{\Omega}|}=\frac{\pi}{\sqrt{26}}\approx 0.62$;
\item[] $\mathbf{Case ~\uppercase\expandafter{\romannumeral2}}$:~ $\bm{\Omega} = (10,15,23)^\top$  with period $T=\frac{2\pi }{\sqrt{854}     } \approx 0.22$;
\item[] $\mathbf{Case ~\uppercase\expandafter{\romannumeral3}}$:~ $\bm{\Omega} = (20,40,40)^\top$  with period $T=\frac{\pi}{30} \approx      0.10$.
\end{itemize}
\end{exmp}

\begin{figure}[h]
    \centering
  \includegraphics[width=6cm, height=0.6cm]{./figure/four_shear/legend.pdf}\\
  \centering
  \includegraphics[width=4.4cm, height=3.3cm]{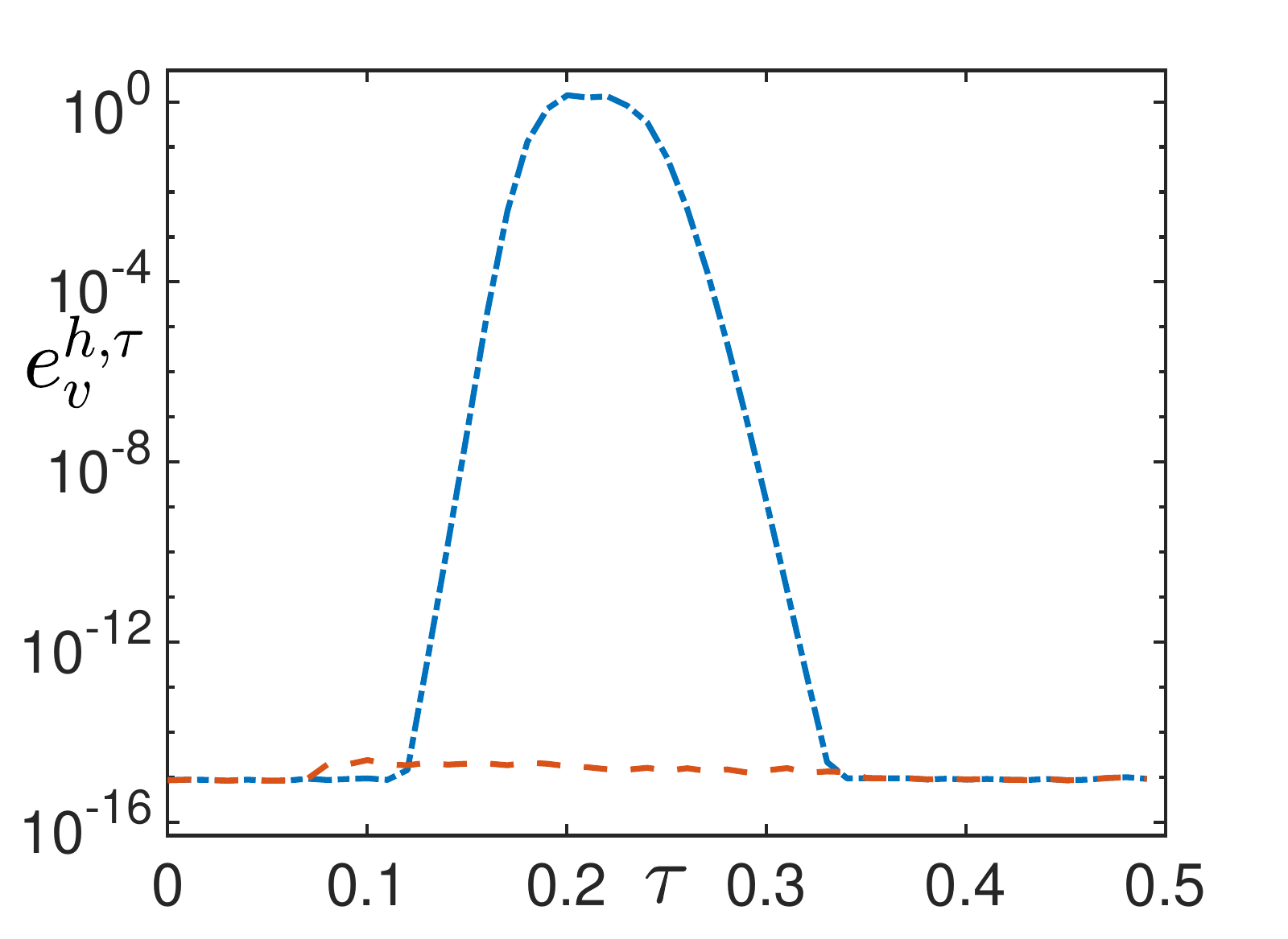}
  \hspace{-0.4cm}
    \includegraphics[width=4.4cm, height=3.3cm]{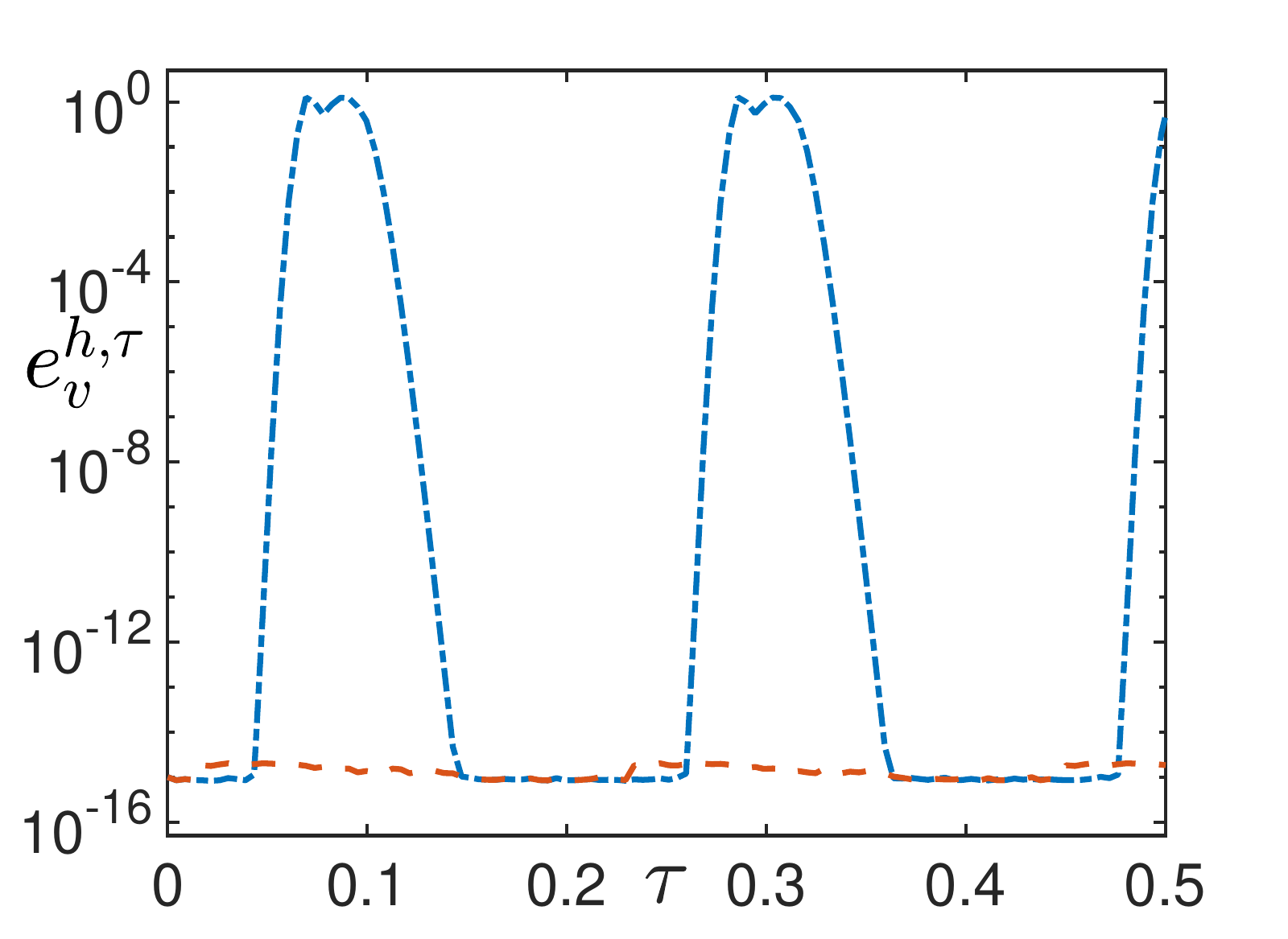}
  \hspace{-0.4cm}
  \includegraphics[width=4.4cm, height=3.3cm]{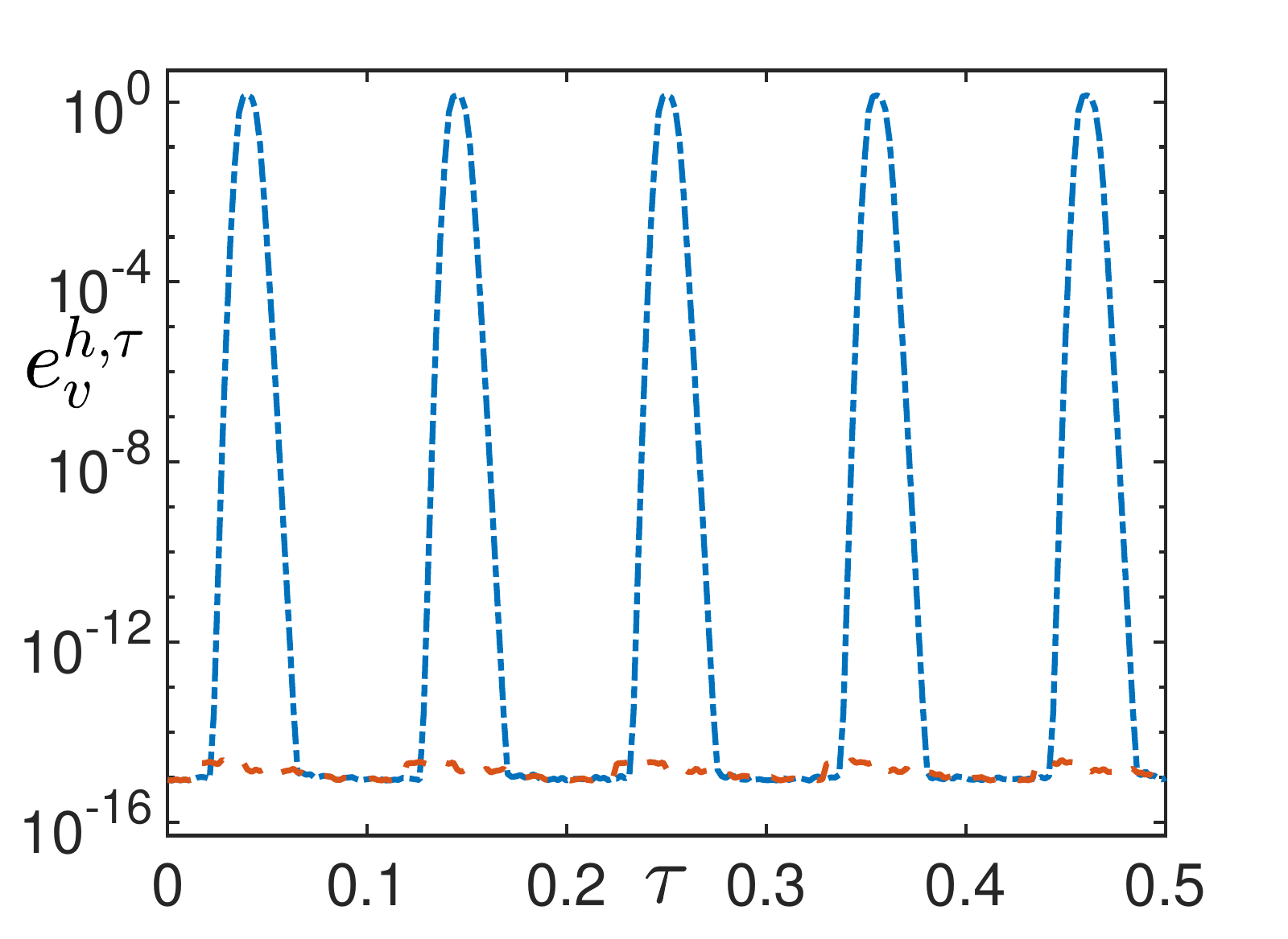}
  \caption{The errors of  \textbf{\textrm{Case \uppercase\expandafter{\romannumeral1}-\uppercase\expandafter{\romannumeral3}}} (left to right) in \textbf{Example \ref{Exmp:shear}}. }
  \label{compare_error}
\end{figure}

In simulations, the time step is typically not too large (e.g., below $1$), so we take $\tau \in (0, 0.5]$ here.
Within this interval, the above three cases correspond to approximately  $1$, $2$ and $5$ complete periods, respectively.
The computational domain is chosen as $[-32, 32]^3$ with $h = 1/8$ and the final time is $t=\tau$.
Figure \ref{compare_error} presents the errors for $\mathbf{Case ~\uppercase\expandafter{\romannumeral1}}$-$\mathbf{\uppercase\expandafter{\romannumeral3}}$,
from which we conclude that Algorithm \ref{LinearSol} can be applied with any time step.

\subsection{Exact integration of the nonlinear subproblem}
The nonlinear subproblem \eqref{NonLinearSubprob} can be solved exactly in physical space.
Noting  that $f_{\alpha}~ (\alpha=x,y,z)$ and $\mathcal{H}_{\rm{non}}$ are Hermitian matrices and satisfy the commutator relations
$[f_x,f_y] := f_x f_y -f_y f_x =\im f_z$, $[f_y,f_z]=\im f_x$ and $[f_z,f_x]=\im f_y$,
we obtain
\beas
\partial_t {F}_{\alpha}(t)
& =&\partial_t(\Psi^{\mathsf{H}}f_{\alpha} \Psi)
=\im \Psi^{\mathsf{H}} \mathcal{H}_{\rm{non}}(t) f_{\alpha} \Psi  - \im \Psi^{\mathsf{H}} f_{\alpha} \mathcal{H}_{\rm{non}}(t)  \Psi=  \im c_1 \Psi^{\mathsf{H}} [\mathbf{F}\cdot \mathbf{f}, f_{\alpha} ] \Psi\\[0.1em]
&=&  \im c_1 \big({F}_x \Psi^{\mathsf{H}} [f_x, f_{\alpha}] \Psi + {F}_y \Psi^{\mathsf{H}} [f_y, f_{\alpha}] \Psi
+ {F}_z \Psi^{\mathsf{H}} [f_z, f_{\alpha}] \Psi \big) =0,
\eeas
which implies that
${F}_{\alpha}(t) \!\equiv\! {F}_{\alpha}(t_n)
\!:=\! {F}^n_{\alpha}$.
In addition, it is easy to prove that the density $\rho(t)$ is time-independent, i.e.,
$\rho(t) \equiv \rho(t_n):=\rho^n$.
Thus,  the nonlinear subproblem \eqref{NonLinearSubprob} is equivalent to
\beas
\left\{\begin{array}{l}
i\partial_t\Psi(\bx,t)
=\left( \big[ V +c_0\rho^n \big]  I_3 +c_1 \mathbf{F}^n\cdot \mathbf{f}\right)\Psi(\bx,t) \quad t_n \le t \le t_{n+1}, \\[0.5em]
\Psi(\bx, t_n) =\Psi^{n}, \quad \bx \in \mathbb{R}^3, \\
\end{array}\right.
\eeas
where $\mathbf{F}^n =({F}_x^n, {F}_y^n,{F}_z^n)^\top$.
It is clear that the exact solution to \eqref{NonLinearSubprob} reads as
\be \label{NonSolMatr}
\Psi(\bx, t)  = e^{-\im  (V +c_0 \rho^n  ) (t-t_n)}   e^{-\im c_1 (t-t_n) \mathbf{F}^n\cdot \mathbf{f}} \Psi^n.
\ee
Fortunately, $\Psi^n$ is an eigenvector corresponding to the eigenvalue
$|\mathbf{F}^{n}|^2$ of  $\left(\mathbf{F}^n\cdot \mathbf{f}\right)^2$, i.e.,
$$
(\mathbf{F}^n\cdot \mathbf{f})^2 \Psi^n
= |\mathbf{F}^n|^2 \Psi^n.
$$
Then the solution \eqref{NonSolMatr} can be reformulated as
\begin{small}
\be\label{NonLinExtSolu}
\Psi(\bx, t) = e^{-\im   (V +c_0 \rho^n) (t-t_n)} \Big[\cos\big(c_1 (t\!-\!t_n) |\mathbf{F}^{n}|\big) \Psi^n
- \im \frac{\sin\big(c_1(t\!-\!t_n) |\mathbf{F}^{n}|\big)}{|\mathbf{F}^{n}|}  (\mathbf{F}^n\cdot \mathbf{f}) \Psi^n \Big].
\ee
\end{small}

\
As both subproblems of the spin-1 BECs with ArbRot term are exactly solvable in time, a numerical scheme can be easily constructed using a standard splitting procedure, such as the classical second-order Strang splitting method or fourth-order Yoshida splitting method.
To sum up, from time $t_n$ to $t_{n+1}$,
we can solve the CGPEs \eqref{GGPE} by using  Algorithm \ref{alg:dynamics}.

\begin{algorithm}
\caption{High-order splitting Fourier spectral methods}
\label{alg:dynamics}
\begin{algorithmic}[0]
\Require  Given $\Psi^n$ at time $t_n$.
\vspace{0.5em}
\State Compute the dynamics of \eqref{GGPE} via splitting scheme \eqref{splitting}:
\vspace{0.1cm}
\State \hspace{0.5cm}\textbf{1.} Solve the linear subproblem \eqref{LinearSubprob} by Algorithm \ref{LinearSol}.
\vspace{0.1cm}
\State \hspace{0.5cm}\textbf{2.} Solve the nonlinear subproblem \eqref{NonLinearSubprob} via \eqref{NonLinExtSolu}.
\vspace{0.5em}
\Ensure The updated wave function $\Psi^{n+1}$ at time $t_{n+1}$.
\end{algorithmic}
\end{algorithm}

\begin{remark}
In cases where real-time dynamics are not required,
we can first introduce the rotating coordinate transformation $\tilde{\bx}=\mathcal{R}(t)\bx$ into \eqref{GGPE}, which yields the transformed equation (without ArbRot term)
\be
\label{GGPE_noRot}
i\partial_t \Phi( \tilde{\bx},t) = \Big[\Big(-\frac{1}{2}\Delta  + W( \tilde{\bx},t) + c_0 \tilde{\rho}\Big)  I_{3}+c_1 \tilde{\mathbf{F}}\cdot \mathbf{f} \Big]\Phi(\tilde{\bx},t),
\ee
where $W( \tilde{\bx},t)=V(\mathcal{R}^{-1}(t)\tilde{\textbf{x}})$ is time-dependent,
and $\tilde{\rho},~ \tilde{\mathbf{F}}$ are obtained from $\rho,~ \mathbf{F}$ by replacing $\psi$ with $\phi$.
If $V(\bx)$ is isotropic, it is easy to obtain that $W( \tilde{\bx},t)=V(\tilde{\bx})$ is time-independent.
An analogous splitting scheme can be applied to Eq. \eqref{GGPE_noRot}.
It is worth pointing out that the computationally nontrivial rotating coordinate transformation is required only once in practice, at the final time.

\end{remark}

\subsection{Stability analysis}
\label{sec:stability}
To establish the practicality of our method, we present a rigorous stability analysis  in this subsection.
As shown in Theorem \ref{thm:Stabi},
Algorithm \ref{alg:dynamics} is unconditionally stable.
\begin{thm}
	\label{thm:Stabi}
The mass  is conserved at the discrete level,
i.e.,
for any mesh size $h>0$ and time step size $\tau >0 $,
we have
\beas
\sum_{v=-1}^{1} \| \psi_{v}^n \|^2_2 = \sum_{v=-1}^{1} \| \psi_{v}^0 \|^2_2.
\eeas
Thus we conclude that the splitting Fourier spectral method is unconditionally stable.
Moreover, the  magnetization is also conserved at the discrete level, i.e.,
\beas
\sum_{v=-1}^{1} v \| \psi_{v}^n \|^2_2 = \sum_{v=-1}^{1} v \| \psi_{v}^0 \|^2_2.
\eeas
\end{thm}
\begin{proof}
Considering Eq. \eqref{SolLinear}, we can obtain the following identities
\beas
\| \phi_v^{n+1} \|_{2}^2 = h^3 \sum_{\textbf{j}\in \mathcal{I}_N}|\phi^{n+1}_{v,~\!\textbf{j}}|^2 = \sum_{\textbf{p} \in \mathcal{K}_N}
 |(\widehat{\psi}_{v}^{n})_{\textbf{p}}|^2.
\eeas
Combined with Parseval's identity, we know
\bea
\| \phi_v^{n+1} \|_{2}^2=\| \psi_v^{n} \|_{2}^2.
\eea
Then for the rotating coordinate transformation, we get
\bea
\| \phi_v^{(1)} \|_{2}^2
&=& h^3 \sum_{\textbf{j}\in \mathcal{I}_N} |\phi_{v, ~\!\textbf{j}}^{(1)}|^2
 = h^3 \sum_{\textbf{j}\in \mathcal{I}_N} \left| \sum_{p\in \mathcal{J}} e^{\im (\lambda_1 y_k+\lambda_2 z_l) \nu_p^x} ~ \Big(\widehat{\phi_{v, kl}^{n+1}}\Big)^x_{p}  ~ e^{\im \nu_p^x x_j} \right|^2 \notag \\
&=& h^2 \sum_{k=0}^{N-1} \sum_{l=0}^{N-1} \sum_{p\in \mathcal{J}} \left|\Big(\widehat{\phi_{v, kl}^{n+1}}\Big)^x_{p} \right|^2
= h^3 \sum_{\textbf{j}\in \mathcal{I}_N} |\phi_{v,~\!\textbf{j}}^{n+1}|^2 = \| \phi_v^{n+1} \|_{2}^2.
\eea
In a similar way, we have the following equalities
\bea \label{DisL2Rot}
\| \psi_v^{n+1} \|_{2}^2 = \| \phi_v^{(3)} \|_{2}^2= \| \phi_v^{(2)} \|_{2}^2 =\| \phi_v^{(1)} \|_{2}^2.
\eea
Therefore, for the linear subproblem \eqref{LinearSubprob},
it satisfies the identity
 \beas
 \| \psi_v^{n+1} \|^2_2 =\| \psi_v^{n} \|^2_2.
\eeas

For the nonlinear subproblem \eqref{NonLinearSubprob}, we obtain
\beas
\| \Psi^{n+1} \|_2^2 =\|  e^{-\im  \big(c_0 \rho^n + V(\bx) \big) \tau}  e^{-\im c_1 \tau \mathbf{F}^n\cdot \mathbf{f}} \Psi^{n} \|_2^2 = \| \Psi^{n} \|_2^2,  \notag
\eeas
where $\| \Psi^{n} \|_{2}^2 := \sum_{v=-1}^{1} \| \psi_{v}^n \|^2_{2}$,
and
\beas
\sum_{v=-1}^{1} v \| \psi_{v}^{n+1} \|^2_2
&=& \|\cos(g^n) \psi_1^n - \im \frac{\sin(g^n)}{|\mathbf{F}^{n}|}  (F^n_z\psi_1^n + F^n_{-}\psi_0^n/\sqrt{2})\|_{2}^2
\\
&& -\|\cos(g^n) \psi_{-1}^n - \im \frac{\sin(g^n)}{|\mathbf{F}^{n}|}  (-F^n_z\psi_{-1}^n + F^n_{+}\psi_0^n/\sqrt{2})\|_{2}^2\\
&=& h^3 \!\sum_{\textbf{j}\in \mathcal{I}_N}
\!\big( \!  \cos^2(g^n)F^n_z
\!+ \!\sin^2(g^n)((F^n_z)^2\!+\!F^n_{+}F^n_{-})F^n_z/|\mathbf{F}^{n}|^2   \big)_{\textbf{j}}\\
&=& h^3 \sum_{\textbf{j}\in \mathcal{I}_N}
\big(   F^n_z  \big)_{\textbf{j}}
=  \sum_{v=-1}^{1} v \| \psi_{v}^n \|^2_2
\eeas
with $g^n := c_1 \tau |\mathbf{F}^{n}|$.
 The proof is completed.
\end{proof}

\subsection{Comparison with other methods}
Since no other method has ever been proposed for this specific problem, we adapt a splitting method based on Explicit-Exact-Integrator
(EEI),
which has been successfully applied to  scalar dipolar BECs with ArbRot \cite{EEI7}.
Specifically,
the linear subproblem is solved by
\bea
\label{ExtSplit3D_seven}
e^{\im\tau\left[ \Delta/2+\bm{\Omega} \cdot \bm{L} \right]} =
& e^{-\im \tau (\zeta_1 x^2+\zeta_2 y^2 + \zeta_3 z^2)}
e^{-\tau (\xi_1 y+\xi_2 z) \partial_{x}}  e^{-\tau \xi_3 z \partial_{y}} \\ \notag
& e^{\im \tau \nabla^{\top} K \nabla}
e^{-\tau \eta_1 x \partial_{y}} e^{-\tau (\eta_2 x +\eta_3 y)\partial_{z}} e^{-i\tau \bx^{\top}  P \bx},
\eea
where $\zeta_1$, $\zeta_2$, $\zeta_3$, $\xi_1$, $\xi_2$, $\xi_3$,
 $\eta_1, \eta_2, \eta_3 \in \mathbb R$
and $K$, $P \in \mathbb{R}^{3\times 3}$, with the coefficients given explicitly in \cite{EEI7}.
The nonlinear subproblem is solved by \eqref{NonLinExtSolu}.
For convenience, we denote the second-order/fourth-order EEI-based methods as \textbf{ECS2}/\textbf{ECS4},
and denote the proposed second-order/fourth-order splitting Fourier spectral method in this paper as $\textbf{CS2}$/$\textbf{CS4}$.
Since both the linear and nonlinear subproblems can be solved exactly using either method,
we will only compare their computational complexity.

\begin{table}[htb!]
    \centering
    \caption{Computational complexity of linear subproblem \eqref{LinearSubprob} per time step.}
    \label{comparsioncomplexity}
    \begin{tabular}{ccccc}
    \toprule
    &Order &Method  &FFT/iFFT pairs  & Multiplications \\
    \midrule
    &\multirow{2}{*}{Second}   & CS2 & $ 21N^2 $ &  $ 15N^3$ \\[0.3em]
    & & ECS2  & $ 21N^2$ &$ 21N^3$\\
    \midrule
    &\multirow{2}{*}{Fourth}  & CS4 & $ 63N^2 $ &$ 45N^3$ \\[0.3em]
    & & ECS4 &  $ 63N^2$ &$ 63N^3$ \\
    \bottomrule
    \end{tabular}
\end{table}
		
Table \ref{comparsioncomplexity} presents the computational complexity for the linear subproblem \eqref{LinearSubprob} of CS2/CS4 and ECS2/ECS4,
where ``FFT/iFFT pairs" refers to one-dimensional FFT/iFFT pairs.
From this table, we can see that the computational complexity of ECS2/ECS4 is slightly larger than that of our methods.

\section{Numerical results}
\label{sec:NumeResu}
In this section, we first  verify the spatial and temporal  accuracy as well as efficiency for the wave function evolution. Subsequently, we explore the dynamical laws,
including the evolution of mass, energy, magnetization and center of mass.
Finally, we investigate some interesting phenomena,
such as the dynamics evolution of spin-1 BECs under different physical parameters.
In the following examples, we consider the harmonic potential $V(\bx)$ \eqref{harmonicPot} and
take the computational domain $[-8, 8]^3$ unless otherwise specified.
The algorithms were implemented in FORTRAN
and run on a single 3.00 GHz Intel(R) Xeon(R) Gold 6248R CPU with a 36 MB cache in Ubuntu GNU/Linux with the Intel compiler ifort.

\subsection{Performance on accuracy and efficiency}
The reference solution at time $t$ is computed by our method using a very fine mesh size $h_0 = 1/16$ and a small time step $\tau_0 = 10^{-4}$.
Spatial errors are calculated using a fixed
time step $\tau_0$ with varying mesh sizes $h$, and temporal errors are computed using a fixed mesh size $h_0$ with  varying time step $\tau$.

\begin{exmp}[\textbf{Accuracy}]
\label{Exmp:Accu}
Here, we test the spatial and temporal accuracy.
To this end, we choose $\bm{\Omega}=(0.1, 0.2, 0.3)^{\top}$,  $c_0 = 8$, $c_1 = 0.8$, $\bm{\gamma}:= (\gamma_x,~ \gamma_y,~ \gamma_z)^{\top} = (1,~ 1,~ 1)^{\top}$.
The initial value is chosen as
\begin{equation}\label{initial2DNum}
\Psi^0(\mathbf{x}) = \big(\psi_1^0(\mathbf{x}), \psi_0^0(\mathbf{x}), \psi_{-1}^0(\mathbf{x})\big)^\top  =\phi(\mathbf{x})(2,\sqrt{7},3)^\top,
\end{equation}
where
$\phi(\mathbf{x})= \sqrt{0.05}/ \pi^{3/4} e^{-|\mathbf{x}|^2/2}$.
\end{exmp}

Table \ref{Accu3D} shows the temporal and spatial errors of
CS2 and CS4 at time $t=0.3$.
It can be observed that CS2/CS4 achieves spectral accuracy in space and
second/fourth order convergence in time.

\begin{table}[htb!]
    \centering
    \caption{Numerical errors at time $t = 0.3 $  in \textbf{Example \ref{Exmp:Accu}.}}\label{Accu3D}
    \begin{tabular}{cccccc}
    \toprule
    \multirow{13}{*}{CS2}
    & \multicolumn{5}{c} {\textbf{Temporal direction}} \\
    \cmidrule{2-6}
    &$\tau$& $1/80$  & $1/160$ & $1/320$ & $1/640$ \\
    \cmidrule{2-6}
& $e_{1}^{h_0,\tau}$    & 1.8553E-05 & 4.6367E-06 & 1.1582E-06 & 2.8867E-07 \\
& rate                  &            & 2.0005     & 2.0012     & 2.0045     \\
& $e_{0}^{h_0,\tau}$     & 1.9641E-05 & 4.9085E-06 & 1.2261E-06 & 3.0559E-07 \\
& rate                  &            & 2.0005     & 2.0012     & 2.0045     \\
& $e_{-1}^{h_0,\tau}$   & 1.8914E-05 & 4.7269E-06 & 1.1808E-06 & 2.9428E-07 \\
& rate                  &            & 2.0005     & 2.0012     & 2.0045     \\
    \cmidrule{2-6}
    &\multicolumn{5}{c} {\textbf{Spatial direction}} \\
    \cmidrule{2-6}
    &$h$ & $1$ &$1/2$& $1/4$ & $1/8$   \\
   \cmidrule{2-6}
& $e_{1}^{h,\tau_0}$   & 1.5155E-02 & 1.3491E-04 & 7.0937E-10 & 2.0294E-13  \\[0.1em]
& $e_{0}^{h,\tau_0}$   & 1.6069E-02 & 1.4765E-04 & 7.9857E-10 & 2.0283E-13  \\[0.1em]
& $e_{-1}^{h,\tau_0}$  & 1.5459E-02 & 1.3915E-04 & 7.3910E-10 & 2.0253E-13  \\[0.1em]
    \midrule
    \multirow{13}{*} {CS4}
    &\multicolumn{5}{c} {\textbf{Temporal direction}} \\
    \cmidrule{2-6}
    &$\tau$ & $1/80$& $1/160$ & $1/320$ & $1/640$   \\
    \cmidrule{2-6}
& $e_{1}^{h_0,\tau}$    & 4.4758E-08 & 2.8171E-09 & 1.7642E-10 & 1.1103E-11 \\
& rate                    &    & 3.9898     & 3.9972     & 3.9899     \\
& $e_{0}^{h_0,\tau}$    & 4.7503E-08 & 2.9899E-09 & 1.8724E-10 & 1.1781E-11 \\
& rate                  &        & 3.9899     & 3.9972     & 3.9903     \\
& $e_{-1}^{h_0,\tau}$  & 4.5673E-08 & 2.8747E-09 & 1.8003E-10 & 1.1332E-11 \\
& rate                  &           & 3.9898     & 3.9971     & 3.9897     \\
     \cmidrule{2-6}
    &\multicolumn{5}{c} {\textbf{Spatial direction}} \\
    \cmidrule{2-6}
    &$h$ & $1$&$1/2$& $1/4$ & $1/8$  \\
    \cmidrule{2-6}
& $e_{1}^{h,\tau_0}$   & 1.5155E-02 & 1.3491E-04 & 7.0937E-10 & 6.5142E-13  \\[0.1em]
& $e_{0}^{h,\tau_0}$   & 1.6069E-02 & 1.4765E-04 & 7.9857E-10 & 6.5290E-13  \\[0.1em]
& $e_{-1}^{h,\tau_0}$  & 1.5459E-02 & 1.3915E-04 & 7.3910E-10 & 6.5287E-13  \\[0.1em]
    \bottomrule
    \end{tabular}
\end{table}

\begin{exmp}[\textbf{Efficiency}]\label{Exmp:EffiTest}
We test the efficiency by showing the computational time variation versus the total grid number $N_{t}:=N^3$.
To this end, we choose $\bm{\Omega}=(0.6, 0.7, 0.8)^{\top}$, $c_0 = 50$, $c_1=-1$, $\bm{\gamma}= (1,~ 1,~ 1)^{\top}$ and the initial value \eqref{initial2DNum}.
\end{exmp}

\begin{figure}[htbp]
    \centering
    \begin{minipage}[t]{0.48\textwidth}
        \centering
        \vspace{-3.8cm}
        \captionof{table}{Computational time versus  $N_{t}$  in \textbf{Example \ref{Exmp:EffiTest}.} }
        \setlength\extrarowheight{1.5pt}
        \begin{tabular}{m{0.7cm} m{0.6cm} m{0.8cm} m{0.8cm} m{0.9cm}}
        \toprule
          $N_t$ & $64^3$ & $128^3$ & $192^3$ & $256^3~$ \\
        \midrule
         CS2 & 8.15 & 74.81 & 278.61 & 732.93~ \\
         CS4 & 19.57 & 188.45 & 683.86 & 1792.39~ \\
        \bottomrule
        \end{tabular}

        \label{tab:cputime_table}
    \end{minipage}
    \hfill
    \begin{minipage}[t]{0.48\textwidth}
    \centering
    \includegraphics[width=5.3cm,height=4.2cm]{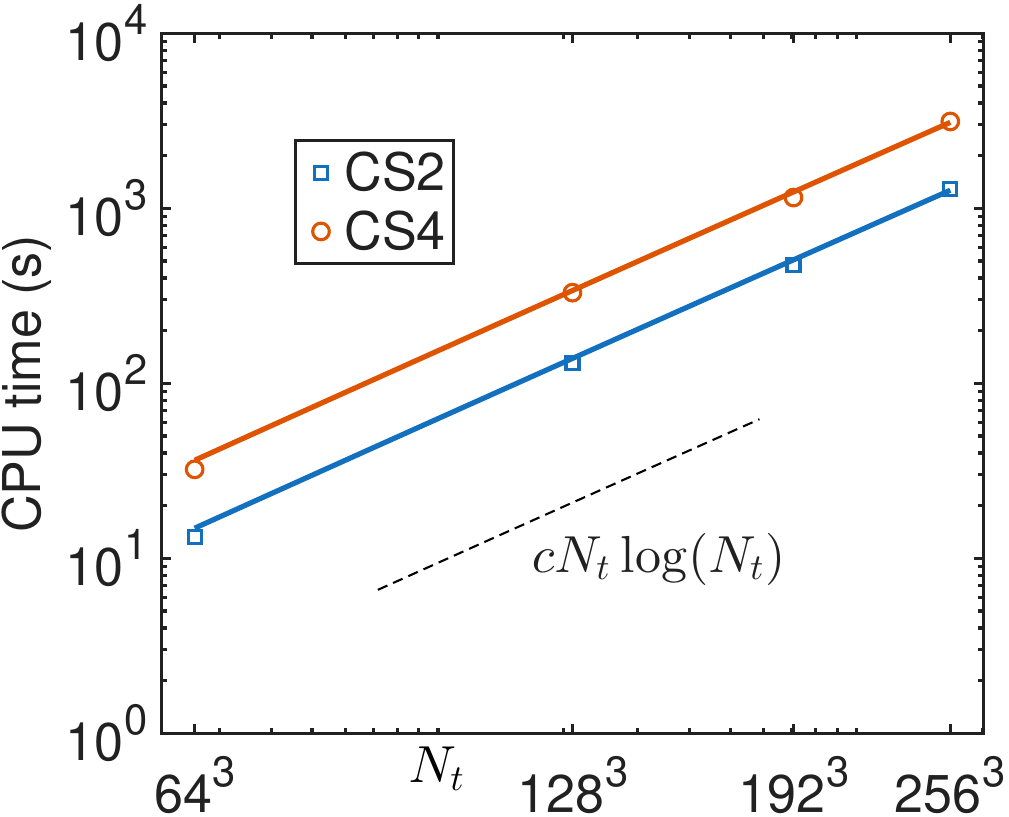}
    \vspace{-0.2cm}
    \caption{Log-log plot of computational time versus $N_{t}$  in \textbf{Example \ref{Exmp:EffiTest}.} }
    \label{effiTestloglogEEI}
    \end{minipage}
\end{figure}

Table \ref{tab:cputime_table} and Figure \ref{effiTestloglogEEI} present the computational time (in seconds) as a function of the total grid number $N_{t}$,
elapsed from $t=0$ to $t=0.1$ with time step $\tau = 10^{-3}$.
We can see that our methods are efficient, and the CPU time scales roughly as $C N_{t} \log(N_{t})$,
which agrees well with our expectation.

\subsection{Dynamics of physical quantities} \label{DynaProp}

Here, we investigate the dynamical laws,
including the evolution of mass, energy, magnetization and center of mass.

\begin{exmp}[\textbf{Dynamics of mass, energy and magnetization}]
\label{Exmp:DynaPropAngular}
For Eq. \eqref{GGPE}, the initial value is given by
\begin{equation}\label{initialDyn}
\Psi^0(\mathbf{x}) = \big(\psi_1^0(\mathbf{x}), \psi_0^0(\mathbf{x}), \psi_{-1}^0(\mathbf{x})\big)^\top  = \phi(\mathbf{x})(3,2,\sqrt{7})^\top,
\end{equation}
where
$\phi(\mathbf{x})=\sqrt{0.05}/ \pi^{3/4} e^{-|\mathbf{x}|^2/2}$.
We then consider the following cases:
\begin{itemize}[itemsep=4pt, topsep=4pt]
\item[] $\mathbf{Case ~ \uppercase\expandafter{\romannumeral1}}:$  $\bm{\Omega} =(0.1,~ 0.3, ~0.5)^{\top}$,  $\bm{\gamma} = (1,~ 1,~ 1)^     {\top}$, $c_0 = 80$, $c_1 = -1$;
\item[] $\mathbf{Case ~ \uppercase\expandafter{\romannumeral2}}:$  $\bm{\Omega} =(0.2,~ 0.6, ~1.0)^{\top}$,  $\bm{\gamma} = (1.005,~ 1,~ 1)^{     \top}$, $c_0 = 8$, $c_1 =0.8$;
\item[] $\mathbf{Case ~ \uppercase\expandafter{\romannumeral3}}:$  $\bm{\Omega} =(0.1,~ 0.3, ~0.5)^{\top}$,  $\bm{\gamma} = (1,~ 1.005,~ 1)     ^{\top}$, $c_0 = 100$, $c_1 = 5$;
\item[] $\mathbf{Case ~ \uppercase\expandafter{\romannumeral4}}:$  $\bm{\Omega} =(0.2,~ 0.6, ~1.0)^{\top}$,  $\bm{\gamma} = (1,~ 1,~ 1.005)     ^{\top}$, $c_0 = 20$, $c_1 = -1$.
\end{itemize}
\end{exmp}

\begin{figure}[h!]
 \centering
 \includegraphics[width=6.3cm, height=5cm]{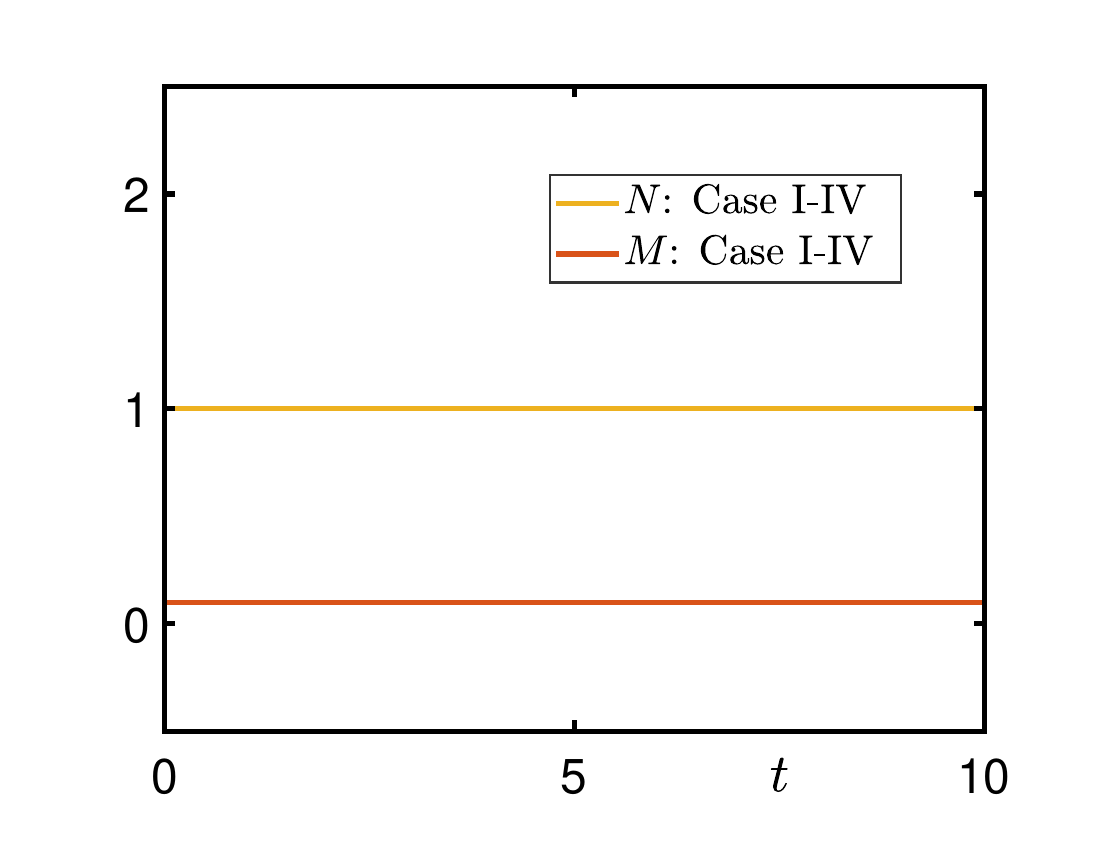}
 \hspace{0.1cm}
 \includegraphics[width=6.3cm, height=5cm]{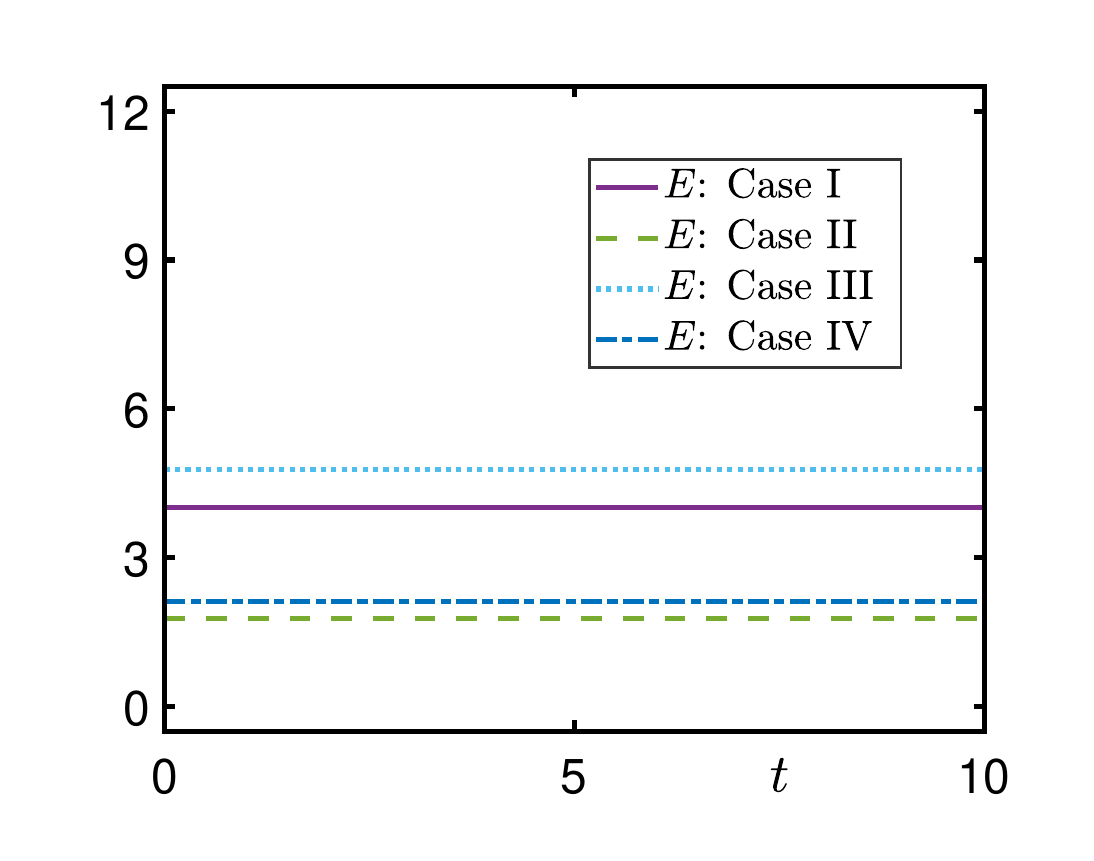}
 \\[0em]
 \caption{Evolution of energy ${E}(t)$, mass ${N}(t)$ and magnetization $M(t)$ for  \textbf{\textrm{Case \uppercase\expandafter{\romannumeral1}-\uppercase\expandafter{\romannumeral4}}} in   \textbf{Example \ref{Exmp:DynaPropAngular}}.}
 \label{MassEnergyMag}
\end{figure}
\vspace{-0.2cm}
In Figure \ref{MassEnergyMag}, we show the dynamics of energy ${E}(t)$, mass ${N}(t)$ and magnetization ${M}(t)$
obtained numerically by CS4 with $h=1/8$ and $\tau = 10^{-3}$ over the interval $t \in [0, 10]$.
Our method conserves the mass and magnetization very well during the simulation, and
these numerical observations align well with the analytical results (Theorem \ref{thm:Stabi}).

\begin{exmp}[\textbf{Dynamics of center of mass}]
\label{Exmp:DynaPropCenter}
The initial value is taken as
\begin{equation}
\Psi^0(\mathbf{x}) = \big(\psi_1^0(\mathbf{x}), \psi_0^0(\mathbf{x}), \psi_{-1}^0(\mathbf{x})\big)^\top  =\phi(\mathbf{x}-\mathbf{x}_0)(\sqrt{6},2\sqrt{2},\sqrt{6})^\top,
\end{equation}
where
$\phi(\mathbf{x})= \sqrt{0.05}/ \pi^{3/4} e^{-|\mathbf{x}|^2/2}$ and $\mathbf{x}_0 = (1.2,0.8,0.2)^\top $.
We take $c_0=100$, $c_1 = -1$
and study the following three cases:
\begin{itemize}[itemsep=4pt, topsep=4pt]
  \item[] $\mathbf{Case ~  \uppercase\expandafter{\romannumeral1}}:$ $\bm{\gamma} = (1, 1, 1)^{\top}$, $\bm{\Omega} =(0, 0, 0)^{\top},  (0, 0, 1)^{\top}, (0, 1.2, 0)^{\top}$;  
\item[] $\mathbf{Case ~  \uppercase\expandafter{\romannumeral2}}:$ $\bm{\gamma} = (1, 1, 1)^{\top}$, $\bm{\Omega} = (\frac{ \pi}{6},\frac{ \pi}{6}, \frac{ \pi}{6})^{\top}, (\frac{ \pi}{2},\frac{ \pi}{6}, \frac{ \pi}{6})^{\top}, (\frac{ \pi}{6}, \frac{ \pi}{2}, \frac{ \pi}{6})^{\top}$; 
  \item[] $\mathbf{Case ~  \uppercase\expandafter{\romannumeral3}}:$ $\bm{\Omega} = (1,1,1)^{\top}$, $\bm{\gamma} =(1, 1, 1)^{\top},  (0.6, 1, 1)^{\top}, (1, 0.6, 1)^{\top}$.
\end{itemize}
\end{exmp}

\begin{figure}[!htpb]
 \centering
 \includegraphics[width=4.2cm, height=4.2cm]{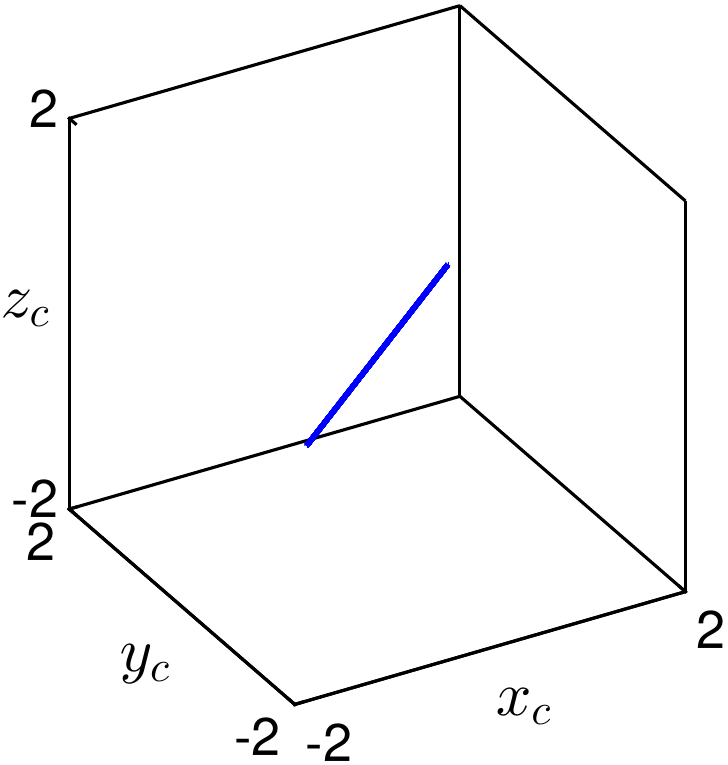}
 \hspace{-0.1cm}
  \includegraphics[width=4.2cm, height=4.2cm]{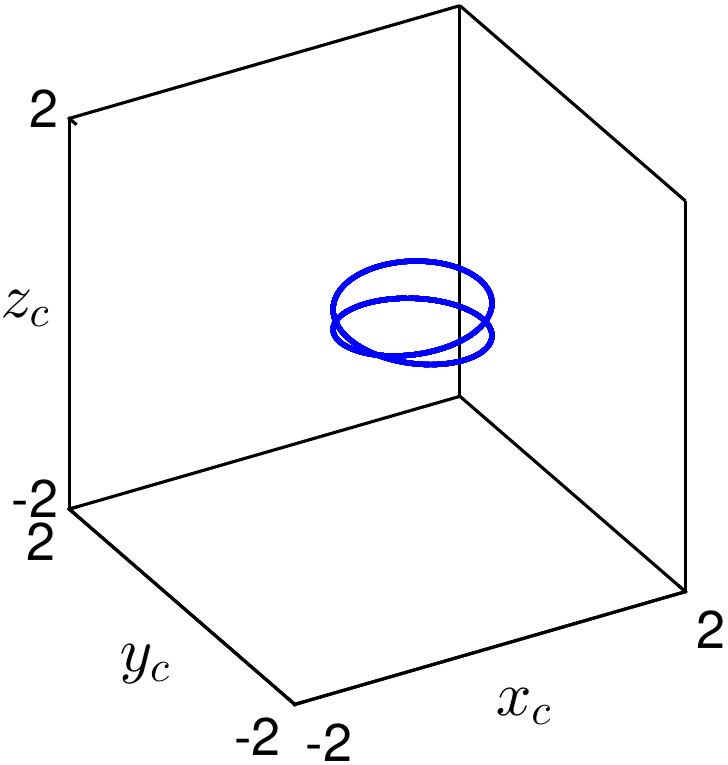}
 \hspace{-0.1cm}
 \includegraphics[width=4.2cm, height=4.2cm]{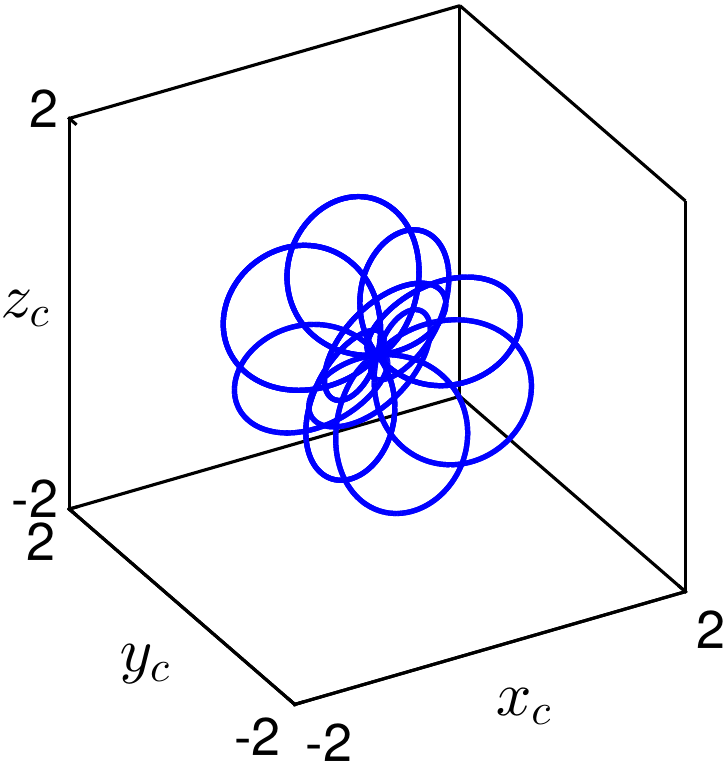}
 \\[0.5em]
  \includegraphics[width=4.2cm, height=4.2cm]{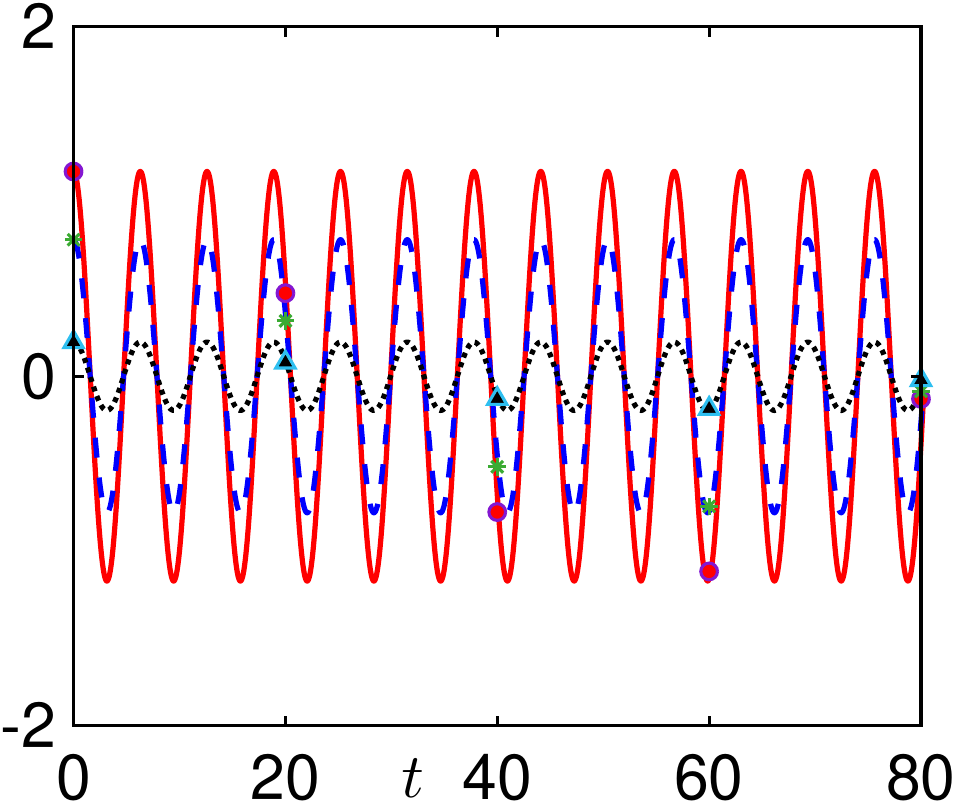}
 \hspace{-0.1cm}
  \includegraphics[width=4.2cm, height=4.2cm]{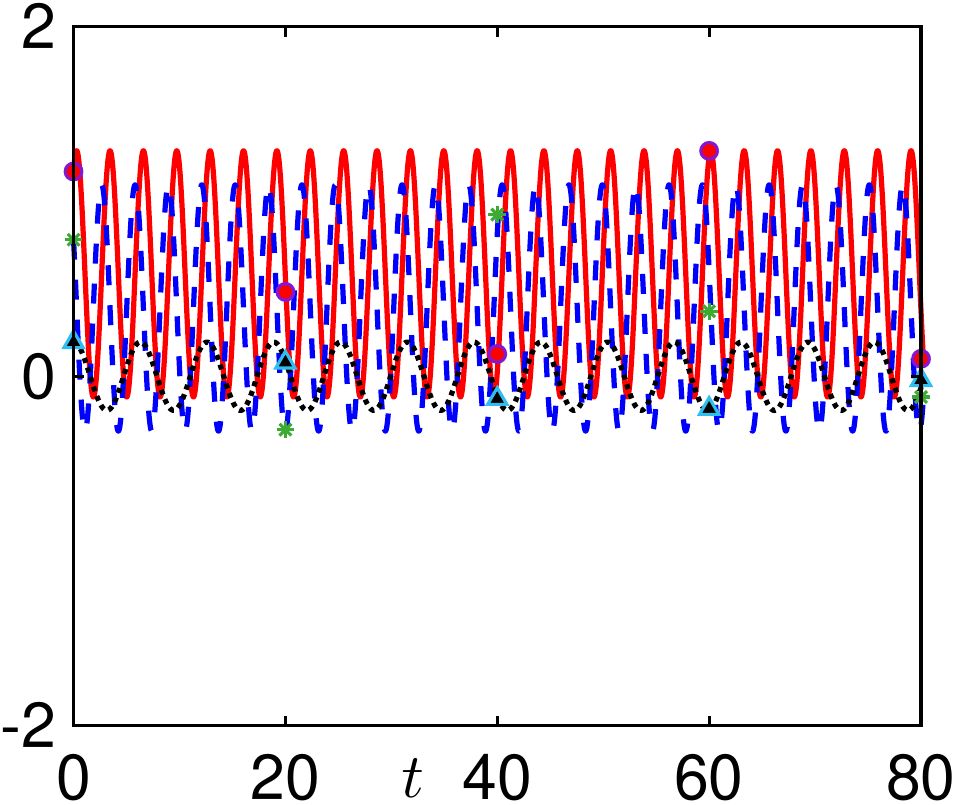}
 \hspace{-0.1cm}
 \includegraphics[width=4.2cm, height=4.2cm]{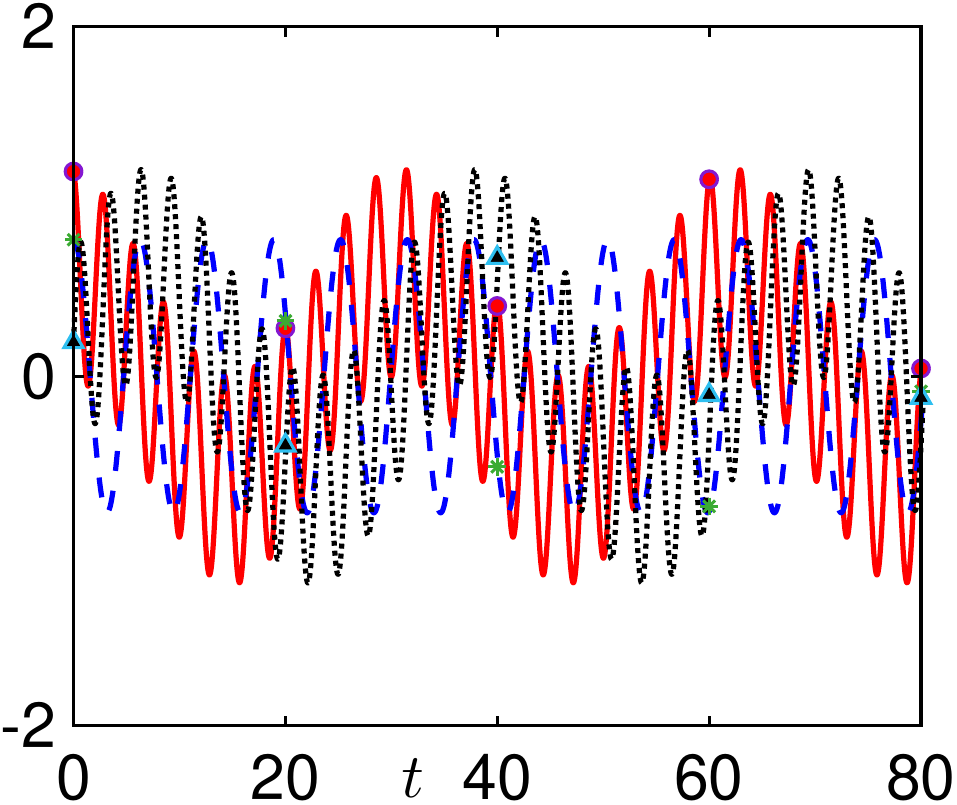}
 \\[0.5em]
 \includegraphics[width=4.2cm, height=4.2cm]{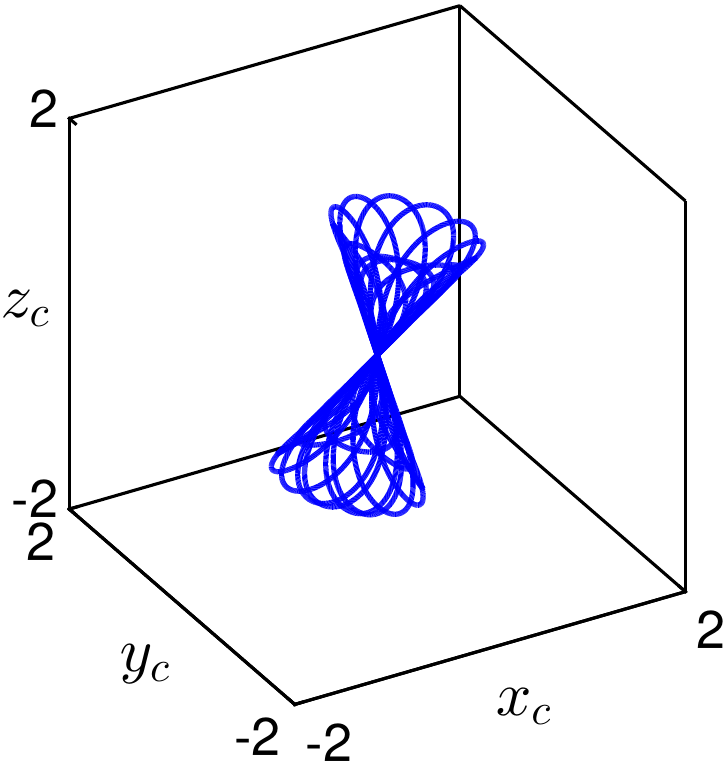}
 \hspace{-0.1cm}
  \includegraphics[width=4.2cm, height=4.2cm]{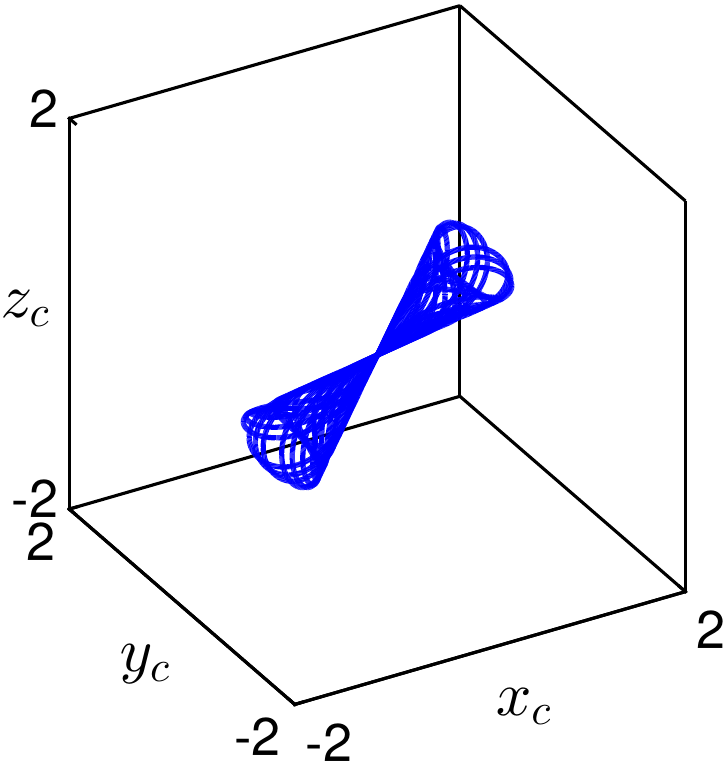}
 \hspace{-0.1cm}
   \includegraphics[width=4.2cm, height=4.2cm]{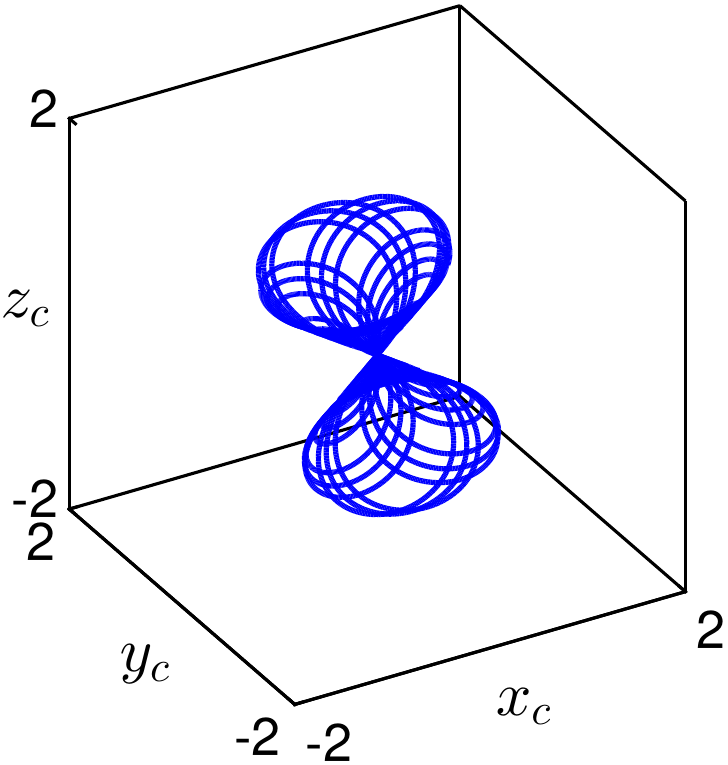}
 \\[0.5em]
  \includegraphics[width=4.2cm, height=4.2cm]{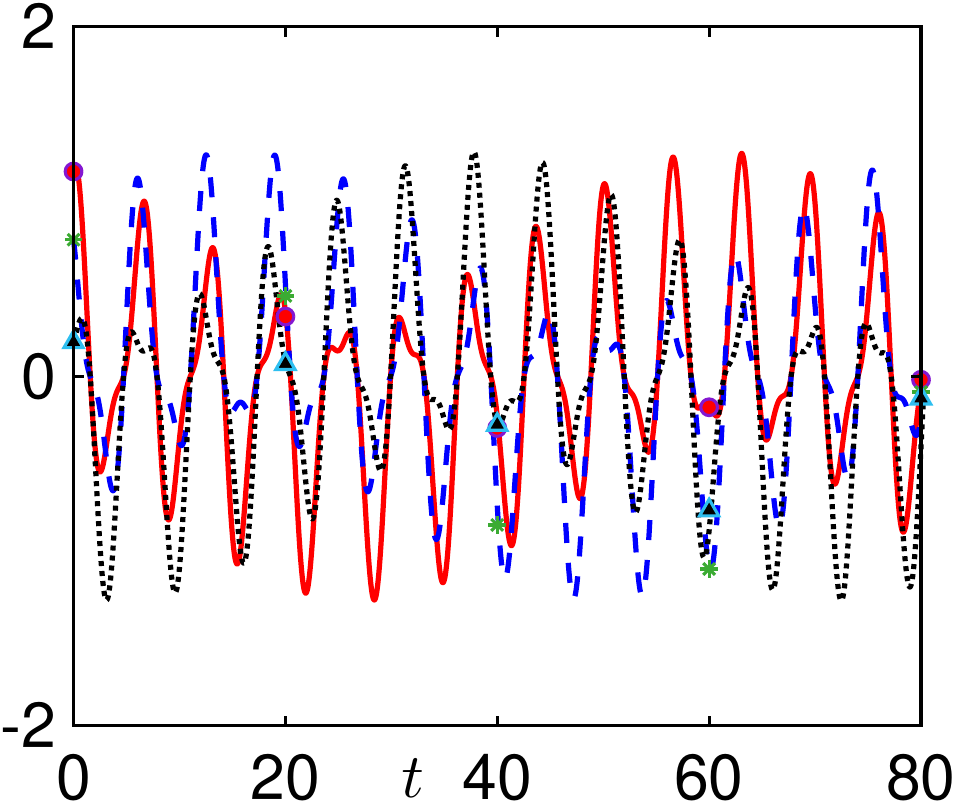}
 \hspace{-0.1cm}
  \includegraphics[width=4.2cm, height=4.2cm]{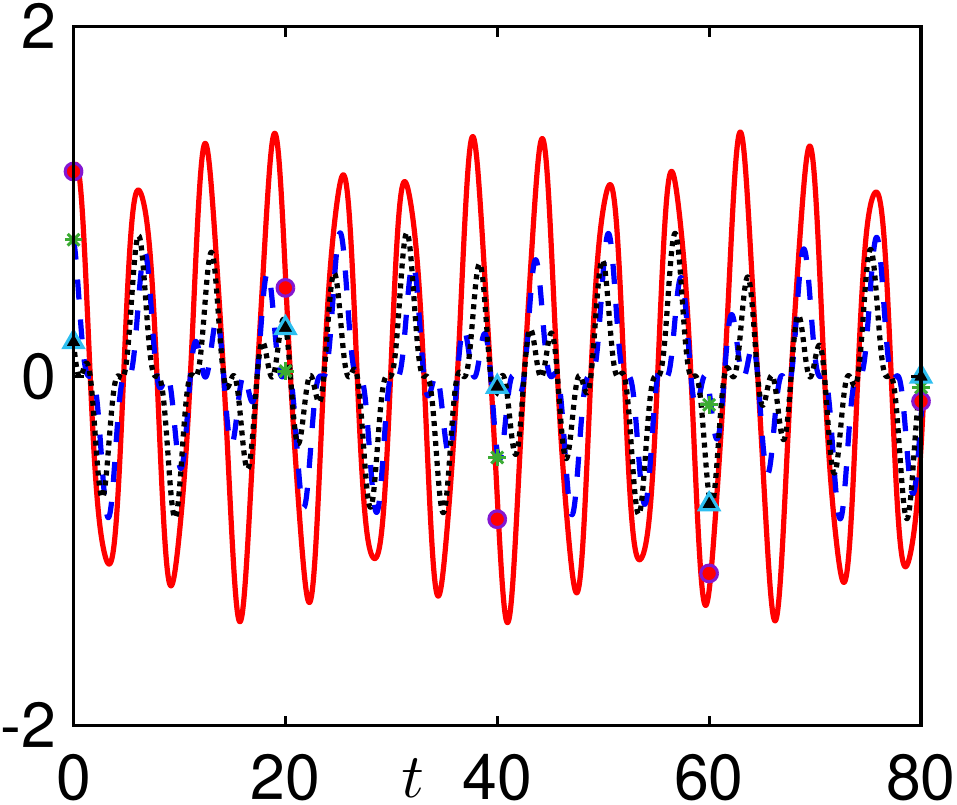}
 \hspace{-0.1cm}
 \includegraphics[width=4.2cm, height=4.2cm]{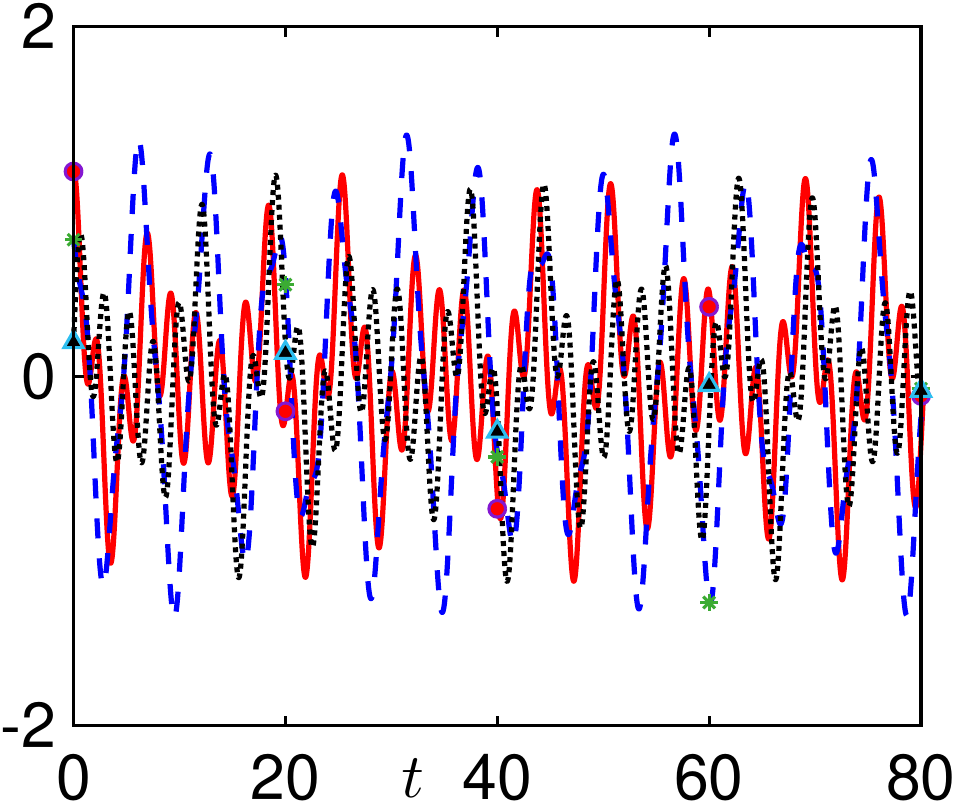}
 \\
 \caption{Motion of the center of mass $\bx_c(t) = (x_c(t), y_c(t), z_c(t))^{\top}$  for \textbf{\textrm{Case \uppercase\expandafter{\romannumeral1}}} (top two rows) and \textbf{\textrm{Case \uppercase\expandafter{\romannumeral2}}} (bottom two rows)  in \textbf{Example \ref{Exmp:DynaPropCenter}}.
Rows 1 and 3: trajectory $\bx_c(t)$ for $t \in [0, 80]$.
Rows 2 and 4: time evolution of $x_c(t)$ (red solid line),
$y_c(t)$ (blue dot line) and $z_c(t)$ (black  dash line) that are obtained by directly simulating the CGPEs,
where 	
circle($o$), 	
asterisk($\ast$) and triangle($\triangle$) represent the solutions to  ODEs \eqref{centerODE}-\eqref{centerODEInitial}.
The left-to-right order corresponds to the three parameters provided.}
 \label{MassCenterOmega}
\end{figure}

\begin{figure}[!htpb]
 \centering
 \includegraphics[width=4.2cm, height=4.2cm]{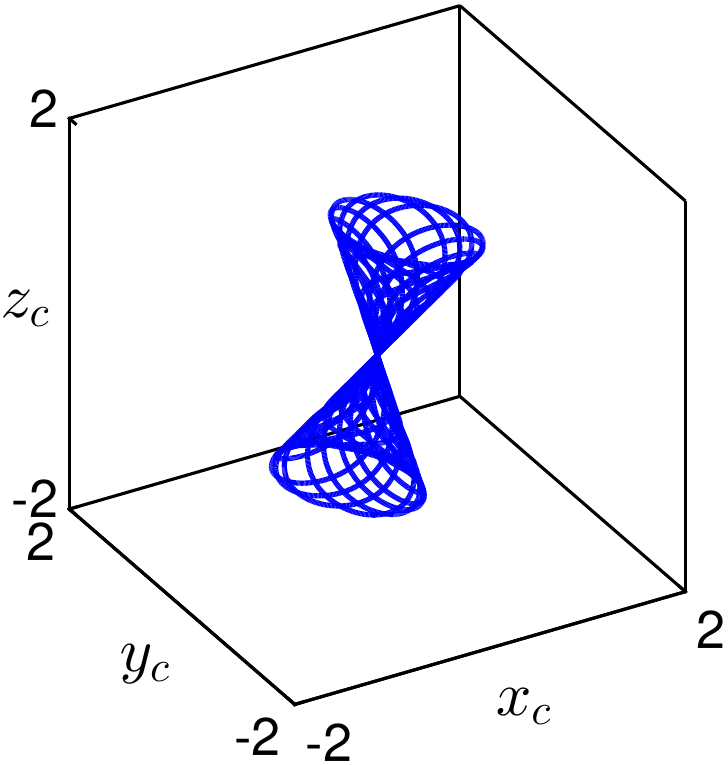}
 \hspace{-0.1cm}
  \includegraphics[width=4.2cm, height=4.2cm]{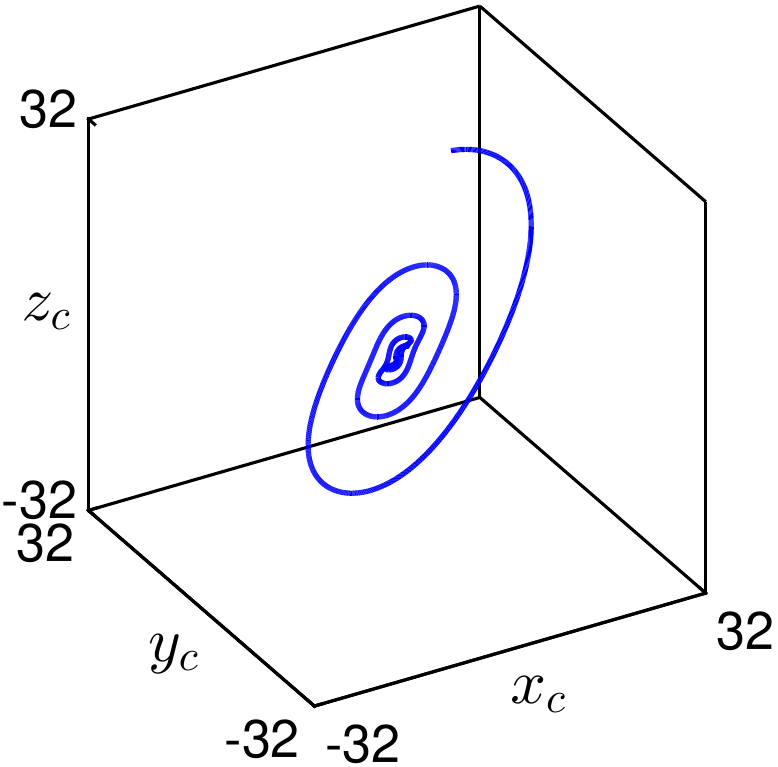}
 \hspace{-0.1cm}
   \includegraphics[width=4.2cm, height=4.2cm]{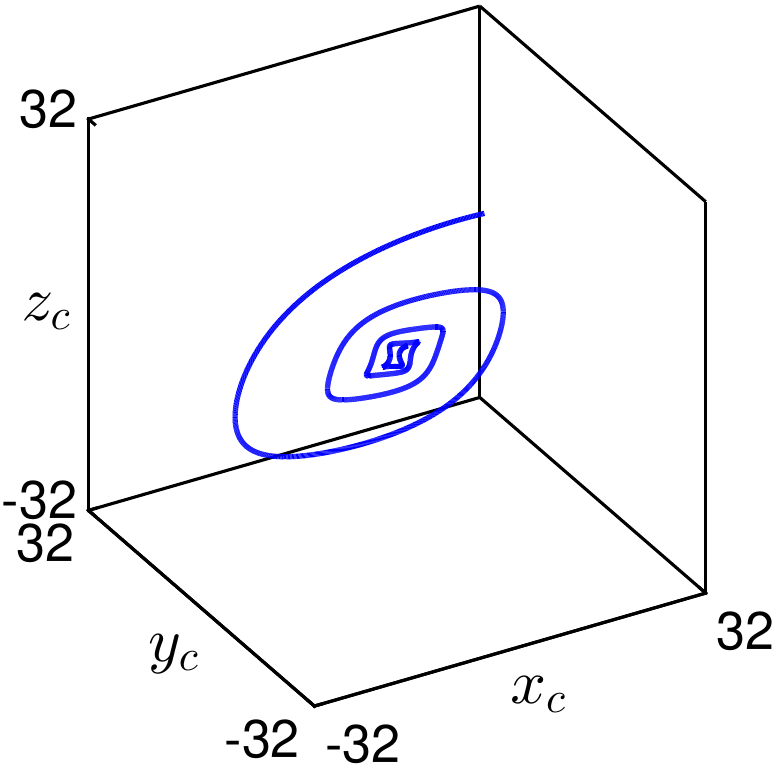}
 \\[0.5em]
  \includegraphics[width=4.2cm, height=4.2cm]{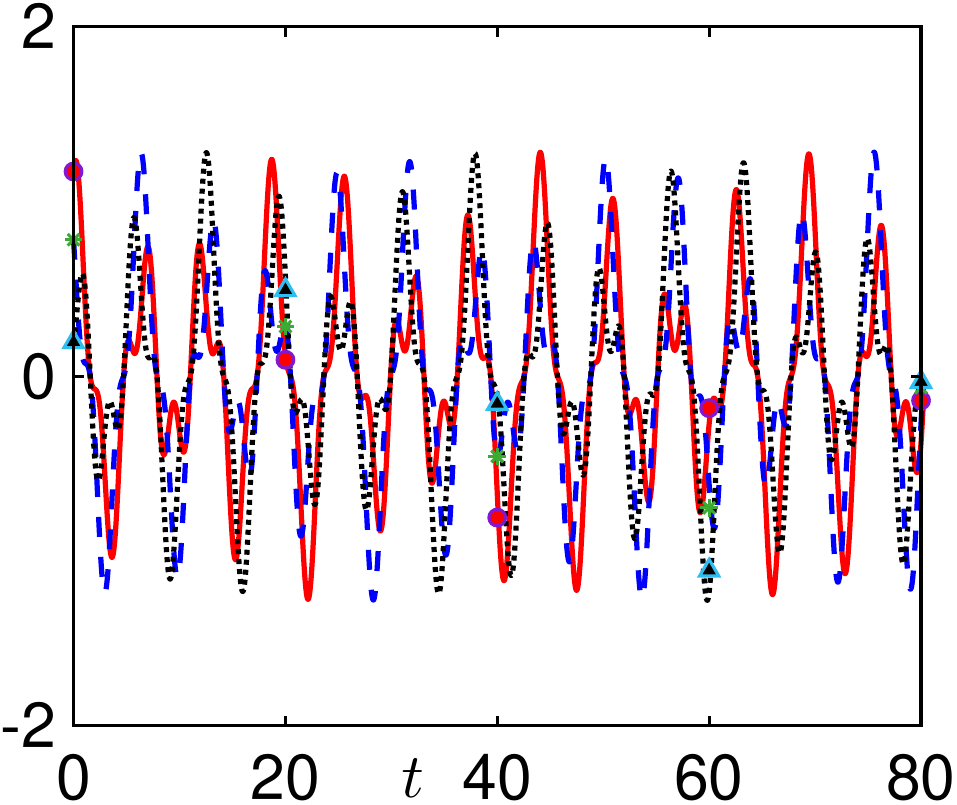}
 \hspace{-0.1cm}
  \includegraphics[width=4.2cm, height=4.2cm]{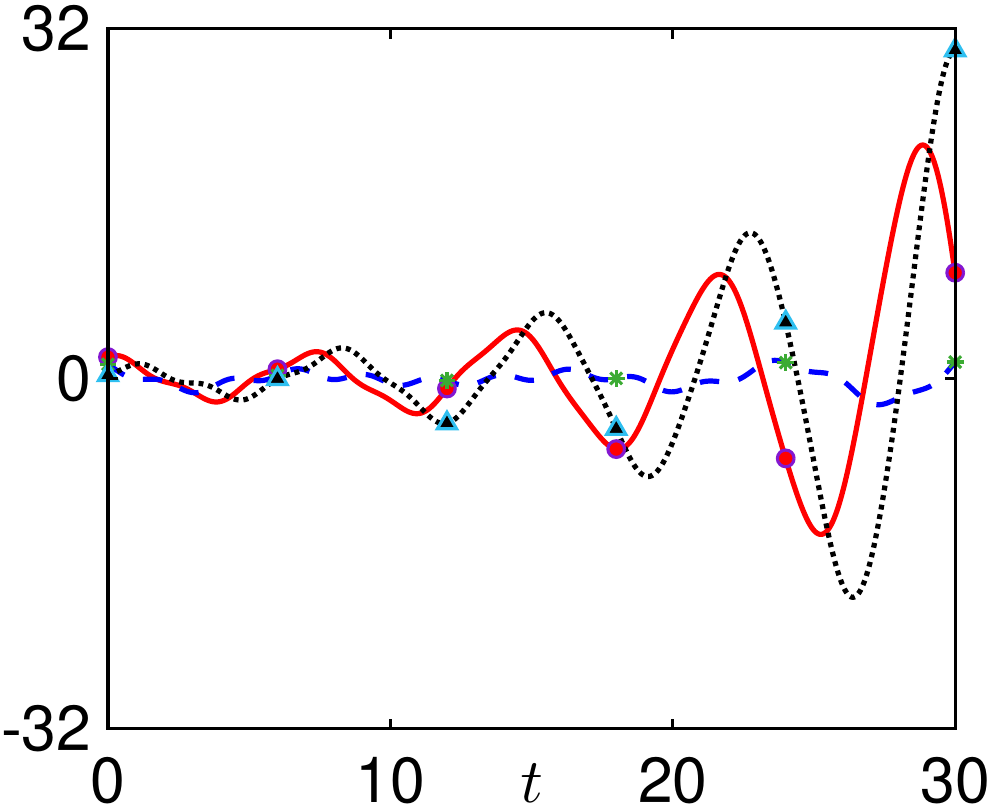}
 \hspace{-0.1cm}
 \includegraphics[width=4.2cm, height=4.2cm]{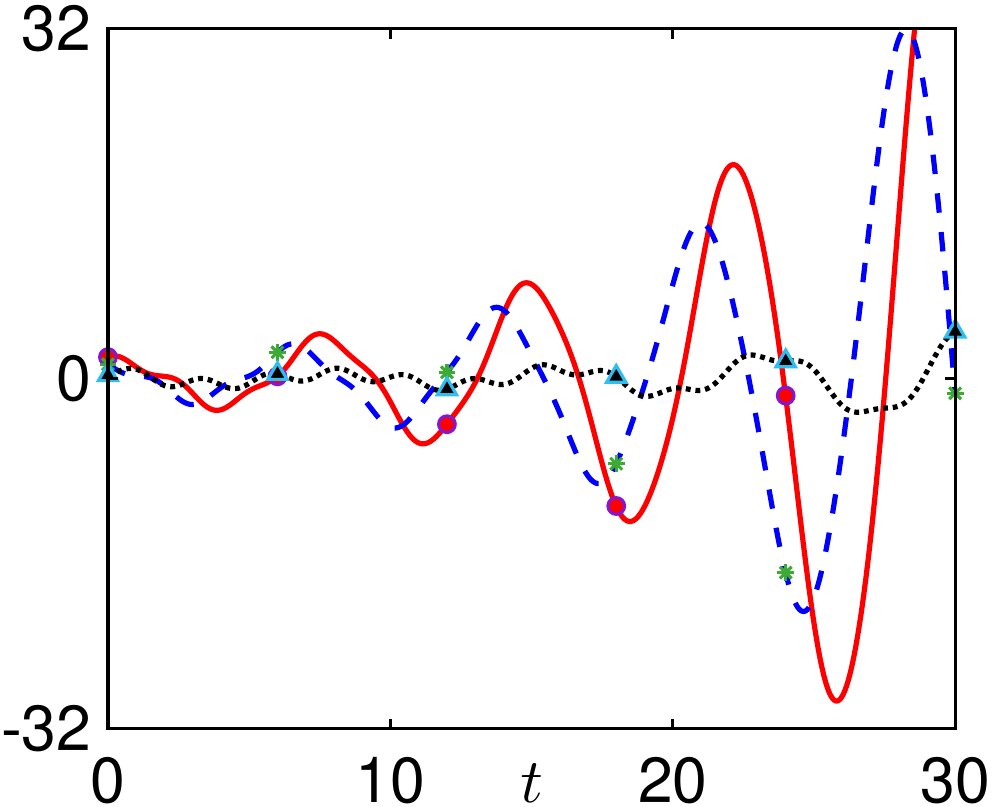}
 \\
 \caption{Motion of the center of mass $\bx_c(t) = (x_c(t), y_c(t), z_c(t))^{\top}$  for \textbf{\textrm{Case \uppercase\expandafter{\romannumeral3}}} in \textbf{Example \ref{Exmp:DynaPropCenter}}.
Top: trajectory $\bx_c(t)$ for $t \in [0, 80]$ (first parameter) and  $t \in [0, 30]$ (last two parameters).
Bottom: time evolution of $x_c(t)$ (red solid line),
$y_c(t)$ (blue dot line) and $z_c(t)$ (black  dash line) that are obtained by directly simulating the CGPEs,
where 	
circle($o$), 	
asterisk($\ast$) and triangle($\triangle$) represent the solutions to the ODEs \eqref{centerODE}-\eqref{centerODEInitial}.
The left-to-right order corresponds to the three parameters provided.}
\label{MassCenterGamma}
\end{figure}
Figures \ref{MassCenterOmega}-\ref{MassCenterGamma} depict the trajectory of the center of mass $\bx_c(t)$ and the
evolution of its coordinates, calculated by CS4 with
$h=1/8$ and $\tau = 10^{-2}$ for different $\bm{\Omega}$ and
 $\bm{\gamma}$.
The computational domain is taken as  $[-32, 32]^3$ for the latter two choices of $\bm{\gamma}$ in \textbf{Case \uppercase\expandafter{\romannumeral3}}.
We can see that
if the external trap is symmetric, i.e., $\gamma_x = \gamma_y = \gamma_z$,
the center of mass remains confined within a bounded region
and the motion is periodic when the angular velocity
$|\bm{\Omega}|$ is rational.
When the external trap is not symmetric,
the center of mass exhibits more complicated dynamics.
In addition, these simulations agree with our theoretical analysis,
and some numerical results not shown here  for brevity confirm that the dynamics of center of mass are independent of $c_0$ and $c_1$.

\subsection{Applications} \label{QuanVortex}
In this subsection, we apply our method to study the dynamics of quantized vortices in spin-1 BECs with ArbRot under different rotation vectors.
\begin{exmp}[\textbf {Dynamics of quantized vortices}]
\label{Exmp:Vortex}
Here, we set   $c_0=800$, $c_1=-100$ and $\bm{\gamma} = (1, 1, 1)^{\top}$. The initial value is taken as a stationary vortex state computed
by the preconditioned conjugate gradient method \cite{gs} for the above parameters with $\bm{\Omega} = (0, 0, 4)^{\top}$ and $k=0.1$ (the spin-orbit coupling strength), i.e., including SOC initially.
Then the dynamics of vortices are investigated for the following two cases
\begin{description}
\twoitems
{$~~\qquad\mathbf{Case ~\uppercase\expandafter{\romannumeral1}}:~~\bm{\Omega} = (0, 0, 4)^{\top}$;}
{$~~\qquad\qquad\mathbf{Case ~\uppercase\expandafter{\romannumeral2}}:~\bm{\Omega} = (1, 2, 2)^{\top}$.}
\end{description}
\end{exmp}

\begin{figure}[!htpb]
\centering
 \includegraphics[trim={4cm 0.2cm 3cm 3cm}, clip, width=4.2cm, height=4cm]{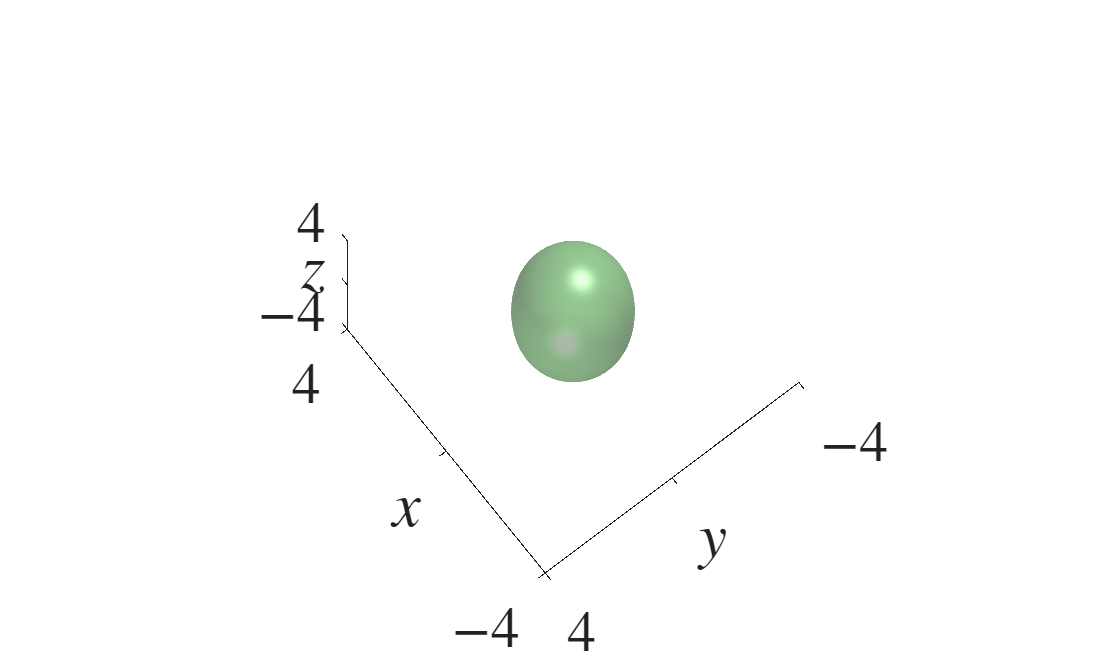}
 \hspace{-0.2cm}
 \includegraphics[trim={4cm 0.2cm 3cm 3cm}, clip, width=4.2cm, height=4cm]{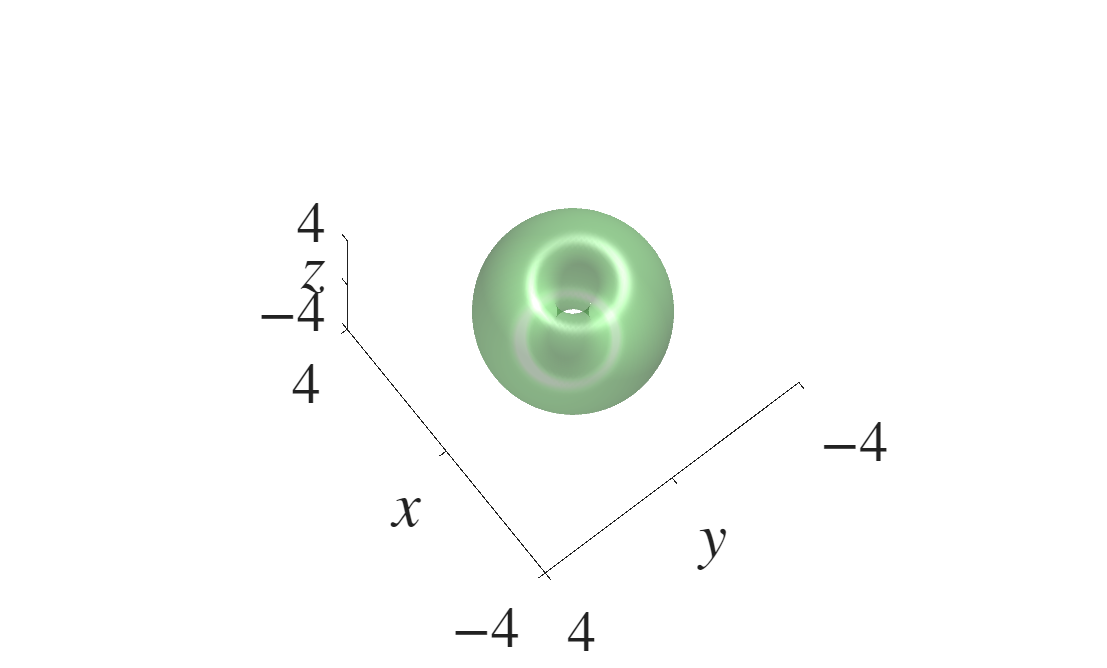}
 \hspace{-0.2cm}
 \includegraphics[trim={4cm 0.2cm 3cm 3cm}, clip, width=4.2cm, height=4cm]{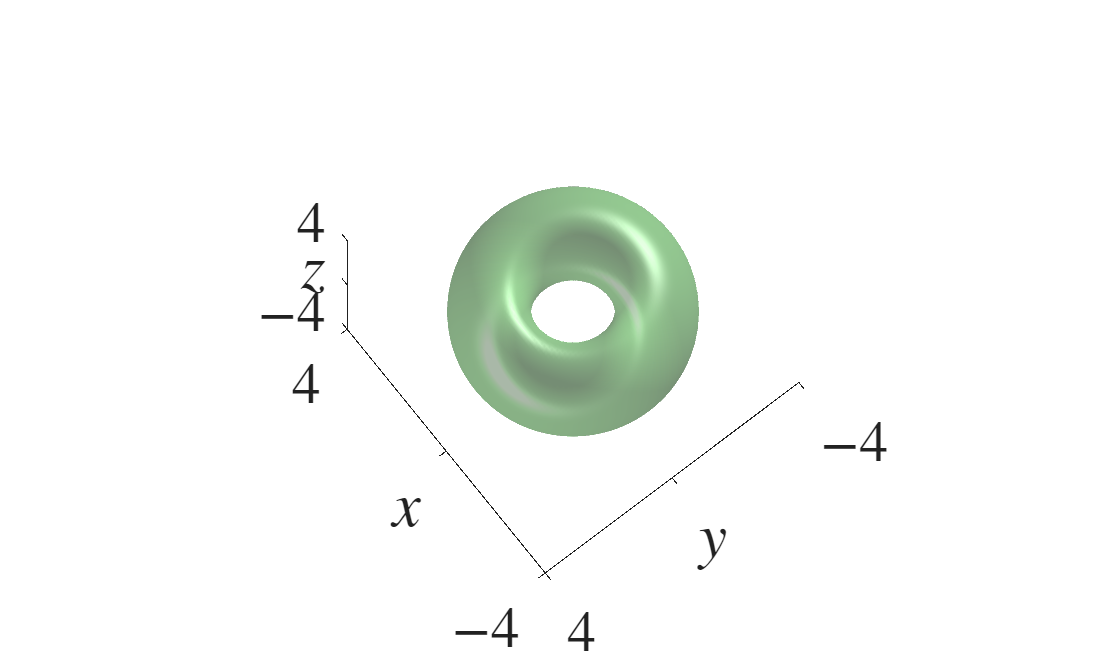}
 \\[1.5em]
 \includegraphics[trim={4cm 0.2cm 3cm 3cm}, clip, width=4.2cm, height=4cm]{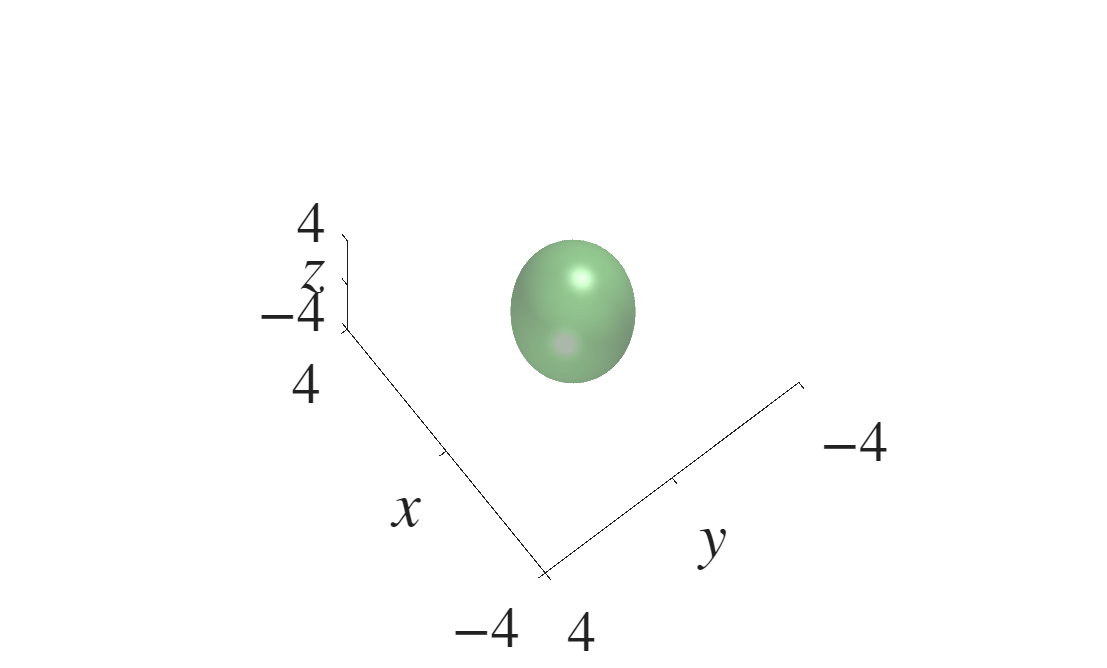}
 \hspace{-0.2cm}
 \includegraphics[trim={4cm 0.2cm 3cm 3cm}, clip, width=4.2cm, height=4cm]{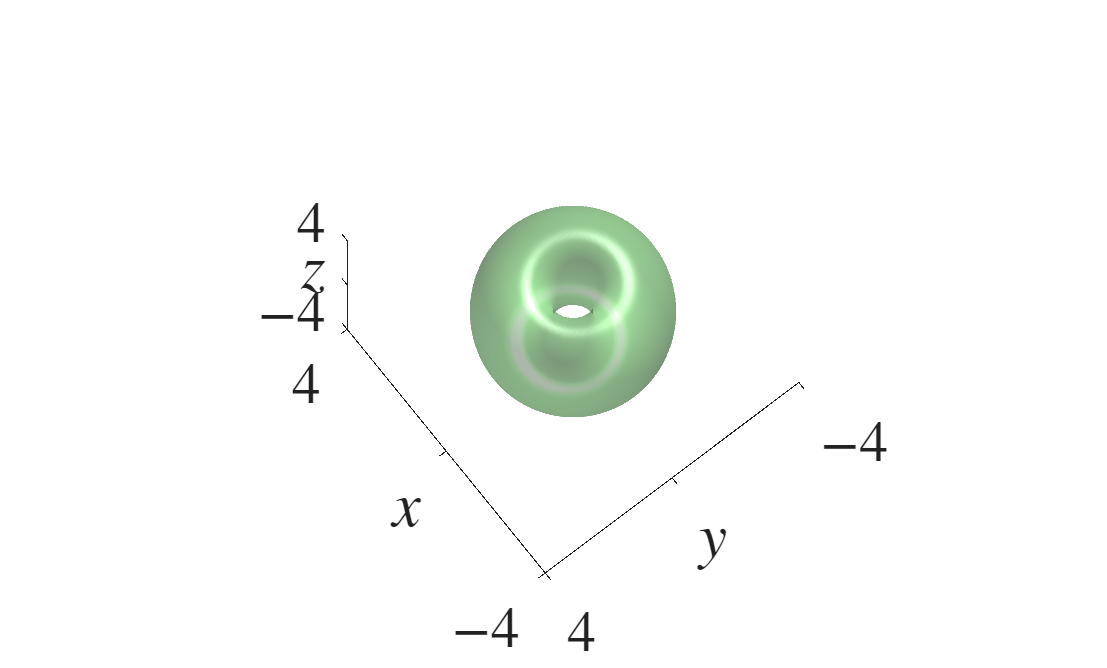}
 \hspace{-0.2cm}
 \includegraphics[trim={4cm 0.2cm 3cm 3cm}, clip, width=4.2cm, height=4cm]{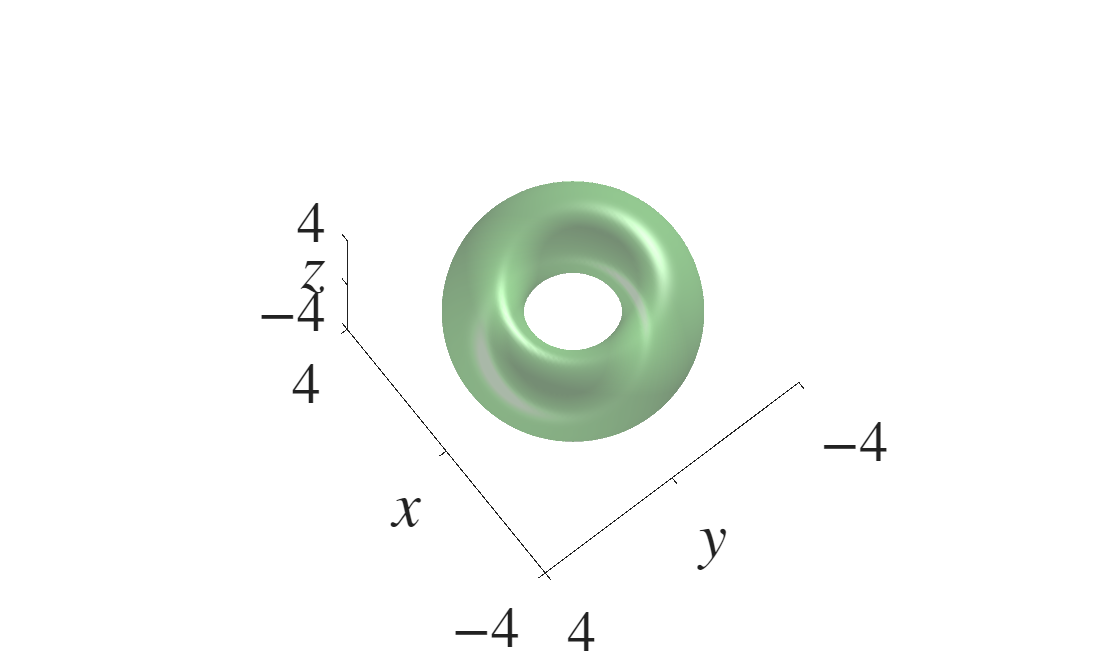}
 \\[1.5em]
 \includegraphics[trim={4cm 0.2cm 3cm 3cm}, clip, width=4.2cm, height=4cm]{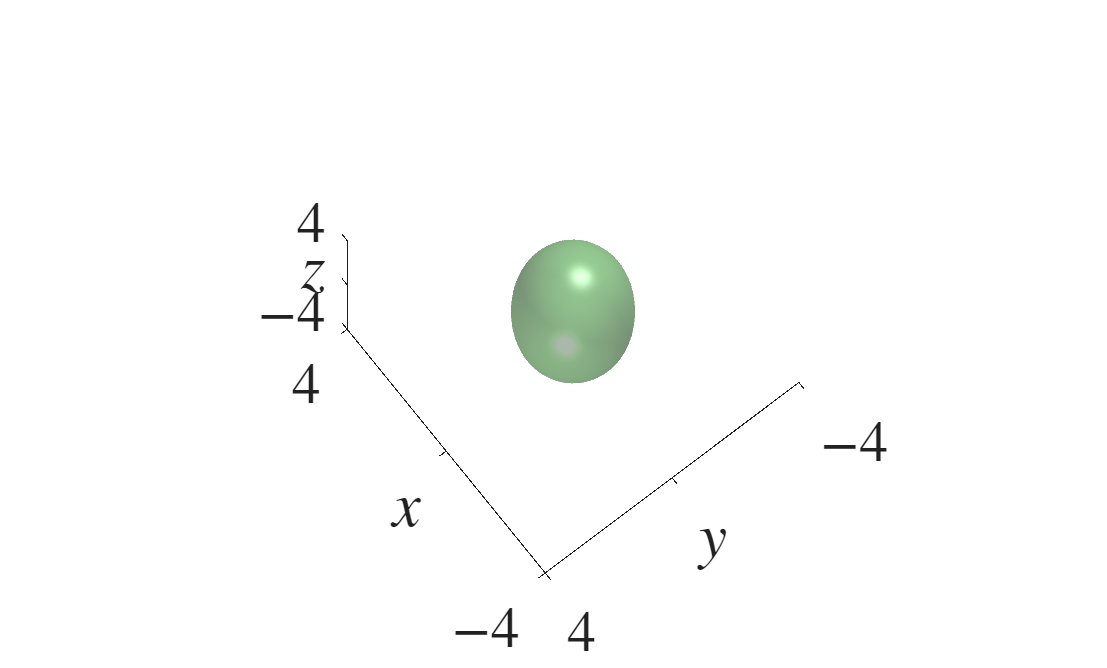}
 \hspace{-0.2cm}
 \includegraphics[trim={4cm 0.2cm 3cm 3cm}, clip, width=4.2cm, height=4cm]{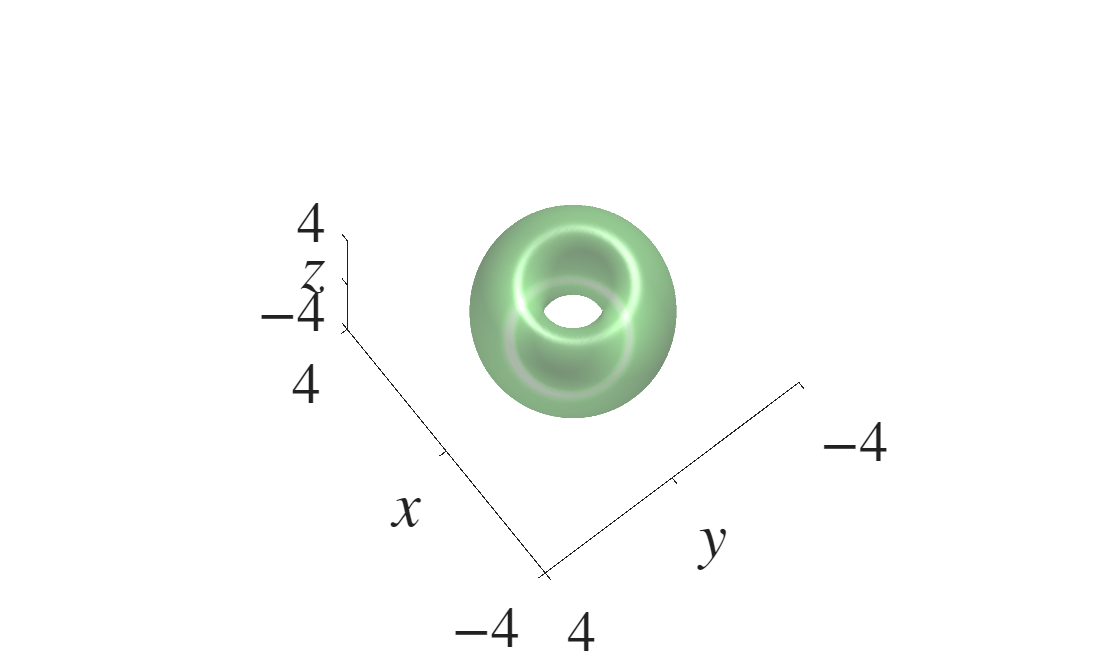}
 \hspace{-0.2cm}
 \includegraphics[trim={4cm 0.2cm 3cm 3cm}, clip, width=4.2cm, height=4cm]{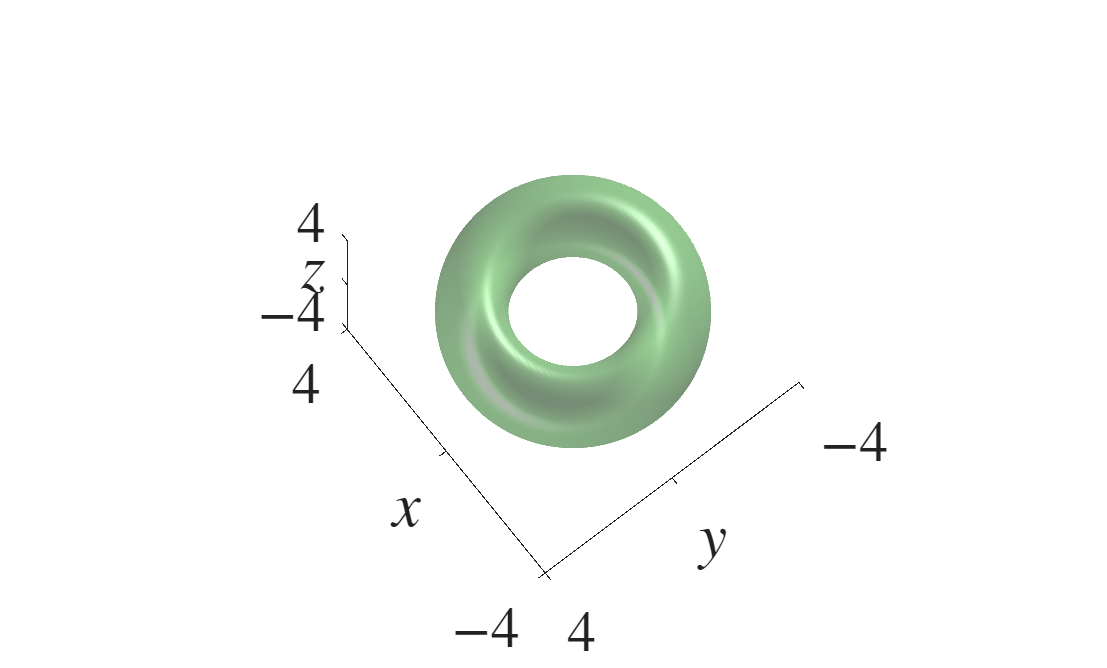}
 \\[1.5em]
 \includegraphics[trim={4cm 0.2cm 3cm 3cm}, clip, width=4.2cm, height=4cm]{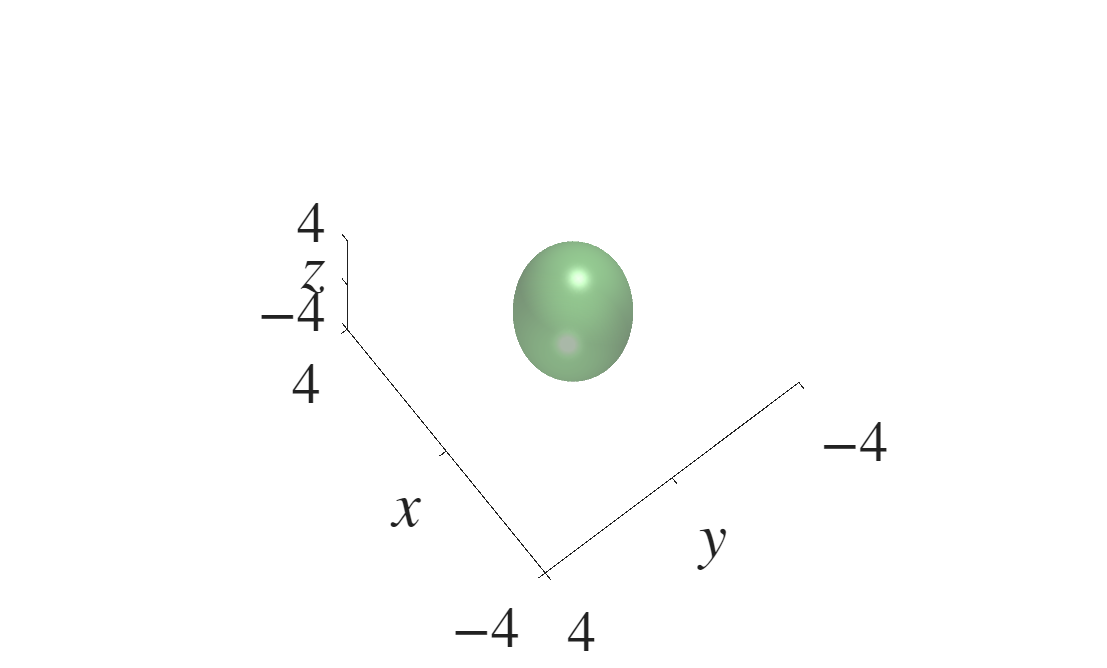}
 \hspace{-0.2cm}
 \includegraphics[trim={4cm 0.2cm 3cm 3cm}, clip, width=4.2cm, height=4cm]{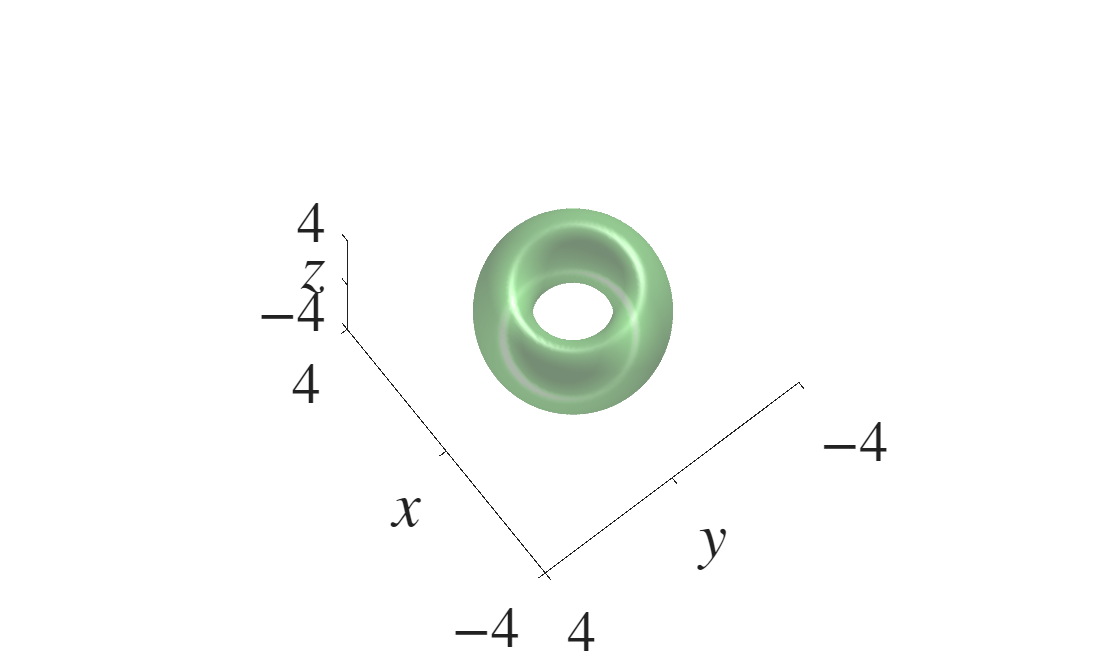}
 \hspace{-0.2cm}
 \includegraphics[trim={4cm 0.2cm 3cm 3cm}, clip, width=4.2cm, height=4cm]{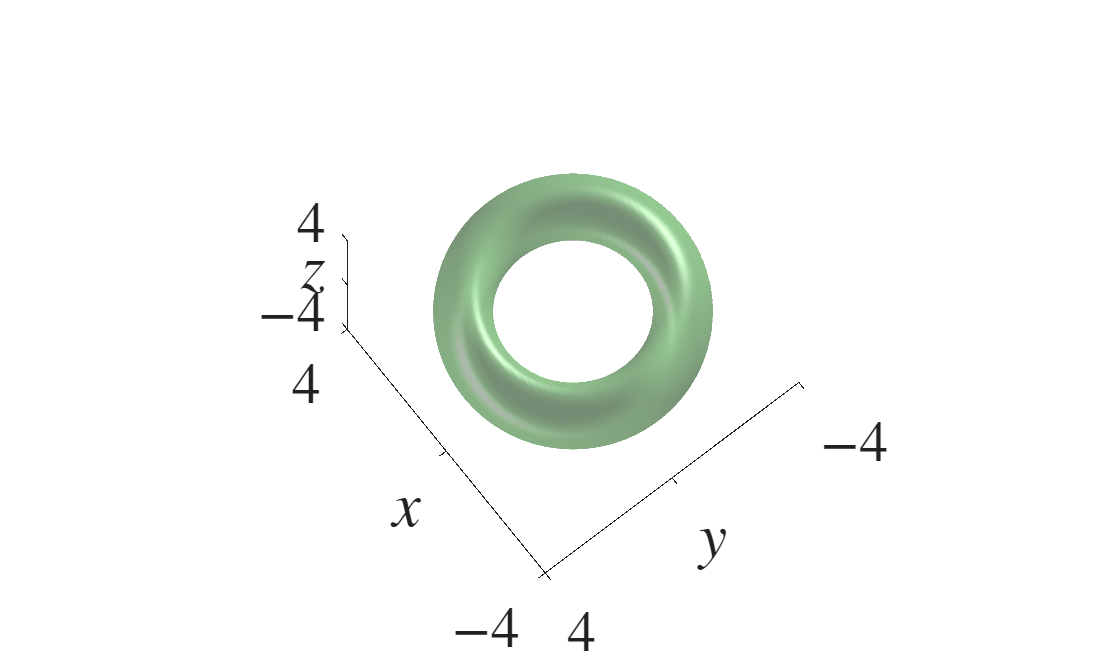}
\\
\caption{Isosurface plots of $|\psi_v|^2 = 0.002$ ($v=1, 0, -1$, left to right) at  time $t=0$, $0.2$, $0.35$ and $0.45$ (top to bottom) for \textbf{\textrm{Case \uppercase\expandafter{\romannumeral1}}} with $\bm{\Omega} = (0, 0, 4)^{\top}$ in \textbf{Example \ref{Exmp:Vortex}}.
  }
\label{VortexRot}
\end{figure}

\begin{figure}[!htpb]
\centering
 \includegraphics[trim={5cm 0.2cm 3cm 3cm}, clip, width=4.2cm, height=4.cm]{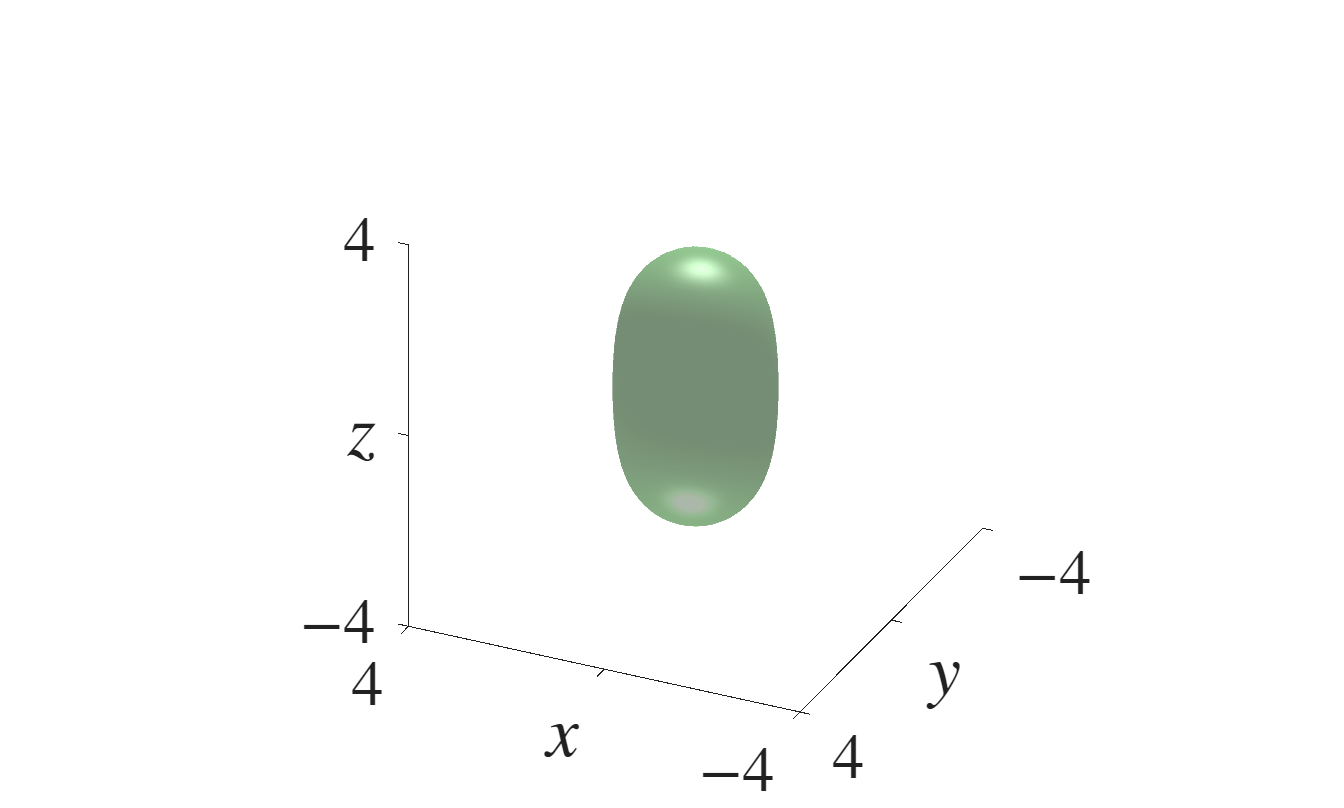}
 \hspace{-0.25cm}
 \includegraphics[trim={5cm 0.2cm 3cm 3cm}, clip, width=4.2cm, height=4.cm]{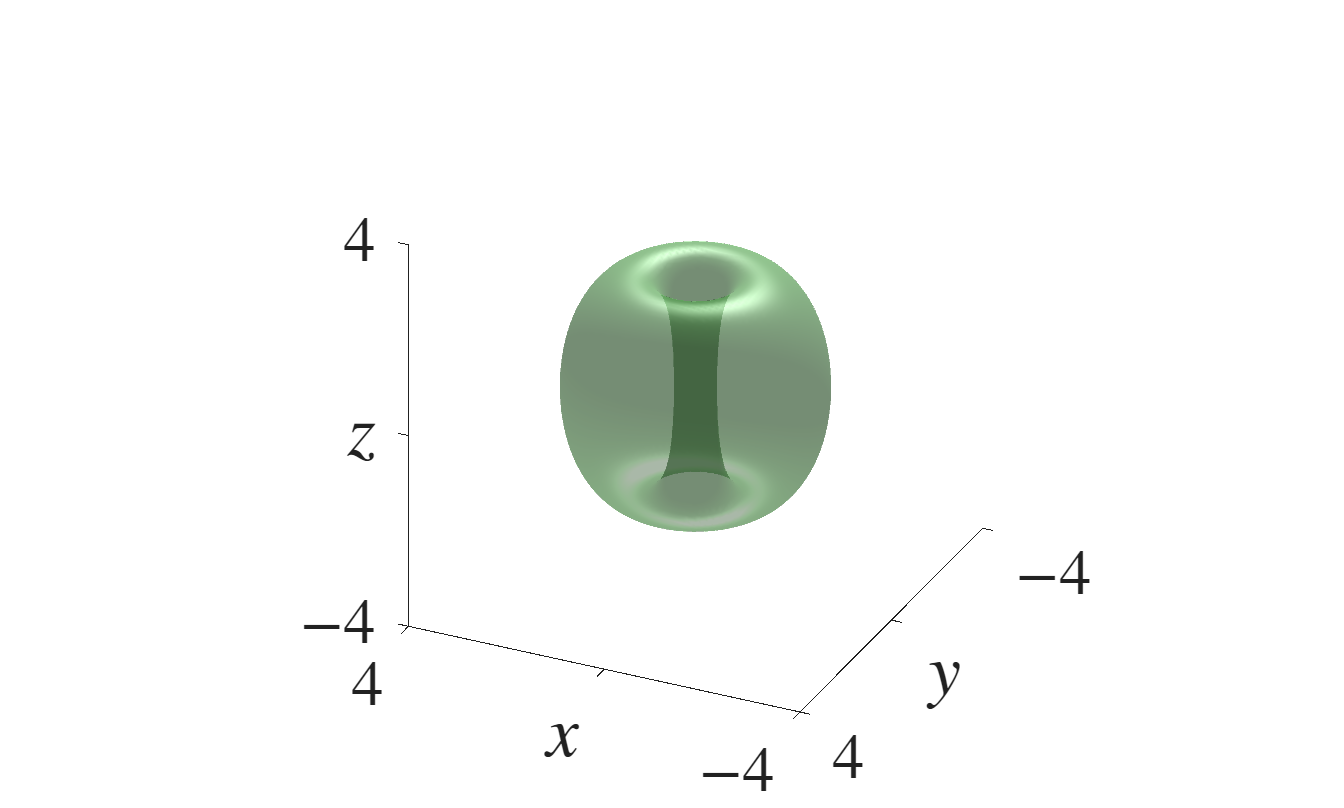}
 \hspace{-0.25cm}
 \includegraphics[trim={5cm 0.2cm 3cm 3cm}, clip, width=4.2cm, height=4.cm]{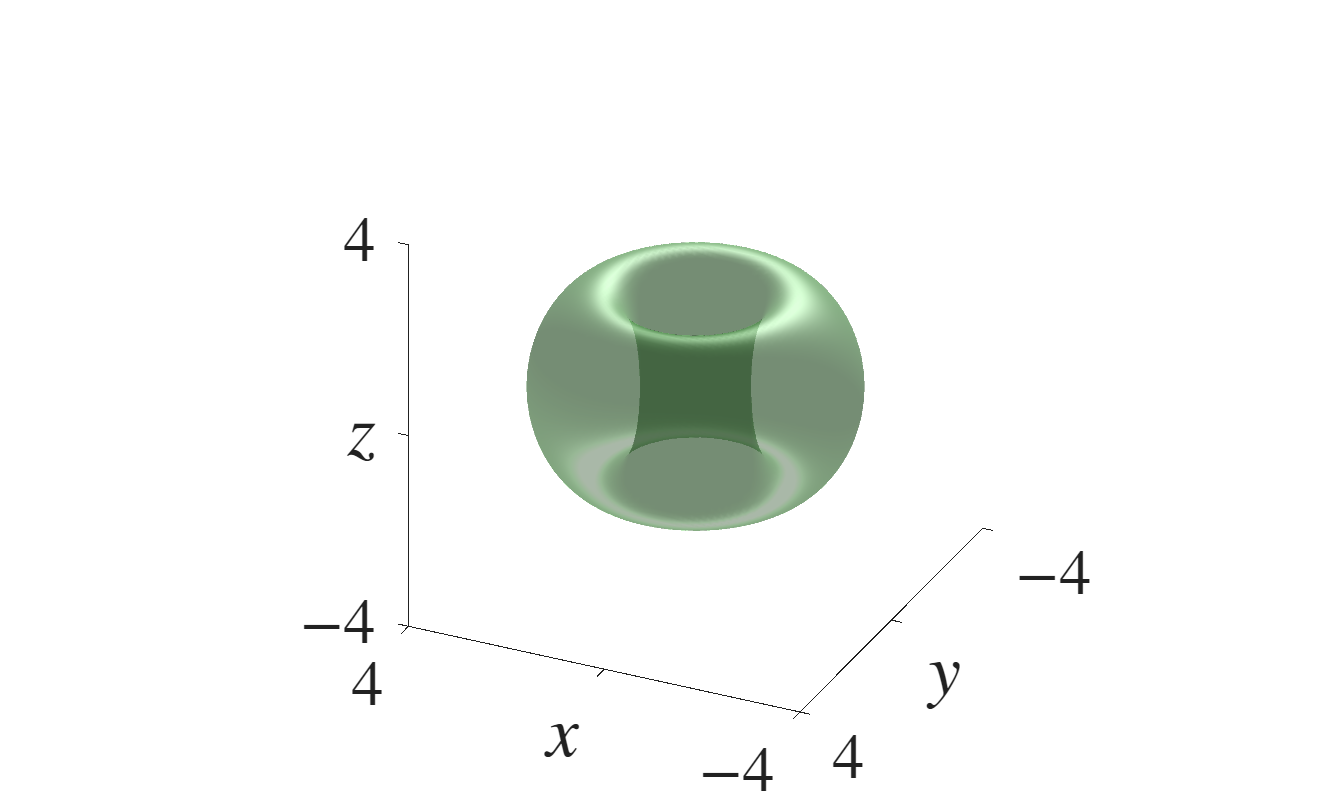}
 \\[1.5em]
 \includegraphics[trim={5cm 0.2cm 3cm 3cm}, clip, width=4.2cm, height=4.cm]{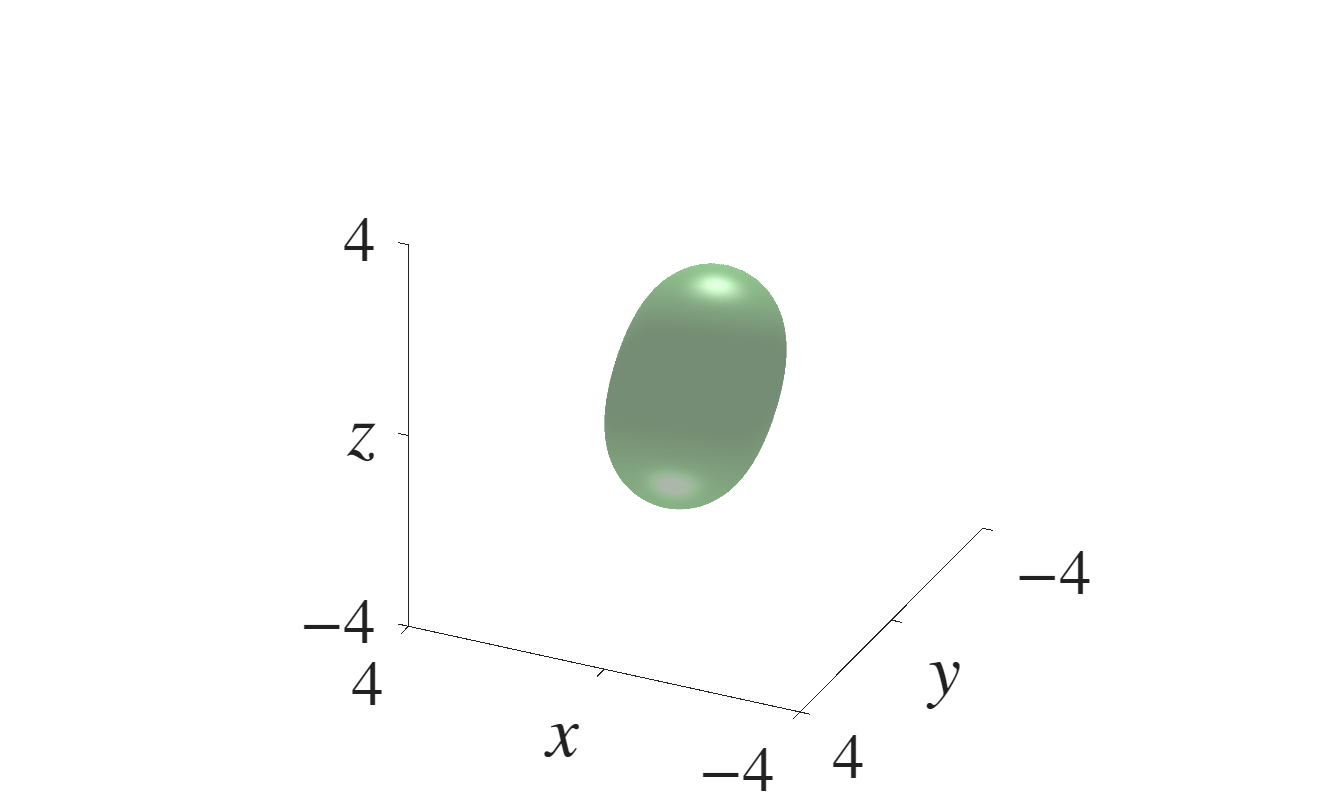}
 \hspace{-0.25cm}
 \includegraphics[trim={5cm 0.2cm 3cm 3cm}, clip, width=4.2cm, height=4.cm]{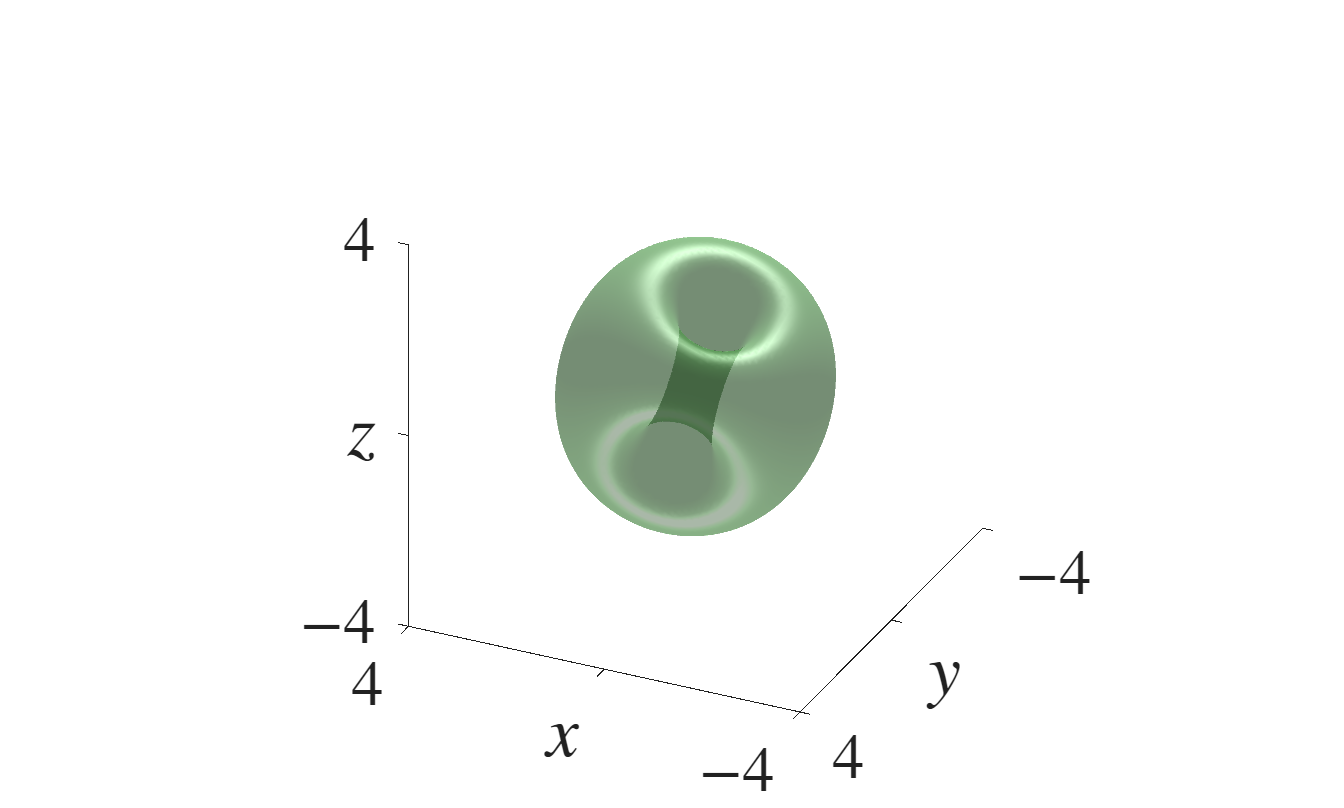}
 \hspace{-0.25cm}
 \includegraphics[trim={5cm 0.2cm 3cm 3cm}, clip, width=4.2cm, height=4.cm]{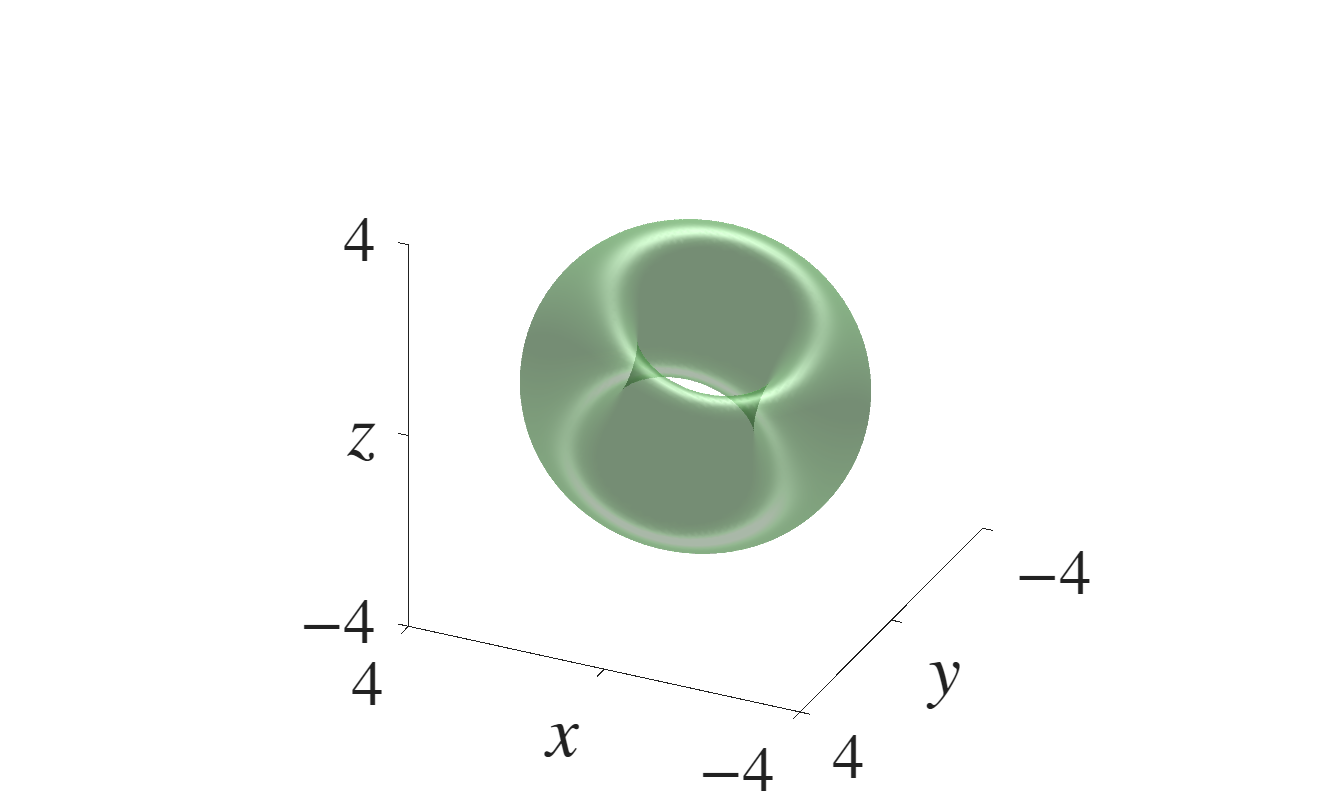}
 \\[1.5em]
 \includegraphics[trim={5cm 0.2cm 3cm 3cm}, clip, width=4.2cm, height=4.cm]{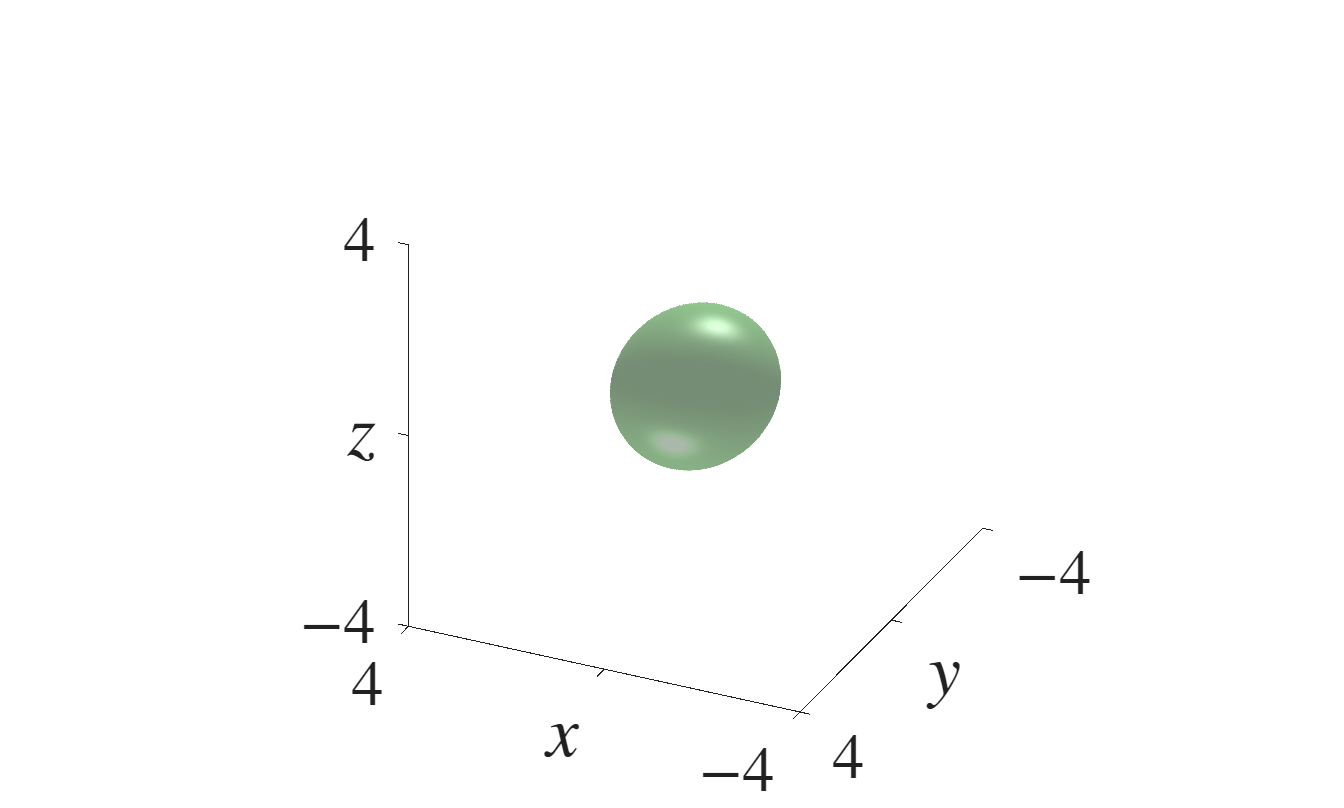}
 \hspace{-0.25cm}
 \includegraphics[trim={5cm 0.2cm 3cm 3cm}, clip, width=4.2cm, height=4.cm]{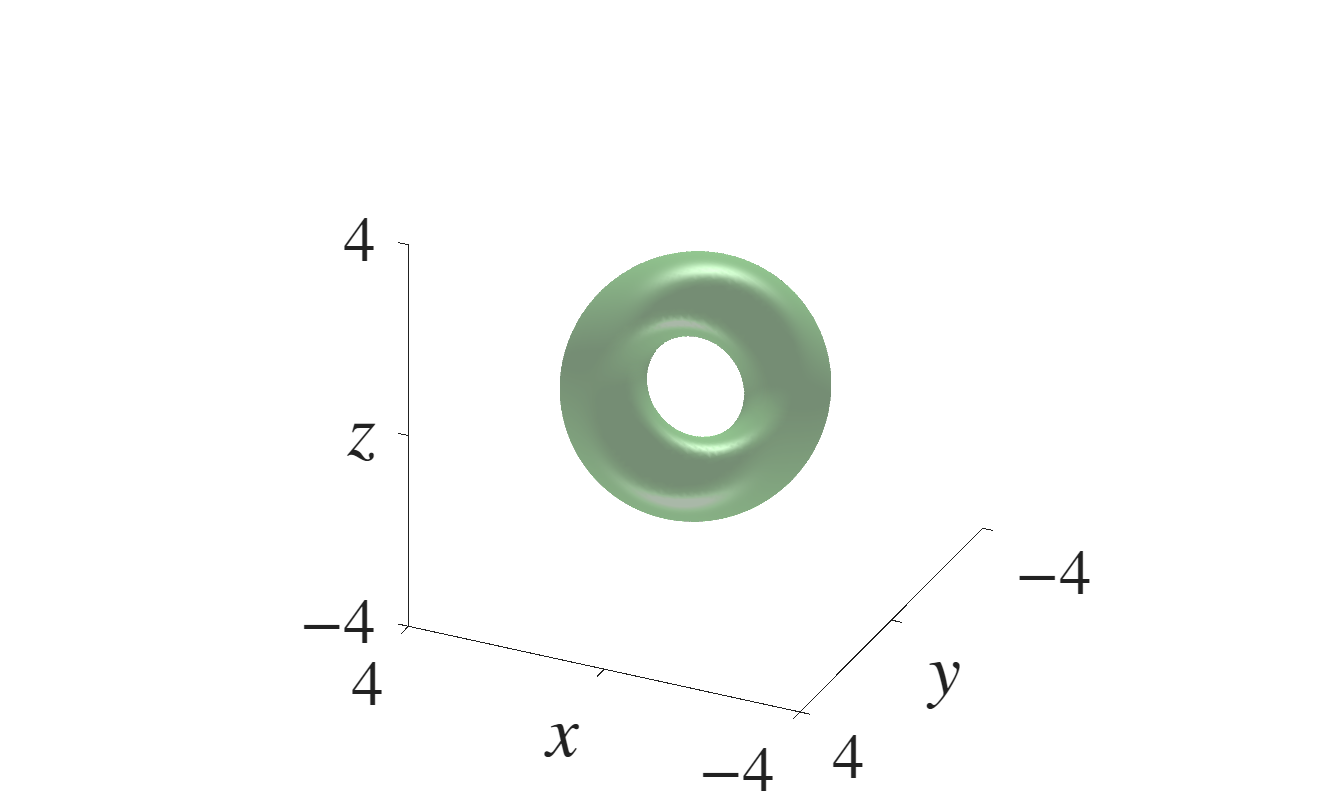}
 \hspace{-0.25cm}
 \includegraphics[trim={5cm 0.2cm 3cm 3cm}, clip, width=4.2cm, height=4.cm]{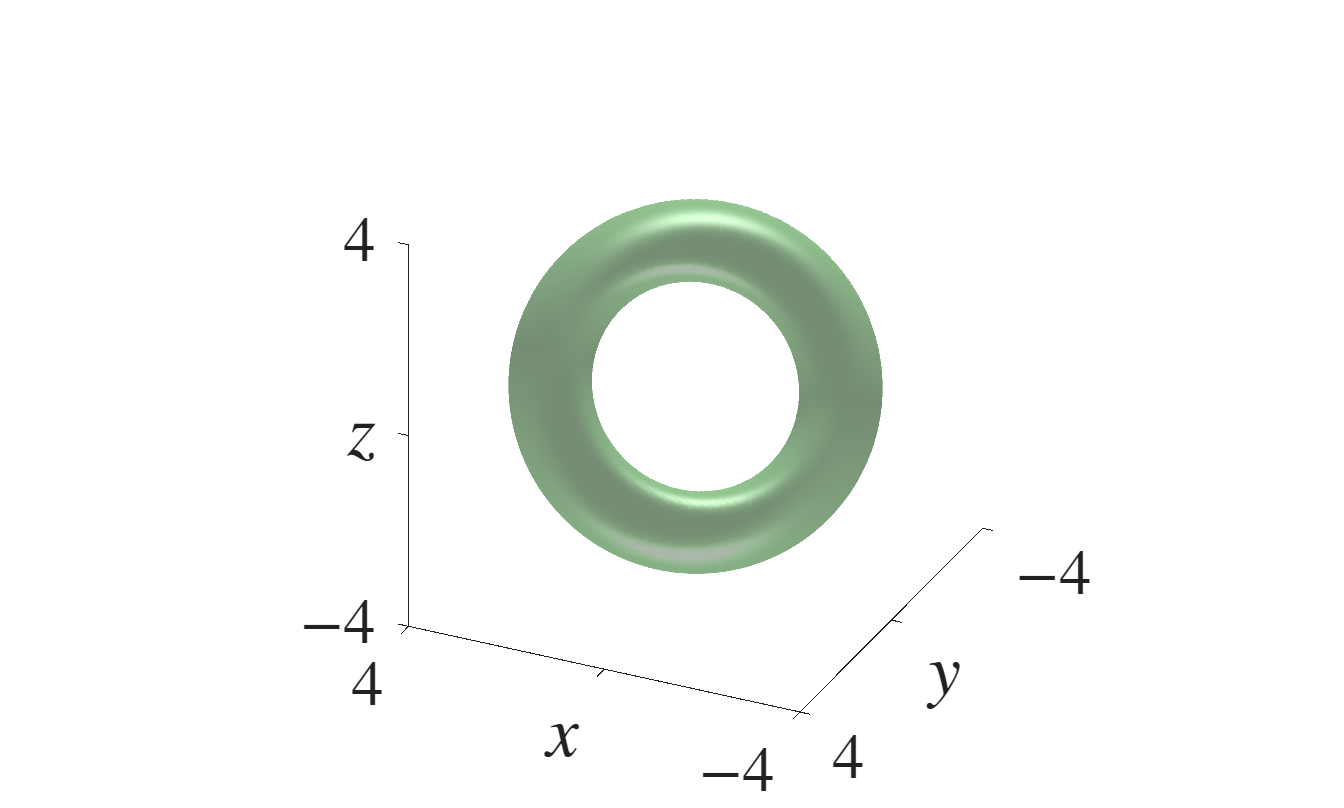}
 \\[1.5em]
 \includegraphics[trim={5cm 0.2cm 3cm 3cm}, clip, width=4.2cm, height=4.cm]{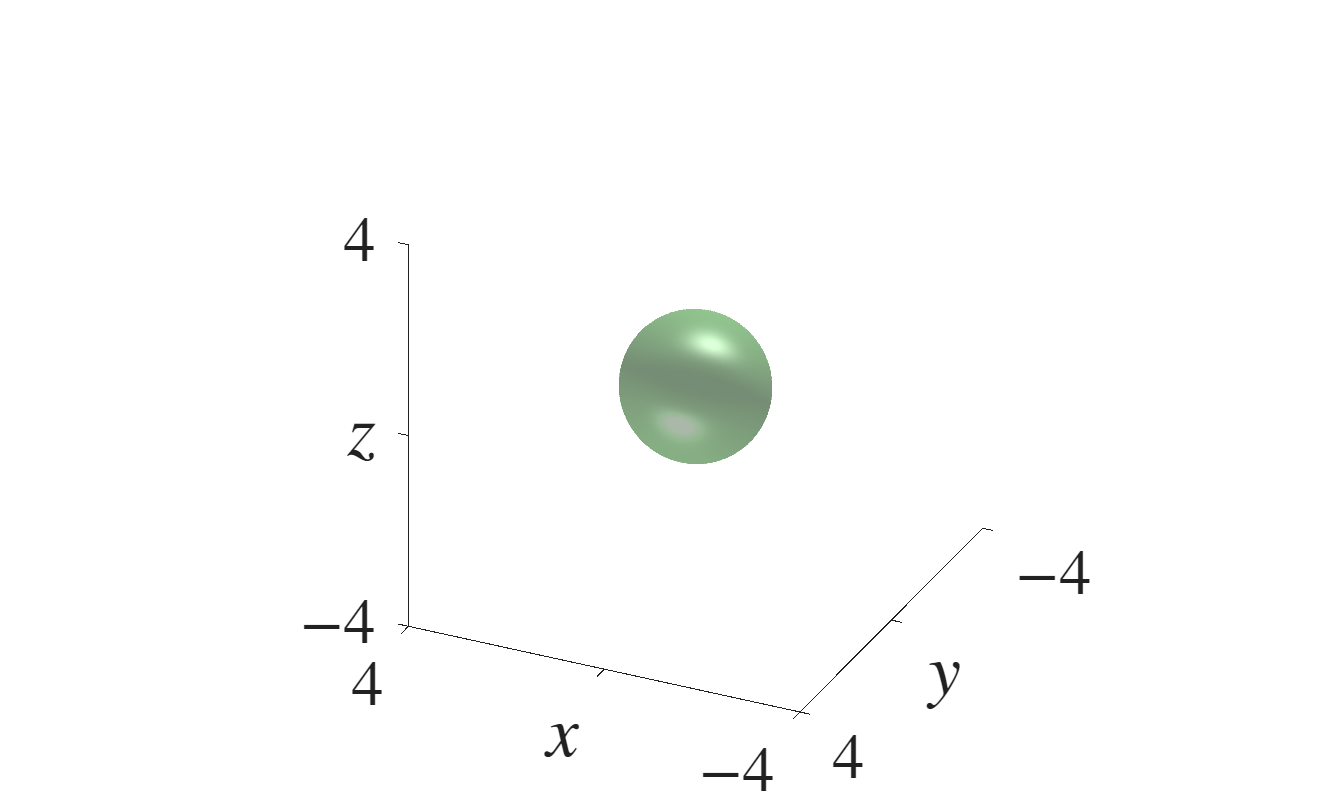}
 \hspace{-0.25cm}
 \includegraphics[trim={5cm 0.2cm 3cm 3cm}, clip, width=4.2cm, height=4.cm]{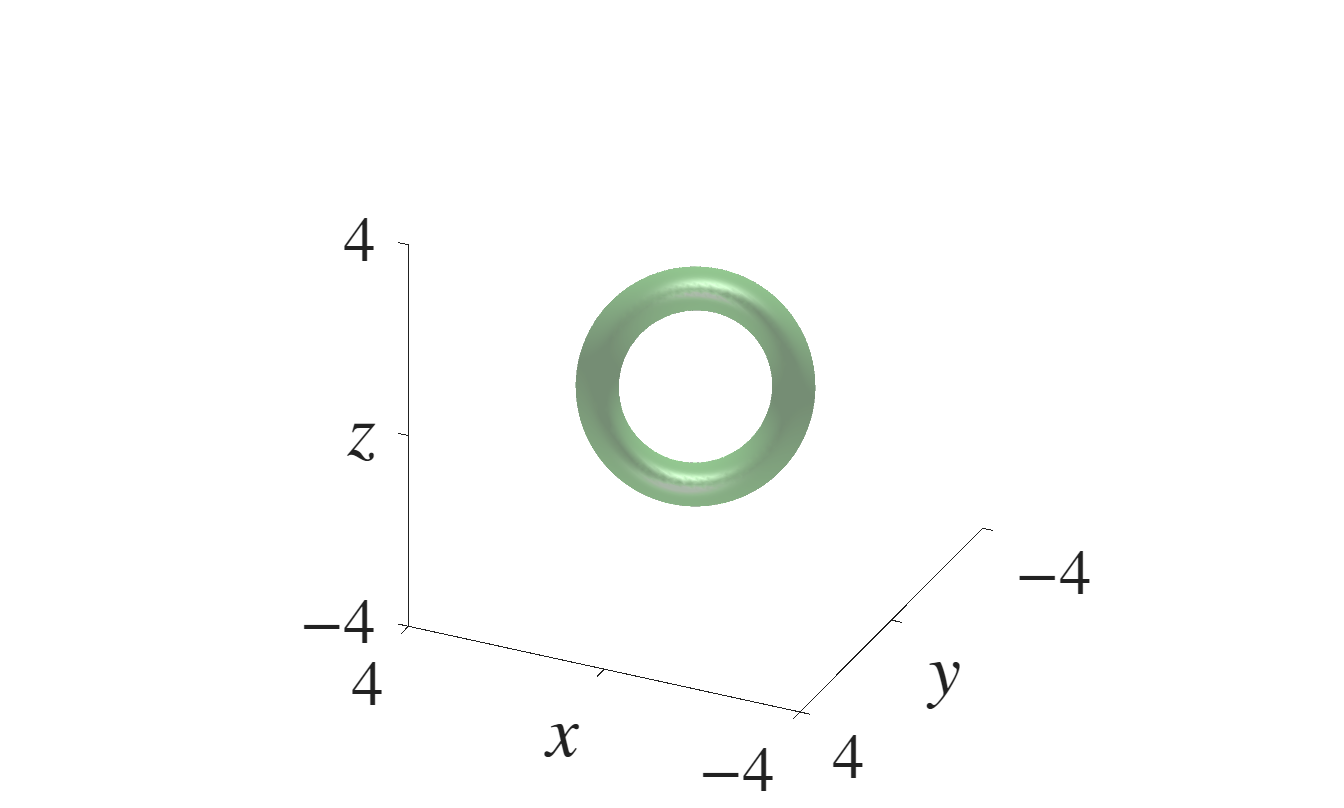}
 \hspace{-0.25cm}
 \includegraphics[trim={5cm 0.2cm 3cm 3cm}, clip, width=4.2cm, height=4.cm]{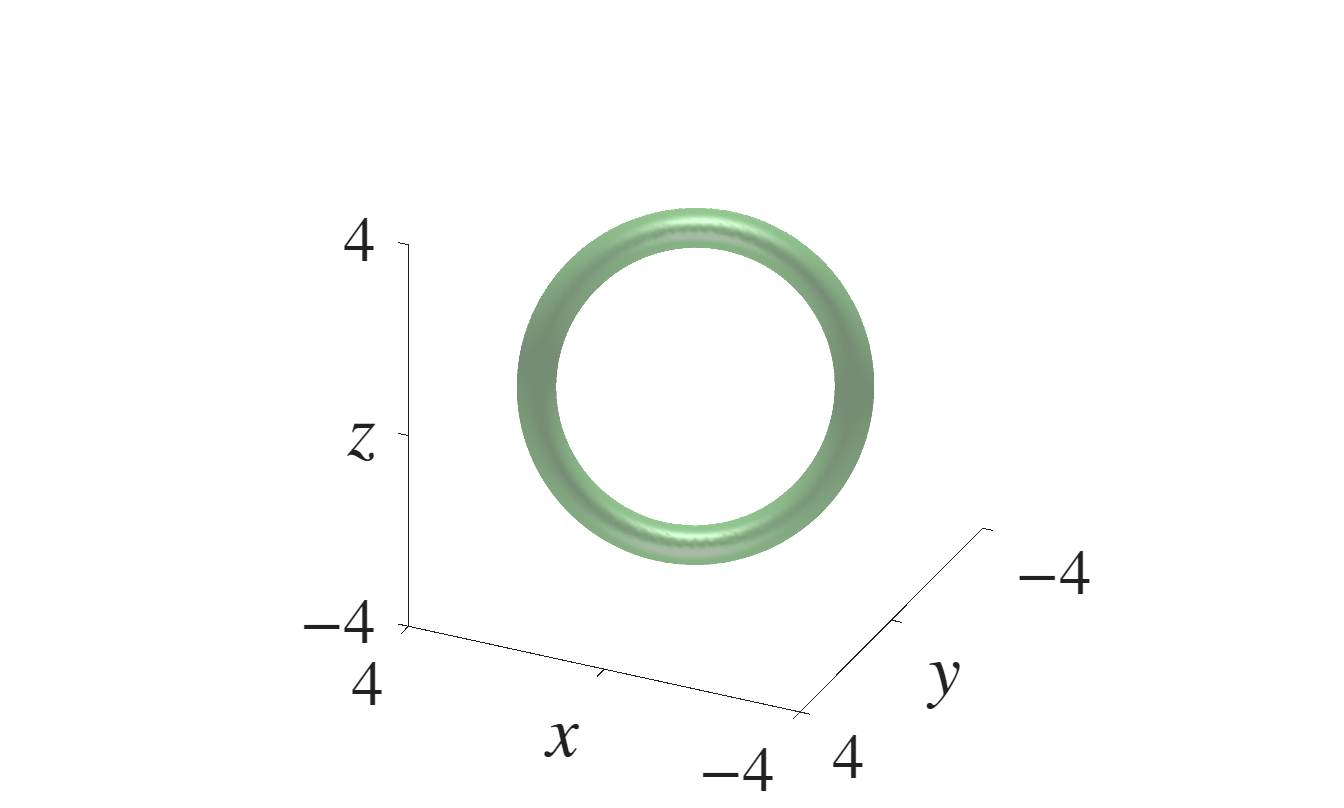}
\\
\caption{Isosurface plots of $|\psi_v|^2 = 0.002$ ($v=1, 0, -1$, left to right) at  time $t=0$, $0.2$, $0.45$ and $0.55$ (top to bottom) for \textbf{\textrm{Case \uppercase\expandafter{\romannumeral2}}} with $\bm{\Omega} = (1, 2, 2)^{\top}$ in \textbf{Example \ref{Exmp:Vortex}}.
  }
\label{VortexGenRot}
\end{figure}

In simulations, we take computational domain $[-12,12]^3$, $N=256$ and $\tau = 10^{-3}$.
Figures \ref{VortexRot}-\ref{VortexGenRot} present the isosurface plots of
the densities $|\psi_v|^2 = 0.002$ computed by CS4 at different time.
We observe that the vortices form different patterns for different rotation vectors $\bm{\Omega}$.
Since the vortex line tends to align parallel to the rotation axis, when its initial direction differs from $\bm{\Omega}$,
the vortices rotate during the dynamics until alignment is achieved.

\section{Conclusions} \label{conclusions}%
We proposed high-order splitting Fourier spectral methods for simulating the dynamics of spin-1
BECs with arbitrary-angle rotation.
The Hamiltonian is split into a linear and a nonlinear part: the linear part consists of the Laplacian term and arbitrary-angle rotation term,
while the latter contains the remaining terms and is exactly solvable in physical space.
For the linear subproblem,
we introduced a rotating coordinate transformation
to eliminate the rotation term.
The resulting subproblem is solved exactly in Fourier space.
For the fast implementation of the coordinate transformation, we developed a spectrally accurate approach based on shear decomposition.
This method achieves the computational efficiency of one-dimensional FFT/iFFT operations
and remains valid for arbitrary time step.
Thus the operator splitting method was easily constructed, which is efficient, unconditionally stable, spectrally accurate in space and of second/fourth-order in time.
Ample numerical results were provided to validate the accuracy and efficiency, and to illustrate the application in confirming dynamical properties and capturing several interesting physical phenomena.
Our methods can be extended to more general systems, such as rotating  dipolar spin-F BECs and spin-orbit-coupled BECs, and we shall
report them in future work.

\appendix
\section{Ensuring non-negative rotation components}\label{substitution}
Here we present the specific variable substitutions for problem \eqref{LinearSubprob} (see Table \ref{tab_substitution}), which reformulate it into an equivalent problem with a transformed rotation vector $\tilde{\bm{\Omega}}$ whose components are all non-negative.
\begin{table}[h!]
   \centering
   \caption{The variable substitutions for various $\bm{\Omega}$ ($\omega_i\geq0, i=1,2,3$).}\label{tab_substitution}
   \begin{tabular}{lll}
      \toprule
      $\qquad\bm{\Omega}$                           & variable substitution                    & $\qquad\tilde{\bm{\Omega}}$ \\
      \midrule
      $(\omega_1, \omega_2, -\omega_3)^{\top}$  & $\tilde{x}=-y,~ \tilde{y}=-x,\!~ \tilde{z}=z$
      & $(\omega_2, \omega_1, \omega_3)^{\top}$  \\
      $(\omega_1, -\omega_2, \omega_3)^{\top}$  & $\tilde{x}=-z,~ \tilde{y}=y,~~~ \tilde{z}=-x$
      & $(\omega_3, \omega_2, \omega_1)^{\top}$  \\
      $(-\omega_1, \omega_2, \omega_3)^{\top}$  & $\tilde{x}=x,~~~ \tilde{y}=-z,~ \tilde{z}=-y$
      & $(\omega_1, \omega_3, \omega_2)^{\top}$  \\
      \midrule
      $(\omega_1, -\omega_2, -\omega_3)^{\top}$  & $\tilde{x}=-x,~ \tilde{y}=y,~~~ \tilde{z}=z$
      &
      \multirow{3}{*}{$(\omega_1, \omega_2, \omega_3)^{\top}$}  \\
      $(-\omega_1, \omega_2, -\omega_3)^{\top}$  & $\tilde{x}=x,~~~ \tilde{y}=-y,~ \tilde{z}=z$
      &\\
      $(-\omega_1, -\omega_2, \omega_3)^{\top}$  & $\tilde{x}=x,~~~ \tilde{y}=y,~~~ \tilde{z}=-z$
      & \\
      \midrule
      $(-\omega_1, -\omega_2, -\omega_3)^{\top}$  & $\tilde{x}=y,~~~ \tilde{y}=x,~~~ \tilde{z}=z$
      & $(\omega_2, \omega_1, \omega_3)^{\top}$  \\
      \bottomrule
   \end{tabular}
\end{table}


\end{document}